\documentclass[noinfoline]{imsart}

\usepackage[a4paper,top=3cm, bottom=3cm, left=2.5cm, right=2.5cm]{geometry}
\usepackage{fancyhdr}

\RequirePackage{amsthm,amsmath,amsfonts,amssymb, mathtools}

\allowdisplaybreaks

\usepackage{mathrsfs}

\usepackage{graphicx}
\graphicspath{{Images/}}
\usepackage{float}
\usepackage[figuresright]{rotating} 
\usepackage{subcaption}
\usepackage{pdfpages}
\usepackage{wrapfig}
\usepackage{tikz}
\usepackage{pgfplots}
\usepgfplotslibrary{groupplots}
\pgfplotsset{compat=1.18}

\RequirePackage[colorlinks,citecolor=blue,urlcolor=blue]{hyperref}
\usepackage{array}
\usepackage{booktabs}
\usepackage{enumerate}
\usepackage{enumitem} 
\usepackage{multicol}
\usepackage{algorithm}
\usepackage[noend]{algpseudocode}

\renewcommand{\b}[1]{\mathbb{#1}}
\newcommand{\m}[1]{\mathbf{#1}}
\renewcommand{\c}[1]{\mathcal{#1}}

\newcommand{\defeq}{\vcentcolon=}

\DeclareMathOperator*{\argmin}{\arg\min}

\DeclareMathOperator{\tr}{tr}

\DeclareMathOperator{\diag}{Diag}
\DeclareMathOperator*{\spann}{span}
\DeclareMathOperator{\vect}{vec}

\newcommand{\vertiii}[1]{{\left\vert\kern-0.25ex\left\vert\kern-0.25ex\left\vert #1
    \right\vert\kern-0.25ex\right\vert\kern-0.25ex\right\vert}} 

\newcommand{\vertii}[1]{{\vert\kern-0.25ex\vert\kern-0.25ex\vert #1
    \vert\kern-0.25ex\vert\kern-0.25ex\vert}} 

\usepackage[capitalize]{cleveref}

\theoremstyle{plain}
\newtheorem{theorem}{Theorem}[section]

\newtheorem{corollary}[theorem]{Corollary}

\newtheorem{proposition}[theorem]{Proposition}

\newtheorem{definition}{Definition}[section]

\newtheorem{remark}{Remark}[section]

\begin{document}

\pagestyle{fancy}
\fancyhead{}
\renewcommand{\headrulewidth}{0pt}
\fancyhead[CE]{Rajita Chandak and Ho Yun}
\fancyhead[CO]{Subspace Alignment and Dimension Selection}

\begin{frontmatter}
  \title{{\large Cross-Subspace Alignment and Dimension Selection: \\ Asymptotic Geometry and Phase Transitions}}

  \begin{aug}
    \author[A]{\fnms{Rajita} \snm{Chandak}\ead[label=e2]{rchandak3@wisc.edu}}
    \and
    \author[B]{\fnms{Ho} \snm{Yun}\ead[label=e1]{hy3005@columbia.edu}}

    \address[A]{Department of Statistics, University of Wisconsin-Madison, \printead{e2}}
    \address[B]{Department of Statistics, Columbia University, \printead{e1}}
\end{aug}

\begin{abstract}
Principal Component Analysis (PCA) is a popular model for dimension reduction. While the asymptotic alignment between empirical and population subspaces is classically well-understood, the geometric relationship between independent empirical subspaces—a scenario central to cross-validation and distributed learning—remains largely unformalized. In this work, we derive the exact deterministic macroscopic limit of the chordal distance between independent empirical subspaces under data splitting. We establish a novel split-sample phase transition. We prove that sample splitting artificially inflates the local dimensionality ratio, thereby demanding a strictly stronger spectral gap for cross-subspace recovery compared to the full-sample ideal, and prove that a symmetric half-split is optimal, yielding an explicit closed-form boundary for recovery.
\end{abstract}


\end{frontmatter}

\maketitle

\section{Introduction}
\begin{figure}[h!]
  \centering
  \begin{subfigure}[b]{0.45\textwidth}
    \centering
    \includegraphics[width=\textwidth]{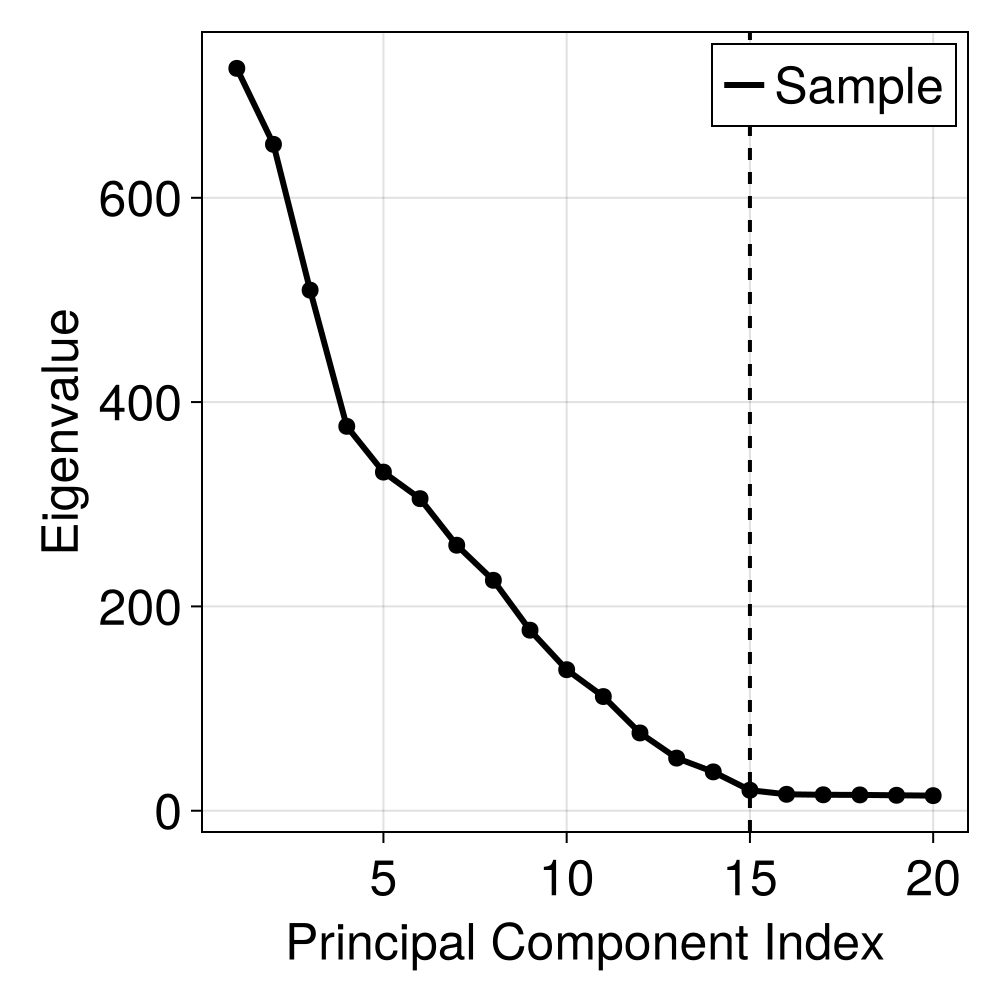}
    \caption{Scree plot.\\ \ }
    \label{fig:bad_scree}
  \end{subfigure}
  \hfill
  \begin{subfigure}[b]{0.45\textwidth}
    \centering
    \includegraphics[width=\textwidth]{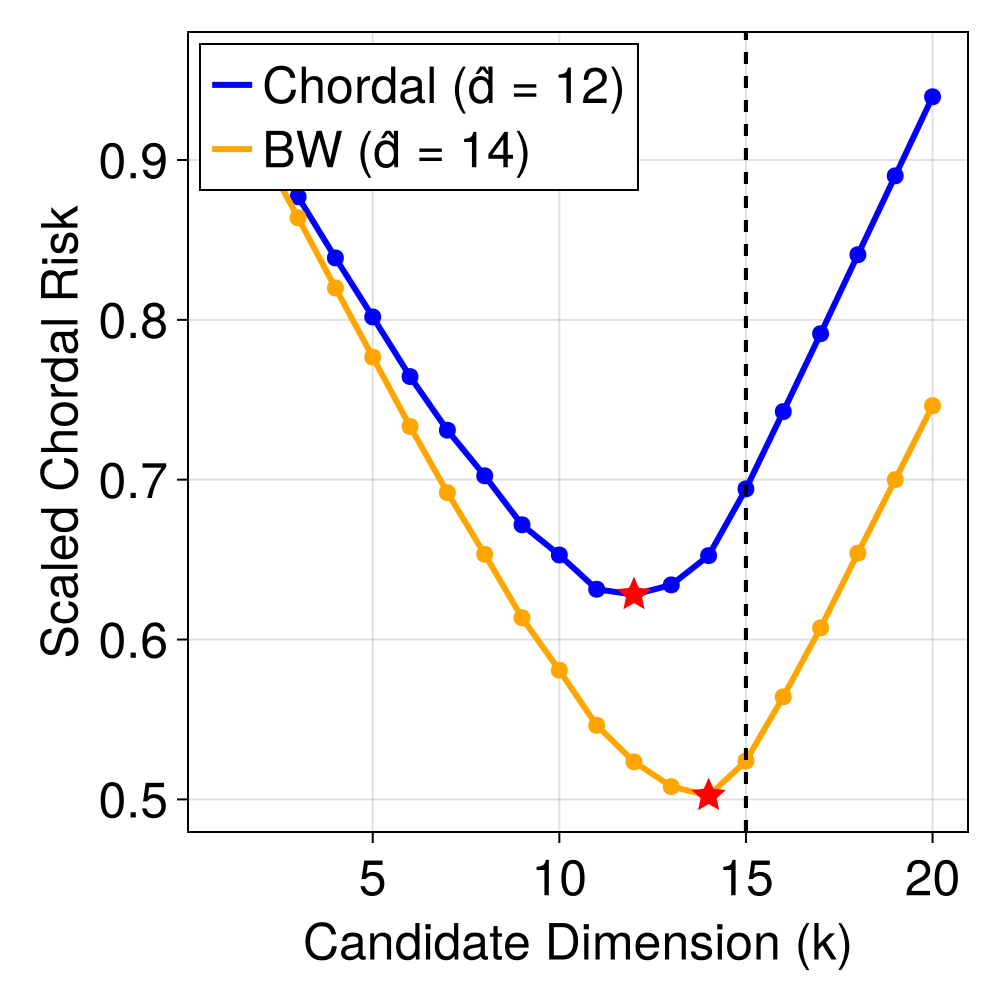}
    \caption{Two proposed risk metrics of the present work: Chordal and
      Bures-Wasserstein (BW).}
    \label{fig:risk-curve}
  \end{subfigure}
    \caption{The scree plot is not always an informative measure for identifying
    the low-dimensional signal. Sample risk curves and Scree plot are generated
    from a sample size of $n=100$ with $p=1000$ features and an SNR of $15$. The true dimensionality
    is $d^*=15$ while the reference dimension allowed for the risk curves is
    $D=20$. The Chordal risk estimates $\hat{d} = 12$, the Bures-Wasserstein
    estimates $\hat{d}=14$ while the Scree plot has no clear elbow for selecting
  the dimension.}
    \label{fig:badex}
\end{figure}

Principal Component Analysis (PCA) \cite{hotelling1933analysis, pearson1901lines}
is often cited as one of the most popular tools for dimension reduction not just
within statistics, but also across a wide range of fields from
genetics \cite{jombart2010discriminant} and climate modeling
\cite{jewson2020alternative} to financial analysis \cite{pasini2017principal}
and computer graphics \cite{bouwmans2014robust}.
Much of the existing theory provides guarantees on the asymptotic behaviour of
the empirical components derived from PCA \cite[and references
therein]{zou2006sparse, warmuth2008randomized, guan2009sparse, journee2010generalized, jolliffe2016principal, simonov2019refined}. However, with
high-dimensional datasets becoming abundant, it is unclear whether the
classical approach of using expert knowledge or some prior understanding of the
low-dimensional representation of the data can be relied upon. Furthermore, with
the goal of \textit{allowing the data to speak for itself}, in recent times,
there has been a shift to removing prior expert knowledge that can induce
unquantifiable bias. As a result, accommodating for incorrect priors,
it is of importance to provide guarantees on the consistent extraction of
low-dimensional projections. A particularly relevant example is of the Scree
plot~\cite{cattell1966scree} being used to provide an initial guess on the true
underlying sub-space of dimension $d^{*}$, a standard approach in the
literature. As illustrated in Figure~\ref{fig:badex}, however, the Scree plot
can become uninformative with high-dimensional data.
In Figure~\ref{fig:bad_scree}, there is no clear ``elbow'' and so researchers may
disagree on the right cut-off dimension. Particularly, given the nature of
such guess-work involved in pin-pointing the ``elbow'' of a scree plot, one can often
end up with a incorrect specification of the low-dimensional subspace that the
data is projected onto, leading to a loss in predictive power of any downstream
estimation or inference task. This paper aims to provide a solution to such
high-dimensional PCA applications (see Figure~\ref{fig:risk-curve}), providing two
alternative risk metrics that have robust guarantees in detecting informative
low-dimensional projection spaces.
Beyond the problem of subspace selection, one can also treat this
high-dimensional PCA problem
within the adaptive estimation or the classical machine learning
framework of using test-train splits of the data to validate the estimated
signal space. For instance, consider the notion of two distinct subspaces,
$\c{U}$ and $\c{V}$, that capture the low-rank signal of a dataset $\m{X} \in
\b{R}^{n\times p}$. Here, classical PCA would only be able to compare subspaces of a
fixed, known dimension, $d_{\c{U}} = d_{\c{V}} = d$ to identify the right or
empirically ``best'' projection matrix with respect to some loss metric.
Modern statistical learning, however, necessitates metrics
that robustly handle dimensional mismatch ($d_{\c{U}} \neq d_{\c{V}}$).
Furthermore, in practice, the true effective dimension of the population
covariance is rarely known and must be estimated from the data, yielding an
empirical dimension $\hat{d}_{\c{U}} \neq \hat{d}_{\c{V}}$ that may or may not
include the true signal subspace. Consequently, establishing rigorous
statistical consistency requires quantifying both the
angular alignment of the subspaces and the rank estimation error simultaneously.
To motivate this work, we highlight here four particular statistical tasks in
which such a dimension mismatch problem arises naturally:
\begin{itemize}
\item \textbf{Adaptive Dimension Reduction and Model Selection}: When deciding
how many principal components to keep (e.g., looking at a scree plot or using an
information criteria), we are inherently comparing models of different complexities.
Measuring the distance between a 3-dimensional PCA
subspace and a 4-dimensional PCA subspace, for example, allows for a rigorous quantification
of exactly how much geometric information is gained or lost by changing the
rank threshold. Thus, it is natural to consider a nested candidate model
$\hat{\c{V}}_{[k]}$ of dimension $1 \leq k < d_{max} < p$ (forward nested reduction).

\item \textbf{Subspace Tracking in Streaming or Distributed Data}: When
  attempting to use PCA in online learning, where the covariance matrix is
  constantly updated based on new data points, the
  effective rank of the data manifold can expand or contract over time as new
  patterns emerge. Similarly, in distributed or federated PCA, different local
  servers might estimate different local ranks based on their private datasets. To
  aggregate or compare these local PCA models, we need a distance metric that
  can deal with a situation where Server A estimates a $d_A$-dimensional subspace and
  Server B estimates a $d_B$-dimensional subspace for $d_A \neq d_B$.

\item \textbf{Sequential Partitioning Schemes:}
  PCA is often cited as a tool to help simplify the construction of data-driven
  partitioning schemes when optimizing for a reduction in loss, like in the
  construction of generalized (oblique) decision trees \cite{tomita2020sparse, cattaneo2024convergence,
    menze2011oblique}. In such scenarios, the dataset changes
  at each sequential partitioning in the tree construction, requiring the
  flexibility of generating principal components of varying dimensionality to
  fit a given partition of the support. When selecting optimal partitioning rules
  in the construction of such tree estimators, it is necessary to work with a
  loss function that adapts to the dimensionality of the projection across various
  partitions.

\item \textbf{Differential privacy:} In addition to high-dimensional datasets,
  the concern around user privacy has grown across various communities. The use
  of PCA within differentially private models has become a recent point of
  interest in the literature
  \cite{yun2025high, he2025differentially, donhauser2024certified}.
  In order to maintain privacy measures it may be necessary to modify the subspace
  dimensions in a manner that requires analysts to compare subspace of different
  ranks in order to conduct inference downstream of the privacy-protected estimates.

\end{itemize}
By explicitly accommodating $d_{\c{U}} \neq d_{\c{V}}$, the ensuing distance
metrics should be able to naturally penalize an estimated subspace for entirely
missing variance components that reside in the orthogonal complement but belong
to the true signal space, providing a comprehensive and geometrically unified
measure of discrepancy. This leads us to ask two key questions:
\begin{enumerate}
\item Does there exist a risk function that operates on projections of two
  different dimensionalites that correctly identifies the true prediction loss?
\item Can we consistently measure the loss of a dimension mismatch?
 \end{enumerate}

We start by providing an overview of some recent work in the intersection of PCA
and high-dimensional statistical guarantees.
\subsection{Literature overview}
The idea of Principal Component Analysis (PCA) was first introduced by
\cite{pearson1901lines} and then subsequently studied in
\cite{hotelling1933analysis}. The most commonly referenced version of the
standard PCA, however, was first presented in \cite{jolliffe2002principal}. As
discussed earlier, the main advantage of PCA was its ability simplify a large
collection of covariates into a handful of distinct signals that captured a
sufficiently large proportion of the ``\textit{explainability}'' of the data,
mathematically defined in terms of the variance.

As the algorithm gained popularity, its applicability also evolved. In
particular, datasets became larger, data types became more varied and the desire
to have \textit{interpretable} models based on a finite collection of covariates
became more important. The technical concepts behind the PCA algorithm are based
on eigenvalue and eigenvector estimation that allow for dimension reduction
in a way that minimally impacts a pre-specified loss function. In this manner,
the explainability of the observations comes from the data itself, as opposed to
through assumptions or expert opinions that are not only difficult to quantify
in a methodological manner but can also induce an unknown bias into the system.
We refer interested readers to reviews of PCA
\cite[and references
therein]{jolliffe2016principal, hasan2021review, gewers2022principal}
for a more detailed history and
initial contributions to the methodology.
We focus the remainder of our literature review on recent advances in PCA-type
methods that are most relevant to our work.

The treatment of PCA approximations as a measure of subspace distances is
not new in itself.
\cite{shawe2005eigenspectrum, blanchard2007statistical} have shown bounds on the
empirical risk of projections using kernel PCA. In particular
\cite{blanchard2007statistical} have shown that the rate of convergence of these
projections can be
faster than $n^{-1/2}$ under certain spectral-gap conditions.
\cite{reiss2020nonasymptotic} show a $n^{-1}$ rate control by way of the
Davis-Kahan $\sin \Theta$ theorem. With some additional assumptions on the data
(like sub-Gaussianity), the authors also discuss asymptotic behavior of the
risk measure (see \cite[Proposition 2.14]{reiss2020nonasymptotic}).
\cite{reiss2020nonasymptotic, jirak2020perturbation, milbradt2020high} show
control in the infinite-dimensional setting (i.e., exponential or polynomial
decay of eigenvalues). The key difference in our work with respect to these
earlier efforts is that we show a non-trivial generalization of these results to
risk minimization with projections of different ranks, to allow the error to
capture not just the noise from finite-sample estimates, but also a potential
mismatch in rank estimation. We discuss further our asymptotic and
infinite-dimensional results in comparison to existing literature in
Section~\ref{sec:main_results}.

More recently, the perspective of parallel analysis, which began with
\cite{horn1965rationale}, has garnered interest from theoreticians on its ability
to effectively conduct factor analysis above a noise threshold. Recent
modifications of this procedure include deterministic parallel analysis (DPA)
which involves estimating the distance between singular values from the data
and corresponding singular values if the true data distribution was i.i.d.
Gaussian. This
comparison is used to identify if the estimated singular value from the data
behaves sufficiently different from pure noise. DPA involves direct
construction of the random noise singular values based on the
Marchenko-Pastur distribution \cite{dobriban2015efficient,
zhou2018eigenvalue}. Further improvements have been made on this approach in
\cite{dobriban2019deterministic, vitale2017selecting} to avoid making arbitrary
hyperparameter selections that can influence the identification of the true
number of singular values that form the subspace in which the data belongs.

Switching towards the high-dimensional literature wherein the dimension and
sample size grow proportionally, much work has been done to specifically
understand the spiked covariance model \cite{nadler2008finite,
paul2007asymptotics, hucker2023note}, with a focus on asymptotic behavior. Some
of these papers have also studied the bias that arises in the high-dimensional
setting which requires careful correction when the principal components are used
for subsequent regression analysis (PCR).
One the other hand, \cite{abiyounes2025lanczos} approach the spiked covariance
without direct estimation of the eigenvalues, rather using a Lanczos
algorithm type approach.

A recent alternative approach for dealing with the high-dimensional issues of the scree plot
directly has been proposed in \cite{ke2023estimation}. The authors
propose an approach (BEMA) to estimating the spiked covariance model in high
dimensions. The approach here is specifically intended to provide a clearer
understanding of scree plots, like the one shown in \cref{fig:bad_scree}, by
combining it with parallel analysis tools and quantile matching. While this
method has gained traction amongst the community, the specifications to the
spiked covariance model and the necessary parametric distribution assumption on
the noise can be limiting factors in its applicability.
\cite{wang2026analysis} have also recently shown perturbation bounds for low
rank matrices by extending the Davis-Kahan-Wedin theorem
\cite{wedin1972perturbation}. This approach, as most others,
only works for comparing matrices of the same rank.

\subsection{Main contributions}
In context of the state of the field, we provide here the main contributions of
this work to the literature that highlight the specific technical novelties that
we believe can be tools of interest more generally.

Theorem~\ref{thm:kahan:v2_generalized} presents a non-trivial generalization of
the Davis-Kahan theorem for matrices of different ranks. The result allows for a
novel comparison between covariances of different ranks, capturing both the
issue of unmatched dimensions and the leakage of information that exists outside
the reference subspace. Importantly, the result is independent of the specific
basis choice and rather is a fundamental statement about the projection spaces
with prescribed ranks. We show that the classical Davis-Kahan result in
\cite{yu2015useful} is recovered as a special case of
Theorem~\ref{thm:kahan:v2_generalized}.

Secondly, Theorem~\ref{thm:macro:bias:split} presents a convergence result
on the comparison of two empirical matrices of differing ranks. In particular,
if the reference subspace is now constructed through estimation on a
split-sample dataset, the risk metrics remain consistent and provide a strictly
convex loss function over the support of the reference dimension size. The
result allows us to identify precisely that a 50-50 sample split for
constructing the reference space and the estimated space is optimal
(Proposition~\ref{prop:half:ideal}).

Finally, Theorem \ref{thm:hd:dimension:consistency} shows that under the
conditions of Theorem~\ref{thm:macro:bias:split} and appropriate regularization,
the minimizer of the empirical risk function converges almost surely to the true
target subspace dimension. This theorem also provides the threshold above which
eigenvalue spikes can be detected by the proposed risk functions
(Corollary~\ref{cor:hd:rank:consistency}).

\subsection{Outline}
The remainder of this paper is structured as follows. Section~\ref{sec:prelim}
introduces the two proposed risk functions and some of their basic properties.
Section~\ref{sec:main_results} establishes the main results of this paper that
control the upper bound on the risk functions under eigen-gap conditions and
subsequently establishes their asymptotic distributional convergence.
Section~\ref{sec:high:dim} shows the fully implementable high-dimensional
application of the functions along with establishing their discrete convexity and
optimality of 50-50 sample splitting for the limiting criteria.
Finally Section~\ref{sec:simulations} presents a detailed simulations study on
synthetic data and on a lung cancer dataset comparing the proposed methods to
existing approaches in the literature.
The proofs of the results in Sections~\ref{sec:prelim}, ~\ref{sec:main_results}
and \ref{sec:high:dim} are provided in Appendices~\ref{app:prelim},
~\ref{app:main_proofs} and \ref{app:high:dim} respectively.

\section{Preliminaries}
\label{sec:prelim}
We start by introducing the key notation and conventions that will be
adopted throughout the paper.

\subsection{Notation}\label{ssec:not}
We write $[k] \defeq \{1,\ldots,k\}$ for any nonnegative integer $k \in \b{N}$, with $[0] \defeq \varnothing$. The cardinality of a finite set $\c{I}$ is denoted by $|\c{I}|$. We adopt the conventions that empty sums equal zero, $\inf \varnothing = +\infty$, $0\cdot\infty=0$, $0^{0}=1$, and $1/\infty = 0$.
For vectors, $\|\cdot\|_{2}$ denotes the Euclidean norm, and $\b{S}^{p-1} \defeq \{\m{x} \in \b{R}^{p} : \|\m{x}\|_{2}=1\}$. We write $\spann(\m{M})$ for the column space of a matrix $\m{M}$, with the span of an empty collection understood to be $\{0\}$. For $\m{M} \in \b{R}^{d_{1} \times d_{2}}$, let
$\sigma_{1}(\m{M}) \ge \cdots \ge \sigma_{r}(\m{M}) \ge 0$
be its singular values, where $r = \min(d_{1},d_{2})$.
For $q \in [1,\infty)$, its Schatten $q$-norm is
\begin{align*}
    \vertii{\m{M}}_{q}
    \defeq
    \left(\sum_{j=1}^{r}\sigma_{j}(\m{M})^{q}\right)^{1/q},
\end{align*}
and $\vertii{\m{M}}_{\infty} \defeq \sigma_{1}(\m{M})$
when $r \ge 1$.
All Schatten norms of empty matrices are defined to be zero.
The cases $q=1$, $q=2$, and $q=\infty$ give the trace, Frobenius,
and operator norms, respectively.
For every $q \in [1,\infty]$, we remark that $\vertiii{\m{M}}_{q} = \vertii{\m{M}^{\top}}_{q} = \vertii{(\m{M}^{\top}\m{M})^{1/2}}_{q}$.

For a given dimension $p \in \b{N}$, and an index set $\c{I} \subseteq [p]$, its complement is
$\c{I}^{c} \defeq [p] \setminus \c{I}$.
Let $\m{I}_{p} \in \b{R}^{p \times p}$ denote the identity matrix. We define the spaces of symmetric and non-increasingly ordered diagonal matrices as
\begin{align*}
    \c{S}(p) := \{\m{M} \in \b{R}^{p \times p} \mid \m{M} = \m{M}^{\top}\}, \quad \c{D}(p) := \{\diag(\lambda_{1}, \dots, \lambda_{p}) \in \b{R}^{p \times p} \mid \lambda_{1} \ge \cdots \ge \lambda_{p} \},
\end{align*}
A matrix $\m{M} \in \c{S}(p)$ is positive semidefinite if $\m{x}^{\top}\m{M}\m{x} \ge 0$ for every $\m{x} \in \b{R}^{p}$,
and positive definite if the inequality is strict for every $\m{x} \ne \m{0}$.
For $0 \le d \le p$, define
\begin{align*}
    \c{O}(p,d)
    \defeq
    \{\m{U} \in \b{R}^{p \times d}
        : \m{U}^{\top}\m{U} = \m{I}_{d}\},
\end{align*}
and we say that $\m{U} \in \c{O}(p,d)$ is semi-orthogonal. When $d = p$, we omit the column dimension and simply write $\c{O}(p)$ to denote the standard orthogonal group. For a trivial case where $d=0$, we let $\m{I}_{0}$ to be the empty $0 \times 0$ matrix, and $\c{O}(p,0)$ contains the unique $p \times 0$ matrix. We also set $\c{S}(0) = \c{D}(0) = \{\m{I}_{0}\}$.

For any eigendecomposition of a symmetric matrix, we assume by default that the eigenvalues are ordered non-increasingly. Specifically, given $\m{\Sigma}$ and a symmetric perturbation $\hat{\m{\Sigma}}$, we write $\m{\Sigma} = \m{U} \m{\Lambda} \m{U}^{\top}$ and $\hat{\m{\Sigma}} = \hat{\m{V}} \hat{\m{\Omega}} \hat{\m{V}}^{\top}$, where $\m{U}, \hat{\m{V}} \in \c{O}(p)$ and $\m{\Lambda} = \diag(\lambda_{1},\ldots,\lambda_{p}), \hat{\m{\Omega}} = \diag(\hat{\omega}_{1},\ldots,\hat{\omega}_{p}) \in \c{D}(p)$. The corresponding eigenvectors are denoted by $\m{u}_{j}$ and $\hat{\m{v}}_{j}$, respectively.
In the case where an eigenvalue has algebraic multiplicities, the corresponding eigenbases are fixed whenever needed and, for random matrices, chosen measurably.

\subsection{Space Similarity}
We are interested in characterizing the geometric relationship between projection matrices of varying dimensions. This will allow us to define a loss metric that naturally adapts to subspaces, potentially of different ranks. To this end, we first review the fundamental properties of projection matrices that underpin our proposed metrics.
Given a semi-orthogonal matrix $\m{U} \in \c{O}(p, d_{\c{U}})$ with $d_{\c{U}} \in [p]$, we denote its column space by $\c{U} = \spann (\m{U}) \subseteq \b{R}^{p}$. The corresponding orthogonal projection matrix onto this subspace is defined as $\m{\Pi}_{\c{U}} := \m{U} \m{U}^{\top} \in \c{S}(p)$.
We use the subscript $\c{U}$ rather than $\m{U}$ to emphasize that the projection depends solely on the subspace itself and is invariant to the choice of orthonormal basis. Hence, given a non-trivial subspace $\c{U} \subseteq \b{R}^{p}$ of dimension $d_{\c{U}} \in [p]$, it is natural to define the set of its corresponding semi-orthogonal matrices as
\begin{align*}
    \c{O}(\c{U}) := \{ \m{U} \in \c{O}(p, d_{\c{U}}) \mid \spann(\m{U}) = \c{U} \}.
\end{align*}
We apply this notation \textit{mutatis mutandis} to other subspaces, such as $\c{V}$ generated by $\m{V} \in \c{O}(p, d_{\c{V}})$.

Now, consider two non-trivial subspaces $\c{U}, \c{V} \subseteq \b{R}^{p}$ with dimensions $d_{\c{U}}$ and $d_{\c{V}}$, respectively. Let $\m{U} \in \c{O}(\c{U})$ and $\m{V} \in \c{O}(\c{V})$ be arbitrarily chosen semi-orthogonal basis matrices for these subspaces.
The following proposition relates the singular values of $\m{U}^{\top}\m{V}$ to the corresponding projection matrices.

\begin{proposition}\label{prop:matrix_product}
Let $\c{U},\c{V} \subseteq \b{R}^{p}$ be linear subspaces, and let
$\m{U} \in \c{O}(\c{U})$ and $\m{V} \in \c{O}(\c{V})$.
The nonzero singular values of $\m{U}^{\top}\m{V}$ and
$\m{\Pi}_{\c{U}}\m{\Pi}_{\c{V}}$ coincide, including multiplicities.
Consequently, the singular values of $\m{U}^{\top}\m{V}$ are
independent of the chosen orthonormal bases and belong to $[0,1]$.
Moreover, all its nonzero singular values equal one if and only if
\begin{align}\label{eq:SVD:commute}
    \m{\Pi}_{\c{U}}\m{\Pi}_{\c{V}}
    = \m{\Pi}_{\c{V}}\m{\Pi}_{\c{U}}
    = \m{\Pi}_{\c{U}\cap\c{V}}.
\end{align}
\end{proposition}

The primary consequence of Proposition~\ref{prop:matrix_product} is that the
singular values of the cross-matrix $\m{U}^{\top} \m{V} \in \b{R}^{d_{\c{U}} \times d_{\c{V}}}$ between any two $\m{U}$ and $\m{V}$ is completely determined by the product of their corresponding orthogonal projections. Because any unitarily invariant norm depends exclusively on singular values, this property ensures that geometric distances between subspaces constructed via such norms are strictly independent of the specific choice of basis.
In this regard, we note that \cite[Proposition 32]{wang2026analysis} serves as a counterpart to Proposition~\ref{prop:matrix_product} in the restricted case of equal dimensions ($d_{\c{U}} = d_{\c{V}}$). However, our established invariance holds regardless of dimensional equality, opening the door to applications requiring the comparison of arbitrary subspaces to capture the underlying behavior of the data.
This flexibility is particularly relevant to settings such as distributed principal component analysis (PCA), where local bases are estimated independently and must be systematically compared or aligned.

\begin{definition}\label{def:distances}
Let $\c{U}, \c{V} \subseteq \b{R}^{p}$ be linear subspaces with dimensions $d_{\c{U}}$ and $d_{\c{V}}$, respectively. For some (thus any) semi-orthogonal matrices $\m{U} \in \c{O}(\c{U})$ and $\m{V} \in \c{O}(\c{V})$, we define:
\begin{enumerate}
    \item The \textbf{Bures-Wasserstein distance} is defined as $l_{\c{W}}(\c{U}, \c{V}) := ( d_{\c{U}} + d_{\c{V}} - 2 \vertii{\m{U}^{\top} \m{V}}_{1} )^{1/2}$.
    \item The \textbf{chordal distance} is defined as $l_{\Theta}(\c{U}, \c{V}) := ( d_{\c{U}} + d_{\c{V}} - 2 \vertii{\m{U}^{\top} \m{V}}_{2}^{2} )^{1/2}$.
\end{enumerate}
In the trivial cases where either $d_{\c{U}} = 0$ or $d_{\c{V}} = 0$, we have $l_{\c{W}}(\c{U}, \c{V}) = l_{\Theta}(\c{U}, \c{V}) = ( d_{\c{U}} + d_{\c{V}})^{1/2}$.
\end{definition}

Equivalently, the distances can be expressed through the orthogonal projections as
\begin{align*}
    l_{\c{W}}(\c{U}, \c{V}) = \big( \tr[\m{\Pi}_{\c{V}} + \m{\Pi}_{\c{U}} - 2 (\m{\Pi}_{\c{U}}^{1/2} \m{\Pi}_{\c{V}} \m{\Pi}_{\c{U}}^{1/2})^{1/2}] \big)^{1/2}, \quad l_{\Theta}(\c{U}, \c{V}) = \vertiii{\m{\Pi}_{\c{U}} - \m{\Pi}_{\c{V}}}_{2},
\end{align*}
The first expression is the $2$-Wasserstein distance
squared between the centered, possibly degenerate Gaussian distributions with covariance matrices $\m{\Pi}_{\c{U}}$ and
$\m{\Pi}_{\c{V}}$ \cite{dowson1982frechet,givens1984class}, which is also referred to as the Procrustes
distance \cite{dryden2016statistical}.
Since a subspace is uniquely determined by its orthogonal projection, both expressions define metrics between subspaces. Writing $r \defeq \min(d_{\c{U}},d_{\c{V}})$ and the $j$-th singular value as $\sigma_{j} \defeq \sigma_{j}(\m{U}^{\top}\m{V})$, we have
\begin{align}\label{eq:prin:ang}
    l_{\Theta}(\c{U},\c{V})^{2}
    -l_{\c{W}}(\c{U},\c{V})^{2}
    =2\sum_{j=1}^{r}\sigma_{j}(1-\sigma_{j})
    \ge 0.
\end{align}
Thus, $l_{\c{W}}(\c{U},\c{V}) \le l_{\Theta}(\c{U},\c{V})$.
Equality holds precisely when every singular value belongs to
$\{0,1\}$, or equivalently when the projections commute in \eqref{eq:SVD:commute}.
To build intuition on how these metrics behave, we ilustrate two special cases:
\begin{enumerate}[leftmargin=*]
    \item \textbf{One-dimensional subspaces.}
    If $\c{U}=\spann(\m{u})$ and $\c{V}=\spann(\m{v})$
    for $\m{u},\m{v}\in\b{S}^{p-1}$, then
    \begin{align*}
        l_{\c{W}}(\c{U},\c{V})
        =
        \min_{s\in\{-1,1\}}\|\m{u}-s\m{v}\|_{2}, \qquad
        l_{\Theta}(\c{U},\c{V})
        =
        \sqrt{2}\sin(\theta),
        \qquad
        \theta\defeq\arccos|\m{u}^{\top}\m{v}|.
    \end{align*}

    \item \textbf{Nested subspaces.}
    If $\c{V}\subseteq\c{U}$, the cross-matrix has
    $d_{\c{V}}$ singular values, all equal to one, and
    \begin{align*}
        l_{\c{W}}(\c{U},\c{V})
        =l_{\Theta}(\c{U},\c{V})
        =\sqrt{d_{\c{U}}-d_{\c{V}}}.
    \end{align*}
\end{enumerate}

\begin{remark}[Principal angles]\label{rmk:prin:ang}
For subspaces of arbitrary dimensions, define their
$r=\min(d_{\c{U}},d_{\c{V}})$ principal angles by
\begin{align*}
    \theta_{j}
    \defeq
    \arccos\sigma_{j}(\m{U}^{\top}\m{V}),
    \qquad
    0\le\theta_{1}\le\cdots\le\theta_{r}\le\pi/2,
\end{align*}
and write
$\Theta(\c{U},\c{V})\defeq\diag(\theta_{1},\ldots,\theta_{r})$.
The distances admit the representations
\begin{align*}
    l_{\c{W}}(\c{U},\c{V})^{2}
    =
    |d_{\c{U}}-d_{\c{V}}|
    +4\sum_{j=1}^{r}\sin^{2}(\theta_{j}/2), \quad
    l_{\Theta}(\c{U},\c{V})^{2}
    =
    |d_{\c{U}}-d_{\c{V}}|
    +2\sum_{j=1}^{r}\sin^{2}(\theta_{j}).
\end{align*}
Thus each squared distance separates a dimension penalty from an angular discrepancy.
The projections commute precisely when every principal angle
is either $0$ or $\pi/2$.
When $d_{\c{U}}=d_{\c{V}}=d$, our normalization gives $l_{\Theta}(\c{U},\c{V}) = \sqrt{2} \vertii{\sin\Theta(\c{U},\c{V})}_{2}$,
where the sine is applied to the diagonal entries.
This relates our chordal distance to the usual
$\sin\Theta$ quantities in perturbation theory
\cite{davis1970rotation,yu2015useful}.
\end{remark}

The Bures-Wasserstein distance also admits an orthogonal Procrustes representation for unequal-dimensional subspaces.

\begin{theorem}[Orthogonal Procrustes]\label{thm:procrustes}
Without loss of generality, let $\c{U},\c{V}\subseteq\b{R}^{p}$ be linear
subspaces with $d_{\c{U}}\ge d_{\c{V}}$.
For any $\m{U}\in\c{O}(\c{U})$ and
$\m{V}\in\c{O}(\c{V})$, it holds that
\begin{align*}
    l_{\c{W}}(\c{U},\c{V})^{2}
    =
    \min_{\m{R}\in\c{O}(d_{\c{U}},d_{\c{V}})}
    \vertiii{\m{U}-\m{V}\m{R}^{\top}}_{2}^{2} \ge
    \min_{\m{R}\in\c{O}(d_{\c{U}},d_{\c{V}})}
    \vertiii{\m{U}\m{R}-\m{V}}_{2}^{2}
    =
    l_{\c{W}}(\c{U},\c{V})^{2}
    -(d_{\c{U}}-d_{\c{V}}).
\end{align*}
\end{theorem}

\section{Perturbation Bounds under Dimensional Mismatch}
\label{sec:main_results}
Having established a geometric framework for evaluating subspace similarities, we now turn to perturbation bounds under dimensional mismatch.
Throughout this section, the dimension $p$ is fixed, and $\hat{\m{\Sigma}} \in \c{S}(p)$ represents a perturbed matrix of some deterministic $\m{\Sigma} \in \c{S}(p)$.
We consider the eigendecompositions of $\m{\Sigma},\hat{\m{\Sigma}}\in\c{S}(p)$ as described in \cref{ssec:not}: $\m{\Sigma} = \m{U} \m{\Lambda} \m{U}^{\top}$ and $\hat{\m{\Sigma}} = \hat{\m{V}} \hat{\m{\Omega}} \hat{\m{V}}^{\top}$, where $\m{U}, \hat{\m{V}} \in \c{O}(p)$ and $\m{\Lambda}, \hat{\m{\Omega}} \in \c{D}(p)$.
For index sets $\c{I},\c{J}\subseteq[p]$, write
$d_{\c{U}}\defeq|\c{I}|$ and $d_{\c{V}}\defeq|\c{J}|$, and define
\begin{align*} 
    \m{U}_{\c{I}}
    &\defeq [\m{u}_{i}]_{i\in\c{I}}
    \in\c{O}(p,d_{\c{U}}), &
    \hat{\m{V}}_{\c{J}}
    &\defeq [\hat{\m{v}}_{j}]_{j\in\c{J}}
    \in\c{O}(p,d_{\c{V}}).
\end{align*}
Their column spaces are denoted by $\c{U}_{\c{I}} \defeq\spann(\m{U}_{\c{I}})$ and $\hat{\c{V}}_{\c{J}} \defeq\spann(\hat{\m{V}}_{\c{J}})$, respectively.
We also write $\m{\Lambda}_{\c{I}} \defeq\diag(\lambda_i)_{i\in\c{I}}$ and $\hat{\m{\Omega}}_{\c{J}} \defeq\diag(\hat{\omega}_j)_{j\in\c{J}}$, where
indices are listed in increasing order. These conventions apply to any index set, including complements, e.g., $\c{U}_{\c{I}^{c}}=\c{U}_{\c{I}}^{\perp}$ so that $\m{U}_{\c{I}}\m{U}_{\c{I}}^{\top}+\m{U}_{\c{I}^{c}}\m{U}_{\c{I}^{c}}^{\top}=\m{I}_{p}$.



Classical perturbation theory—from the foundational Davis-Kahan \cite{davis1970rotation} and Wedin \cite{wedin1972perturbation} theorems to modern variants \cite{wang2015singular, yu2015useful, wang2026analysis}—restricts attention to aligned principal subspaces, setting $\c{I} = \c{J} = [d]$ for a shared target rank $1 \le d \le p$. Consequently, upper bounds on the subspace distance between $\c{U}_{[d]}$ and $\c{V}_{[d]}$ (e.g., $\vertii{\sin \Theta (\c{U}_{[d]}, \hat{\c{V}}_{[d]})}_{2}$) are expressed in terms of the perturbation norm $\vertii{\hat{\m{\Sigma}} - \m{\Sigma}}_{2}$ and the spectral gap separating the first $d$ eigenvalues from the remainder of the spectrum.
Restricting the analysis to $\c{I} = \c{J} =
[d]$ inherently limits perturbation theory to quantifying the estimation error
of sample principal components relative to their population counterparts.
In contrast to the literature, we relax the index sets $\c{I}$ and $\c{J}$ to be
of arbitrary size, thus explicitly permitting dimensional mismatch ($d_{\c{U}}
\neq d_{\c{V}}$). While this generalized formulation introduces non-trivial
technical challenges in bounding the norm cross-matrix, it delivers substantial
practical and theoretical utility. By
permitting arbitrary index sets, our framework naturally quantifies the
empirical loss associated with projecting data onto candidate subspaces of
varying dimensions. This flexibility provides the core mathematical foundation
for downstream applications such as adaptive rank selection and subspace
alignment across heterogeneous nodes in distributed statistical learning.

\subsection{Deterministic Worst Bound}

Handling arbitrary and mismatched index sets introduces non-trivial technical
challenges. Specifically, the Frobenius norm of the cross matrix
$\vertii{\m{U}_{\c{I}^{c}}^{\top}\hat{\m{V}}_{\c{J}}}_{2}^{2}$ measures the
total squared component of the selected empirical eigenvectors outside the
reference subspace $\c{U}_{\c{I}}$. Thus, bounding this leakage requires a
careful analysis of the associated generalized Sylvester equation. This
necessitates carefully tracking the asymmetric dimensional scaling factors
across Schatten norms. We first begin by presenting a generalization of the
Davis-Kahan bound under dimensional mismatch, using the separation between the
selected empirical eigenvalues and the complementary population eigenvalues.

\begin{proposition}\label{prop:kahan:v1}
Let $\m{\Sigma},\hat{\m{\Sigma}}\in\c{S}(p)$ and
$\c{I},\c{J}\subseteq[p]$, with
$d_{\c{U}}\defeq|\c{I}|$ and $d_{\c{V}}\defeq|\c{J}|$.
Define
\begin{align*}
    \delta := \min_{i \in \c{I}^{c}} \min_{j \in \c{J}} |\lambda_{i} - \hat{\omega}_{j}|.
\end{align*}
Then, for every $q\in[2,\infty]$, it holds that
\begin{align*}
    \delta
    \vertii{\m{U}_{\c{I}^{c}}^{\top}
        \hat{\m{V}}_{\c{J}}}_{2}
    \le
    \min(p-d_{\c{U}},d_{\c{V}})^{1/2-1/q}
    \vertii{\hat{\m{\Sigma}}-\m{\Sigma}}_{q}.
\end{align*}
\end{proposition}

\begin{remark}[Sharpness of the Bound]\label{rmk:sharpness}
The inequality in Proposition~\ref{prop:kahan:v1} is sharp for every $q\in[2,\infty]$.
To see this, fix $1\le m<p$ and $\delta>0$, and take $\hat{\m{\Sigma}}=\m{I}_{p}$ and
\begin{align*}
    \m{\Sigma}
    =
    \diag(
        \underbrace{1+\delta,\ldots,1+\delta}_{m},
        \underbrace{1,\ldots,1}_{p-m}
    ).
\end{align*}
Choose $\m{U}=\hat{\m{V}}=\m{I}_{p}$,
$\c{I}=[p]\setminus[m]$, and $\c{J}=[m]$.
The spectral gap in Proposition~\ref{prop:kahan:v1} is then
$\delta$, and $\m{U}_{\c{I}^{c}}^{\top}\hat{\m{V}}_{\c{J}}=\m{I}_{m}$, thus both sides of the bound equal $\delta\sqrt{m}$.
\end{remark}

Typically, we consider settings where $d_{\c{U}}, d_{\c{V}} \ll p$, so the explicit dimension factor
reduces to $d_{\c{V}}^{1/2-1/q}$.
Thus, this factor depends on the estimated dimension rather than directly on $p$, although the perturbation norm may still depend on the ambient dimension.
To express the result in terms of chordal distance, observe that
\begin{align*}
    \vertii{
        \m{U}_{\c{I}^{c}}^{\top}\hat{\m{V}}_{\c{J}}
    }_{2}^{2}
    =
    \tr\big[
        \hat{\m{V}}_{\c{J}}^{\top}
        (\m{I}_{p}-\m{U}_{\c{I}}\m{U}_{\c{I}}^{\top})
        \hat{\m{V}}_{\c{J}}
    \big] =
    d_{\c{V}}
    -\vertii{\m{U}_{\c{I}}^{\top}\hat{\m{V}}_{\c{J}}}_{2}^{2}.
\end{align*}
Substitution into Definition~\ref{def:distances} gives the exact identity for the excess risk:
\begin{align}\label{eq:chordal:iden}
    0 \le l_{\Theta}(\c{U}_{\c{I}},\hat{\c{V}}_{\c{J}})^{2} - (d_{\c{U}}-d_{\c{V}})
    = 2\vertii{\m{U}_{\c{I}^{c}}^{\top} \hat{\m{V}}_{\c{J}}}_{2}^{2}.
\end{align}
When $d_{\c{U}} > d_{\c{V}}$, the first term is the unavoidable dimension penalty, and the second measures departure from subspace inclusion.

While Proposition~\ref{prop:kahan:v1} provides a tight deterministic bound, the gap $\delta$ involves both spectra of $\m{\Sigma}$ and $\hat{\m{\Sigma}}$.
For a random perturbation, controlling this gap may therefore
require a separate bound on empirical eigenvalue fluctuations.
In this regard, the following theorem instead uses a gap determined entirely by $\m{\Sigma}$.
Its positivity forces $\c{J}\subseteq\c{I}$, so the theorem
allows unequal ranks while requiring that the empirical target
indices lie within the population reference indices.

\begin{theorem}\label{thm:kahan:v2_generalized}
Let $\m{\Sigma},\hat{\m{\Sigma}}\in\c{S}(p)$ and
$\c{I},\c{J}\subseteq[p]$, with
$d_{\c{U}}\defeq|\c{I}|$ and $d_{\c{V}}\defeq|\c{J}|$. Assume
\begin{align*}
    \delta_{\m{\Sigma}} \defeq \min_{i \in \c{I}^{c}} \min_{j \in \c{J}} |\lambda_i-\lambda_j| > 0.
\end{align*}
Then $\c{J}\subseteq\c{I}$, and, for every $q\in[2,\infty]$, it holds that
\begin{align*}
    0 \le l_{\Theta}(\c{U}_{\c{I}},\hat{\c{V}}_{\c{J}})^{2}-(d_{\c{U}}-d_{\c{V}})
    \le \frac{8d_{\c{V}}^{1-2/q}}{\delta_{\m{\Sigma}}^{2}} \vertii{\hat{\m{\Sigma}}-\m{\Sigma}}_{q}^{2}.
\end{align*}
\end{theorem}

\begin{remark}
The inclusion $\c{J}\subseteq\c{I}$ is not merely for technical convenience, but is strictly necessary for an upper bound of this form. To see this, take $\hat{\m{\Sigma}}=\m{\Sigma}$ with the same eigenbases so that
\begin{align*}
    l_{\Theta}(\c{U}_{\c{I}},\c{U}_{\c{J}})^{2}-(d_{\c{U}}-d_{\c{V}})
    = 2|\c{J}\setminus\c{I}|, \qquad \c{I}, \c{J} \subseteq [p].
\end{align*}
Thus, if $\c{J} \not\subseteq \c{I}$, each target index outside $\c{I}$ contributes to excess loss even though $\vertii{\hat{\m{\Sigma}} - \m{\Sigma}}_{q} = 0$.
\end{remark}

Theorem~\ref{thm:kahan:v2_generalized} separates two sources of subspace discrepancy. The term $d_{\c{U}}-d_{\c{V}}$ is an unavoidable geometric penalty for unmatched dimensions, while the excess squared distance measures leakage outside the reference subspace. Note that repeated eigenvalues can make the selected subspaces non-unique when an index set includes only part of an eigenspace. Despite this ambiguity, Theorem~\ref{thm:kahan:v2_generalized} holds uniformly over all orthonormal eigenbases compatible with the ordered eigendecompositions. Its upper bound is therefore independent of these choices, even when the excess distance itself is not. Furthermore, the principal-angle representations in \eqref{eq:prin:ang} imply
\begin{align*}
    0 \le l_{\c{W}}(\c{U}_{\c{I}},\hat{\c{V}}_{\c{J}})^{2}-(d_{\c{U}}-d_{\c{V}})
    \le l_{\Theta}(\c{U}_{\c{I}},\hat{\c{V}}_{\c{J}})^{2}-(d_{\c{U}}-d_{\c{V}}),
\end{align*}
hence the same upper bounds directly apply for the excess squared Bures--Wasserstein distance.

For a random estimator $\hat{\m{\Sigma}}$, the bound separates the population spectral condition from the stochastic control of $\vertii{\hat{\m{\Sigma}}-\m{\Sigma}}_{q}$, avoiding the need to control an empirical eigengap separately. Combining the cases $q=2$ and $q=\infty$ gives
\begin{align*}
    l_{\Theta}(\c{U}_{\c{I}},\hat{\c{V}}_{\c{J}})^{2}-(d_{\c{U}}-d_{\c{V}})
    \le \frac{8}{\delta_{\m{\Sigma}}^{2}}
    \min\left\{\vertii{\hat{\m{\Sigma}}-\m{\Sigma}}_{2}^{2},
    d_{\c{V}}\vertii{\hat{\m{\Sigma}}-\m{\Sigma}}_{\infty}^{2}\right\}.
\end{align*}
For $\c{I}=\c{J}=\{r,\ldots,s\}$, this recovers the $\sin\Theta$ bound of \cite[Theorem~2]{yu2015useful}. When $d_{\c{V}}=1$, the dimension factors disappear and $q=\infty$ gives the smallest upper bound among the Schatten bounds in both Proposition~\ref{prop:kahan:v1} and Theorem~\ref{thm:kahan:v2_generalized}.

\subsection{Asymptotic Behavior}
We now complement the deterministic bound in Theorem~\ref{thm:kahan:v2_generalized} by identifying the limiting distribution of the excess squared chordal distance. This describes the fluctuations of the subspace discrepancy around the dimension penalty $d_{\c{U}}-d_{\c{V}}$ as the perturbation vanishes. While classical asymptotic theory studies eigenvectors and principal components \cite{tyler1981asymptotic,muirhead1982aspects}, we consider excess subspace distances under dimensional mismatch. Our framework accommodates repeated population eigenvalues subject to an additional spectral separation condition.
The following result imposes no requirement that $\hat{\m{\Sigma}}_{n}$ be a sample covariance matrix or that the limiting noise be Gaussian, thereby encompassing a broad class of heavy-tailed and non-Gaussian matrix estimators.

\begin{theorem}\label{thm:general_asymptotic_dist}
Let $p$ be fixed, let $\m{\Sigma}\in\c{S}(p)$ be deterministic, and let
$\varnothing\ne\c{J}\subseteq\c{I}\subseteq[p]$ have cardinalities
$|\c{I}|=d_{\c{U}}$ and $|\c{J}|=d_{\c{V}}$. Assume
\begin{align*}
    \tilde{\delta}_{\m{\Sigma}}
    \defeq \min_{i \in \c{J}^{c}} \min_{j \in \c{J}} |\lambda_{i} - \lambda_{j}|>0.
\end{align*}
Let $\hat{\m{\Sigma}}_n\in\c{S}(p)$ be symmetric random matrices satisfying $n^{\zeta}(\hat{\m{\Sigma}}_n-\m{\Sigma})\xrightarrow{d}\m{Z}$
for some $\zeta>0$ and a random matrix $\m{Z}\in\c{S}(p)$. Then
\begin{align*}
    n^{2\zeta}\left[l_{\Theta}(\c{U}_{\c{I}},\hat{\c{V}}_{\c{J}})^2-(d_{\c{U}}-d_{\c{V}})\right]
    \xrightarrow{d}
    2\sum_{i\in\c{I}^{c}}\sum_{j\in\c{J}}
    \left(\frac{\tilde{\m{Z}}_{ij}}{\lambda_j-\lambda_i}\right)^2,
\end{align*}
where $\tilde{\m{Z}} := \m{U}_{\c{I}^{c}}^{\top} \m{Z} \m{U}_{\c{J}} \in \b{R}^{(p-d_{\c{U}}) \times d_{\c{V}}}$.
\end{theorem}

\begin{remark}[Convergence failure under a tied perturbation]
In Theorem~\ref{thm:general_asymptotic_dist}, the weaker condition $\delta_{\m{\Sigma}}>0$ is not sufficient to guarantee convergence in distribution, unlike Theorem~\ref{thm:kahan:v2_generalized}. If $\c{J}$ splits a repeated population eigenspace and the first-order perturbation leaves the relevant directions tied, smaller-order terms can determine which directions are selected.
To illustrate this, let $t_n=n^{-\zeta}$ for some $\zeta>0$, and consider
\begin{align*}
    \m{\Sigma}=\diag(2,2,1),
    \qquad
    \hat{\m{\Sigma}}_n=
    \begin{pmatrix}
        2+(-1)^n t_n^{3/2} & 0 & 0\\
        0 & 2-(-1)^n t_n^{3/2} & t_n\\
        0 & t_n & 1
    \end{pmatrix},
    \qquad
    \c{I}=\{1,2\},\quad \c{J}=\{1\}.
\end{align*}
Here $\delta_{\m{\Sigma}}=1$, whereas $\tilde{\delta}_{\m{\Sigma}}=0$. Moreover, the scaled perturbations converge deterministically:
\begin{align*}
    n^{\zeta}(\hat{\m{\Sigma}}_n-\m{\Sigma})
    \, \xrightarrow{d} \,
    \m{Z} =
    \begin{pmatrix}
        0&0&0\\
        0&0&1\\
        0&1&0
    \end{pmatrix}.
\end{align*}
For all sufficiently large even $n = 2m$, the leading eigenvector of $\m{\Sigma}$ is $\m{e}_1$, which lies entirely in $\c{U}_{\c{I}}$. For odd $n = 2m+1$, the leading eigenvector can be chosen as
\begin{align*}
    \hat{\m{v}}_1=(0,\cos\theta_n,\sin\theta_n)^{\top},
    \qquad
    \theta_n=\frac12\arctan\left(\frac{2t_n}{1+t_n^{3/2}}\right).
\end{align*}
Consequently, if we denote the scaled excess risk by $T_n\defeq n^{2\zeta} [l_{\Theta}(\c{U}_{\c{I}},\hat{\c{V}}_{\c{J}})^2-1 ]$,
we obtain
\begin{align*}
    T_{2m}=0, \qquad    T_{2m+1}=2t_{2m+1}^{-2}\sin^2\theta_{2m+1}\longrightarrow 2, \qquad \text{for all sufficiently large }m.
\end{align*}
Thus $T_n$ has no limiting distribution, despite convergence of the scaled perturbations.
\end{remark}

For the sample covariance matrix, finite fourth moments yield a Gaussian perturbation limit through the multivariate central limit theorem, hence Theorem~\ref{thm:general_asymptotic_dist} gives the following specialization.

\begin{corollary}\label{cor:asymp_sample_cov}
Let $p$ be fixed, and let $X_1,\ldots,X_n\in\b{R}^{p}$ be i.i.d. random vectors with mean $\mu$, covariance matrix $\m{\Sigma}$, and finite fourth moments. Define the sample covariance by
\begin{align*}
    \hat{\m{\Sigma}}_n\defeq\frac1n\sum_{r=1}^n(X_r-\bar{X}_n)(X_r-\bar{X}_n)^{\top} \in \c{S}(p), \qquad \bar{X}_n\defeq\frac1n\sum_{r=1}^n X_r.
\end{align*}
Under the assumptions of Theorem~\ref{thm:general_asymptotic_dist}, in particular $\tilde{\delta}_{\m{\Sigma}} = \min_{i \in \c{J}^{c}} \min_{j \in \c{J}} |\lambda_{i} - \lambda_{j}|>0$, we have
\begin{align*}
    n\left[l_{\Theta}(\c{U}_{\c{I}},\hat{\c{V}}_{\c{J}})^2-(d_{\c{U}}-d_{\c{V}})\right]
    \xrightarrow{d}
    2\sum_{i\in\c{I}^{c}}\sum_{j\in\c{J}}
    \left(\frac{\tilde{\m{Z}}_{ij}}{\lambda_j-\lambda_i}\right)^2,
\end{align*}
where $\tilde{\m{Z}}\in\b{R}^{(p-d_{\c{U}})\times d_{\c{V}}}$ is a centered Gaussian random matrix with
\begin{align*}
    \operatorname{Cov}[\operatorname{vec}(\tilde{\m{Z}})]
    =(\m{U}_{\c{J}}^{\top}\otimes\m{U}_{\c{I}^{c}}^{\top})\m{C}(\m{U}_{\c{J}}\otimes\m{U}_{\c{I}^{c}}),\qquad
    \m{C}
    \defeq\operatorname{Cov}\left[\operatorname{vec}\left((X_1-\mu)(X_1-\mu)^{\top}\right)\right].
\end{align*}
The same conclusion holds when the sample covariance is normalized by $n-1$.
\end{corollary}

In general, the entries of $\tilde{\m{Z}}$ may be correlated. Under elliptical symmetry, their covariance is diagonal, and the limiting excess distance becomes an explicit weighted sum of independent chi-squared variables.

\begin{figure}[h]
    \centering
    \includegraphics[width=0.5\textwidth]{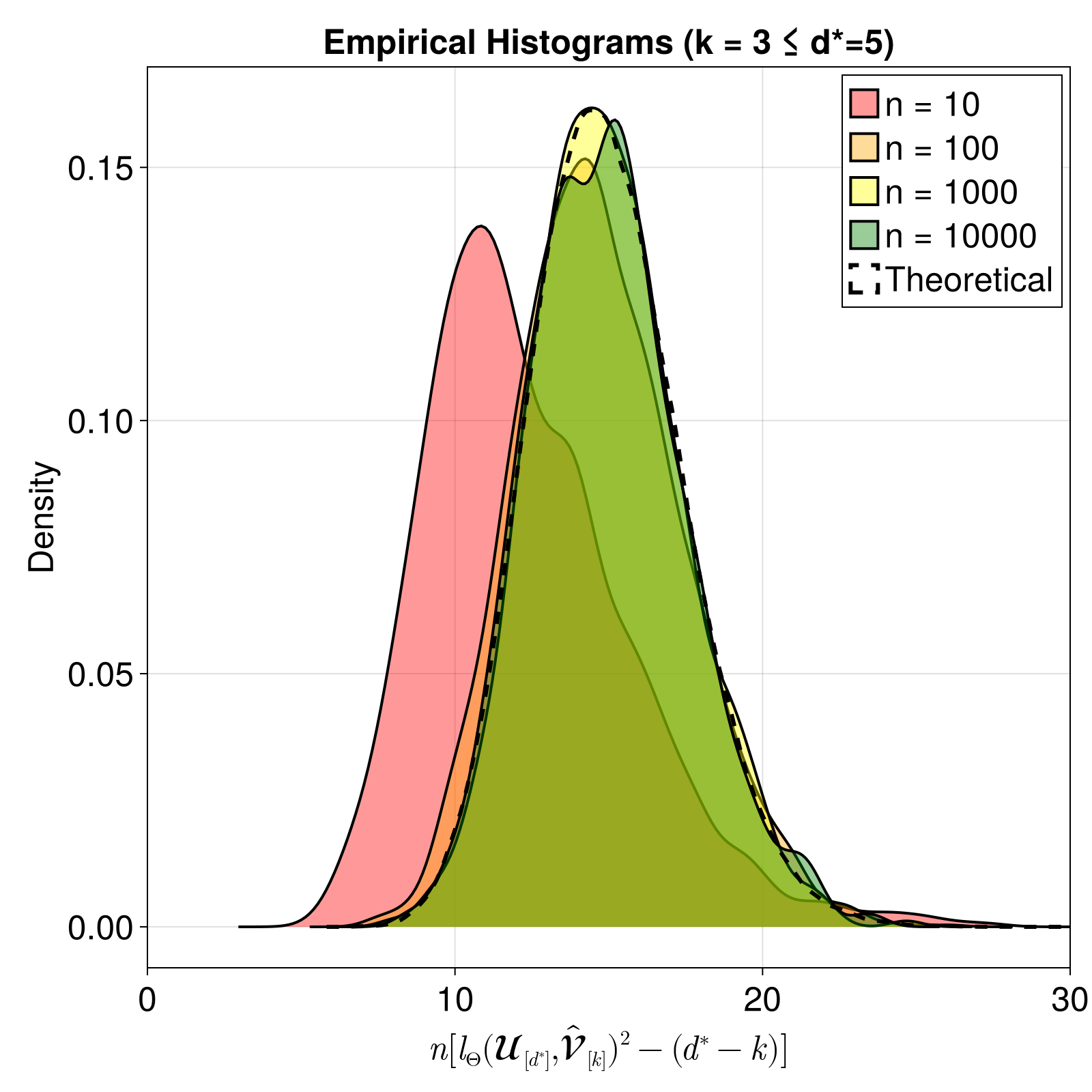}
    \caption{Empirical densities of the scaled excess squared chordal distance and its limiting density in \eqref{eq:chi:limit:ex}, for a Gaussian spiked covariance model with $p=30$, $d^{*}=5$, $k=3$, and population eigenvalues $\m{\Lambda} = \diag(15, 12,
    10, 8, 5, 1, \dots, 1)$.}
    \label{fig:asymp}
\end{figure}

\begin{corollary}\label{cor:asymp:ellip}
  Under the assumptions of Corollary~\ref{cor:asymp_sample_cov}, suppose additionally that the observations are elliptically distributed and $\m{\Sigma}$ is positive definite. Define the elliptical kurtosis parameter by
  \begin{align*}
    \kappa\defeq\frac{\b{E}[R^4]}{p(p+2)}-1,
    \qquad
    R^2\defeq(X_1-\mu)^{\top}\m{\Sigma}^{-1}(X_1-\mu).
  \end{align*}
  Then
  \begin{align*}
    n\left[l_{\Theta}(\c{U}_{\c{I}},\hat{\c{V}}_{\c{J}})^2-(d_{\c{U}}-d_{\c{V}})\right]
    \xrightarrow{d}
    2(1+\kappa)\sum_{i\in\c{I}^{c}}\sum_{j\in\c{J}}
    \frac{\lambda_i\lambda_j}{(\lambda_j-\lambda_i)^2}\chi^2_{1,ij},
    \qquad
    \chi^2_{1,ij}\stackrel{\text{i.i.d.}}{\sim}\chi^2_1.
  \end{align*}
\end{corollary}
Note that the result of Corollary~\ref{cor:asymp:ellip} is a
generalization of \cite[Proposition 2.14]{reiss2020nonasymptotic}, both in terms
of the distributional assumption (finite fourth moment versus Gaussian) on $X$
and on the subsequent distributional convergence for matrices of different ranks.
As an illustration of this result, consider Gaussian
observations under a spiked covariance
model with $1\le d^{*}<p$ signal directions and population eigenvalues
\begin{align*}
    \lambda_1>\lambda_2>\cdots>\lambda_{d^{*}}>
    \lambda_{d^{*}+1}=\cdots=\lambda_p=1.
\end{align*}
For $\c{J}=[k]\subseteq\c{I}=[d^{*}]$, where $1\le k\le d^{*}$, the required separation is $\tilde{\delta}_{\m{\Sigma}}=\lambda_k-\lambda_{k+1}>0$. Since $\kappa=0$ for Gaussian observations, Corollary~\ref{cor:asymp:ellip} gives
\begin{align}\label{eq:chi:limit:ex}
    n\left[l_{\Theta}(\c{U}_{[d^{*}]},\hat{\c{V}}_{[k]})^2-(d^{*}-k)\right]
    \xrightarrow{d}
    2\sum_{j=1}^k\frac{\lambda_j}{(\lambda_j-1)^2}Y_j,
    \qquad
    Y_j\stackrel{\text{i.i.d.}}{\sim}\chi^2_{p-d^{*}}.
\end{align}
Here each $Y_j$ collects the $p-d^{*}$ independent chi-squared terms associated with leakage from the $j$th selected direction into the noise subspace. The discrepancy therefore consists of the dimension penalty $d^{*}-k$ and an excess term of order $n^{-1}$, whose limiting distribution depends on the selected spikes and the dimension of the noise subspace. We illustrate this limit with $p=30$, signal dimension $d^{*}=5$, selected dimension $k=3$, and spike eigenvalues $(15,12,10,8,5)$. Figure~\ref{fig:asymp} compares the empirical densities of the scaled excess distance at increasing sample sizes with the limiting density in \eqref{eq:chi:limit:ex}.

\section{High-Dimensional Spiked Model}\label{sec:high:dim}
Thus far, we have considered the ambient dimension to be fixed. In modern
datasets, however, the number of features $p$ may be comparable to the sample size $n$, motivating an asymptotic regime in which both variables grow together
\cite{sur2019modern, candes2020phase}. In this regime, sampling noise can
obscure population signal directions, and even an empirical direction associated
with a separated eigenvalue can remain imperfectly aligned with the population
signal subspace. We use the distances developed earlier to construct a dimension-selection criterion from the agreement between independently estimated subspaces. Our statistical objective is to estimate the number of supercritical population spikes. The results below characterize which spikes the criterion retains, establish consistency under explicit signal-strength conditions, and identify the additional detection requirement introduced by sample splitting.

Consider the data matrix $\m{X} = \m{Z} \m{\Sigma}^{1/2} \in \b{R}^{n \times p}$ under the following assumptions:
\begin{itemize}[leftmargin=*]
    \item[(A1)] $p = p(n)$ satisfies $p/n \to \gamma \in (0, \infty)$, where $\gamma$ is the limiting aspect ratio.

    \item[(A2)] The population covariance has a fixed, unknown number $d^{*}$ of upper spikes, counted with multiplicity:
    \begin{align}\label{eq:model:spike}
        \m{\Sigma} = \m{I}_{p} + \sum_{j = 1}^{d^{*}} (\lambda_{j} - 1) \m{u}_{j} \m{u}_{j}^{\top}, \qquad \lambda_{1} \ge \cdots \ge \lambda_{d^{*}} > 1.
    \end{align}
    The spike eigenvalues are fixed, the population eigenvectors are deterministic and orthonormal, and the remaining eigenvalues satisfy $\lambda_{d^{*} + 1} = \cdots = \lambda_{p} = 1$. The eigenvectors may depend on $p$, and repeated spike eigenvalues are permitted.

    \item[(A3)] The entries of $\m{Z}$ are independent and identically distributed according to a real distribution that does not depend on $n$, with
    \begin{align*}
        \b{E} [Z_{11}] = 0, \qquad \b{E} [Z_{11}^{2}] = 1, \qquad \b{E} [|Z_{11}|^{m}] < \infty \quad \text{for every integer } m \ge 1.
    \end{align*}
\end{itemize}
When $d^{*} = 0$, we take $\m{\Sigma} = \m{I}_{p}$. The population signal subspace is $\c{U}_{[d^{*}]} \defeq \spann \{ \m{u}_{1}, \ldots, \m{u}_{d^{*}} \}$, with $\c{U}_{[0]} = \{ \m{0} \}$.
The unit noise variance is a normalization. Replacing $\m{X}$ by $\sigma \m{X}$, for any unknown $\sigma > 0$, changes the covariance to $\sigma^{2} \m{\Sigma}$ but leaves all empirical subspaces and the proposed geometric criteria unchanged. Under this scaling, the $\lambda_{j}$ represent population eigenvalues relative to the noise variance.
Under this normalization, the full-sample BBP (Baik-Ben Arous-P\'{e}ch\'{e}) phase transitions
\cite{paul2007asymptotics, benaych2011eigenvalues} occurs at $1 + \sqrt{\gamma}$. Accordingly, define the number of supercritical population spikes by
\begin{align}\label{eq:hd:supercritical:dimension}
    d_{\gamma}^{*} \defeq \# \{ 1 \le j \le d^{*} : \lambda_{j} > 1 + \sqrt{\gamma} \}.
\end{align}
Thus $d_{\gamma}^{*} = d^{*}$ when every population spike is supercritical.
Finally, the moment conditions in (A3) is a convenient sufficient condition for
the spectral results used in our framework.

\subsection{Population Benchmark}\label{ssec:oracle:benchmark}

A dimension-selection criterion must balance the evidence for an additional
signal direction against the possibility of retaining noise. The nonequal-dimensional
subspace distances, both chordal and Bures-Wasserstein, provide a geometric way
to express this balance: a candidate can omit signal directions, include
additional directions, or be poorly aligned with the reference subspace. We
first examine these effects using the population signal subspace as an oracle
reference. This benchmark explains the shape of the resulting loss curves and
why minimizing an non-regularized distance can retain fewer directions than the true
number of supercritical spikes.

For an aspect ratio $\gamma > 0$, define the limiting squared alignment associated with the BBP transition:
\begin{align*}
    g (\lambda ; \gamma) \defeq
    \begin{cases}
        \displaystyle \frac{1 - \gamma / (\lambda - 1)^{2}}{1 + \gamma / (\lambda - 1)}, & \lambda > 1 + \sqrt{\gamma}, \\[6pt]
        0, & 1 \le \lambda \le 1 + \sqrt{\gamma}.
    \end{cases}
\end{align*}
Above the transition, $g (\lambda ; \gamma)$ is strictly increasing in $\lambda$ and takes values in $(0, 1)$. At or below the transition, the limiting alignment vanishes \cite{paul2007asymptotics, benaych2011eigenvalues}.
Let $\hat{\m{\Sigma}} = n^{-1} \m{X}^{\top} \m{X}$ have ordered eigenvalues $\hat{\omega}_{1} \ge \cdots \ge \hat{\omega}_{p}$ and corresponding orthonormal eigenvectors $\hat{\m{v}}_{1}, \ldots, \hat{\m{v}}_{p}$. For $1 \le k \le p$, write
\begin{align*}
    \hat{\m{V}}_{[k]} \defeq [\hat{\m{v}}_{1}, \ldots, \hat{\m{v}}_{k}], \qquad \hat{\c{V}}_{[k]} \defeq \spann \{ \hat{\m{v}}_{1}, \ldots, \hat{\m{v}}_{k} \}.
\end{align*}
When a population spike is repeated, individual empirical eigenvectors need not align with a prescribed population basis. The appropriate limits therefore concern projections onto the entire population eigenspace. For a maximal block $\c{B} \subseteq [d^{*}]$ of equal spike eigenvalues, let $\lambda_{\c{B}}$ denote their common value and define the projection to the corresponding eigenspace as $\m{\Pi}_{\c{B}} \defeq \sum_{j \in \c{B}} \m{u}_{j} \m{u}_{j}^{\top}$. Then, the spiked covariance eigenspace limits give, for every fixed $k \ge 1$,
\begin{align}\label{eq:hd:population:overlap}
    \hat{\m{V}}_{[k]}^{\top} \m{\Pi}_{\c{B}} \hat{\m{V}}_{[k]} \xrightarrow{a.s.} g (\lambda_{\c{B}} ; \gamma) \diag \big( \mathbf{1}_{\{ j \in \c{B} \}} \big)_{j = 1}^{k};
\end{align}
see \cite[Theorems~2.11, 2.16, and~2.17]{bloemendal2016principal}. This statement also applies when $k$ ends inside a repeated spike block. It avoids choosing a limiting orientation within that eigenspace.

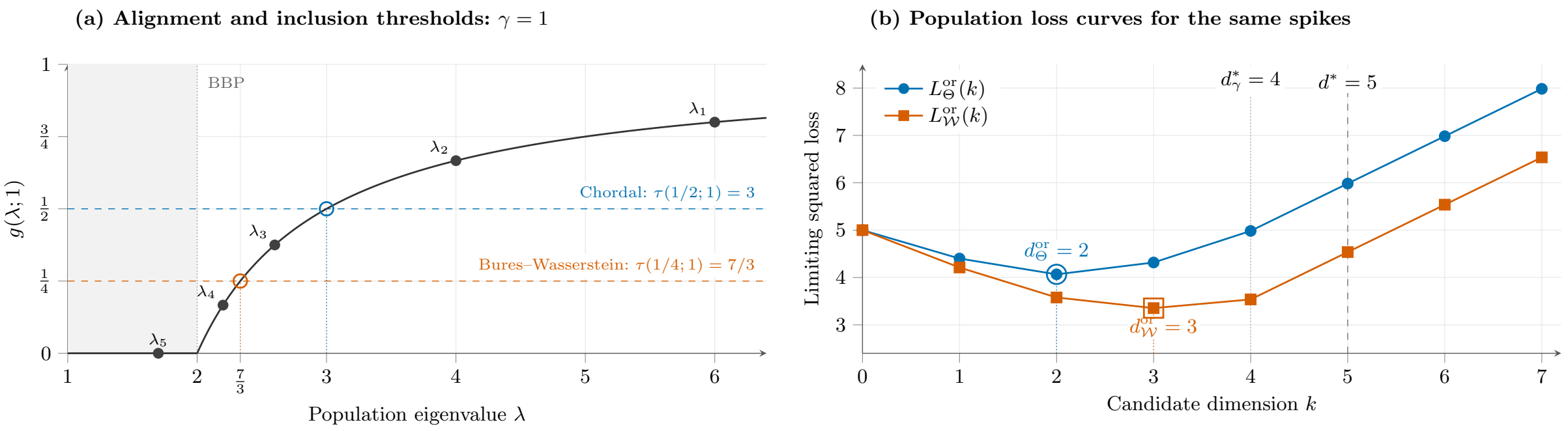
\begin{figure}[t]
\centering
\definecolor{oracleTheta}{HTML}{0072B2}
\definecolor{oracleW}{HTML}{D55E00}

\pgfplotstableread{
k Theta Wasserstein
0 5.000000000 5.000000000
1 4.400000000 4.211145618
2 4.066666667 3.578152456
3 4.316666667 3.353407585
4 4.983333333 3.536911004
5 5.983333333 4.536911004
6 6.983333333 5.536911004
7 7.983333333 6.536911004
}\OracleLossTable

\resizebox{\linewidth}{!}{
\begin{tikzpicture}
\begin{groupplot}[
    group style={group size=2 by 1, horizontal sep=1.5cm},
    width=12.5cm, height=6.1cm,
    axis lines=left,
    tick align=outside,
    tick label style={font=\small},
    label style={font=\small},
    title style={font=\small\bfseries, at={(0,1.05)}, anchor=south west},
    grid=major,
    grid style={black!8},
    axis line style={black!65},
    clip=false
]

\nextgroupplot[
    title={(a) Alignment and inclusion thresholds: $\gamma = 1$},
    xmin=1, xmax=6.4, ymin=0, ymax=1,
    xlabel={Population eigenvalue $\lambda$},
    ylabel={$g (\lambda ; 1)$},
    xtick={1,2,2.333333333,3,4,5,6},
    xticklabels={$1$,$2$,$\frac{7}{3}$,$3$,$4$,$5$,$6$},
    ytick={0,0.25,0.5,0.75,1},
    yticklabels={$0$,$\frac{1}{4}$,$\frac{1}{2}$,$\frac{3}{4}$,$1$}
]
\path[fill=black!5] (axis cs:1,0) rectangle (axis cs:2,1);
\draw[black!45,densely dotted] (axis cs:2,0) -- (axis cs:2,1);
\node[anchor=north west,font=\scriptsize,text=black!60] at (axis cs:2.03,0.98) {BBP};
\draw[oracleTheta,dashed] (axis cs:1,0.5) -- (axis cs:6.4,0.5);
\draw[oracleW,dashed] (axis cs:1,0.25) -- (axis cs:6.4,0.25);
\draw[oracleTheta,densely dotted] (axis cs:3,0) -- (axis cs:3,0.5);
\draw[oracleW,densely dotted] (axis cs:2.333333333,0) -- (axis cs:2.333333333,0.25);
\node[anchor=south east,font=\scriptsize,text=oracleTheta,fill=white,inner sep=2pt] at (axis cs:6.35,0.51) {Chordal: $\tau (1/2 ; 1) = 3$};
\node[anchor=south east,font=\scriptsize,text=oracleW,fill=white,inner sep=2pt] at (axis cs:6.35,0.26) {Bures--Wasserstein: $\tau (1/4 ; 1) = 7/3$};
\addplot[black!80,thick,domain=2:6.4,samples=160] {(x - 2)/(x - 1)};
\addplot[black!80,thick] coordinates {(1,0) (2,0)};
\addplot[only marks,mark=*,mark size=2.2pt,black!75] coordinates {(6,0.8) (4,0.666666667) (2.6,0.375) (2.2,0.166666667) (1.7,0)};
\node[above left,font=\scriptsize] at (axis cs:6,0.8) {$\lambda_{1}$};
\node[above left,font=\scriptsize] at (axis cs:4,0.666666667) {$\lambda_{2}$};
\node[above left,font=\scriptsize] at (axis cs:2.6,0.375) {$\lambda_{3}$};
\node[above left,font=\scriptsize] at (axis cs:2.2,0.166666667) {$\lambda_{4}$};
\node[above,font=\scriptsize] at (axis cs:1.7,0) {$\lambda_{5}$};
\addplot[only marks,mark=o,mark size=3pt,thick,oracleTheta] coordinates {(3,0.5)};
\addplot[only marks,mark=o,mark size=3pt,thick,oracleW] coordinates {(2.333333333,0.25)};

\nextgroupplot[
    title={(b) Population loss curves for the same spikes},
    xmin=0, xmax=7.2, ymin=2.4, ymax=8.5,
    xlabel={Candidate dimension $k$},
    ylabel={Limiting squared loss},
    xtick={0,1,2,3,4,5,6,7},
    ytick={3,4,5,6,7,8},
    legend style={at={(0.02,0.98)},anchor=north west,draw=none,fill=white,font=\small,cells={anchor=west}},
    legend image post style={mark size=2pt}
]
\draw[black!35,densely dotted] (axis cs:4,2.4) -- (axis cs:4,8.3);
\draw[black!55,dashed] (axis cs:5,2.4) -- (axis cs:5,8.3);
\node[font=\small,fill=white,inner sep=2pt] at (axis cs:4,8.15) {$d_{\gamma}^{*} = 4$};
\node[font=\small,fill=white,inner sep=2pt] at (axis cs:5,8.15) {$d^{*} = 5$};
\draw[oracleTheta,densely dotted] (axis cs:2,2.4) -- (axis cs:2,4.066666667);
\draw[oracleW,densely dotted] (axis cs:3,2.4) -- (axis cs:3,3.353407585);
\addplot[oracleTheta,thick,mark=*,mark size=2.2pt] table[x=k,y=Theta] {\OracleLossTable};
\addlegendentry{$L_{\Theta}^{\mathrm{or}} (k)$}
\addplot[oracleW,thick,mark=square*,mark size=2.2pt] table[x=k,y=Wasserstein] {\OracleLossTable};
\addlegendentry{$L_{\c{W}}^{\mathrm{or}} (k)$}
\addplot[only marks,oracleTheta,mark=o,mark size=4.3pt,thick,forget plot] coordinates {(2,4.066666667)};
\addplot[only marks,oracleW,mark=square,mark size=4.3pt,thick,forget plot] coordinates {(3,3.353407585)};
\node[above,font=\small,text=oracleTheta,fill=white,inner sep=1pt] at (axis cs:2,4.3) {$d_{\Theta}^{\mathrm{or}} = 2$};
\node[below,font=\small,text=oracleW,fill=white,inner sep=1pt] at (axis cs:3.1,3.2) {$d_{\c{W}}^{\mathrm{or}} = 3$};
\end{groupplot}
\end{tikzpicture}
}
\caption{\textbf{(a)}: limiting squared alignment, with the BBP threshold at $1 + \sqrt{\gamma} = 2$ and geometric inclusion cutoffs $\tau (1/4 ; 1) = 7/3$ and $\tau (1/2 ; 1) = 3$. Filled markers indicate the population spikes. \textbf{(b)}: we choose a pathological example $(\lambda_{1}, \ldots, \lambda_{5}) = (6, 4, 2.6, 2.2, 1.7)$ so that all the dimensionalities differ: $d_{\Theta}^{\mathrm{or}} = 2$ and $d_{\c{W}}^{\mathrm{or}} = 3$, although $d_{\gamma}^{*} = 4$ spikes are supercritical and $d^{*} = 5$ spikes are present.}
\label{fig:hd:population:benchmark}
\end{figure}

\begin{proposition}[Population benchmark]\label{thm:macro:bias}
Under (A1)--(A3), for every fixed $k \in \b{N}_{0}$,
\begin{align*}
    \widehat{L}_{\Theta}^{\mathrm{or}} (k) &\defeq l_{\Theta} (\c{U}_{[d^{*}]}, \hat{\c{V}}_{[k]})^{2} \xrightarrow{a.s.} L_{\Theta}^{\mathrm{or}} (k) \defeq d^{*} + k - 2 \sum_{j = 1}^{k \wedge d^{*}} g (\lambda_{j} ; \gamma), \\
    \widehat{L}_{\c{W}}^{\mathrm{or}} (k) &\defeq l_{\c{W}} (\c{U}_{[d^{*}]}, \hat{\c{V}}_{[k]})^{2} \xrightarrow{a.s.} L_{\c{W}}^{\mathrm{or}} (k) \defeq d^{*} + k - 2 \sum_{j = 1}^{k \wedge d^{*}} \sqrt{g (\lambda_{j} ; \gamma)}.
\end{align*}
Both convergences are uniform over any fixed finite collection of candidate dimensions.
\end{proposition}
Note that the limiting marginal costs for $1 \le k \le d^{*}$ are given by:
\begin{align}\label{eq:oracle:incre}
    L_{\Theta}^{\mathrm{or}} (k) - L_{\Theta}^{\mathrm{or}} (k - 1) = 1 - 2 g (\lambda_{k} ; \gamma), \qquad
    L_{\c{W}}^{\mathrm{or}} (k) - L_{\c{W}}^{\mathrm{or}} (k - 1) = 1 - 2 \sqrt{g (\lambda_{k} ; \gamma)},
\end{align}
with both increments equal to one for $k > d^{*}$. The marginal cost can be
resolved into two components. Adding a direction incurs a unit dimension
increase but receives an alignment reward for potentially getting closer to the true number
of supercritical spikes. Since the population eigenvalues are ordered and $g
(\lambda ; \gamma)$ is nondecreasing in $\lambda$, the increments are
nondecreasing in $k$. Both limiting loss curves are therefore discretely convex,
with linear increasing tails, see \cref{fig:hd:population:benchmark}. Their
respective minimizers can be identified as
\begin{align}\label{eq:hd:oracle:dimensions}
    d_{\Theta}^{\mathrm{or}} &\defeq \min \argmin_{k \in \b{N}_{0}} L_{\Theta}^{\mathrm{or}} (k) = \# \{ 1 \le j \le d^{*} : g (\lambda_{j} ; \gamma) > 1/2 \}, \\
    d_{\c{W}}^{\mathrm{or}} &\defeq \min \argmin_{k \in \b{N}_{0}} L_{\c{W}}^{\mathrm{or}} (k) = \# \{ 1 \le j \le d^{*} : g (\lambda_{j} ; \gamma) > 1/4 \}. \nonumber
\end{align}
To express alignment thresholds directly in terms of population eigenvalues, define
\begin{align*}
    \tau (t ; \gamma) \defeq 1 + \frac{t \gamma + \sqrt{t^{2} \gamma^{2} + 4 \gamma (1 - t)}}{2 (1 - t)}, \qquad \gamma > 0, \quad 0 < t < 1.
\end{align*}
Then $g (\lambda ; \gamma) > t$ if and only if $\lambda > \tau (t ; \gamma)$. Thus \eqref{eq:hd:oracle:dimensions} can equivalently be written as
\begin{align*}
    d_{\Theta}^{\mathrm{or}} = \# \{ 1 \le j \le d^{*} : \lambda_{j} > \tau (1/2 ; \gamma) \}, \quad
    d_{\c{W}}^{\mathrm{or}} = \# \{ 1 \le j \le d^{*} : \lambda_{j} > \tau (1/4 ; \gamma) \}.
\end{align*}
At equality with the corresponding cutoff, including or excluding the direction gives the same limiting loss.
Consequently, $d_{\Theta}^{\mathrm{or}} \le d_{\c{W}}^{\mathrm{or}} \le
d_{\gamma}^{*} \le d^{*}$. A supercritical spike has positive limiting
alignment, but that alignment can still be insufficient to offset the additional
unit dimension cost. Thus, minimizing subspace loss and counting supercritical
spikes present different statistical objectives. The role of the population
benchmark here is to identify this distinction and motivate an adjustment of the dimension penalty. As  $\c{U}_{[d^{*}]}$ is unknown, we next construct an observable criterion by replacing the population reference with an independent empirical subspace.

\subsection{Independent Empirical Subspaces}\label{ssec:data:split}

Empirical directions representing the same population signal can remain aligned
across independent samples, whereas the components that are orthogonal to the
population signal
subspace become asymptotically orthogonal. This agreement supplies an observable
measure of signal reproducibility. To leverage this fact, we partition the
observations independently of their values into samples $\m{X}_{1} \in
\b{R}^{n_{1} \times p}$ and $\m{X}_{2} \in \b{R}^{n_{2} \times p}$, where $n_{1}
+ n_{2} = n$ and $n_{1}/n \to \alpha \in (0, 1)$. Their respective limiting
aspect ratios are then defined as
\begin{align*}
    \gamma_{1} \defeq \frac{\gamma}{\alpha}, \qquad \gamma_{2} \defeq \frac{\gamma}{1 - \alpha}.
\end{align*}
We use the first sample to construct a reference subspace and the second to construct candidate subspaces. Choose a fixed reference dimension $D \ge \max \{ d^{*}, 1 \}$ and consider candidates $k \in \{ 0, \ldots, D \}$. The condition $D \ge d^{*}$ is an upper-bound assumption on the unknown signal dimension. The results below show that, subject to this condition, the precise value of a fixed $D$ does not affect the limiting minimizing dimension.

For $b \in \{ 1, 2 \}$, let $\hat{\m{v}}_{1}^{(b)}, \ldots, \hat{\m{v}}_{p}^{(b)}$ be orthonormal eigenvectors of $\hat{\m{\Sigma}}_{b} \defeq n_{b}^{-1} \m{X}_{b}^{\top} \m{X}_{b}$, ordered by decreasing eigenvalue, and write
\begin{align*}
    \hat{\m{V}}_{[k]}^{(b)} \defeq [\hat{\m{v}}_{1}^{(b)}, \ldots, \hat{\m{v}}_{k}^{(b)}], \qquad \hat{\c{V}}_{[k]}^{(b)} \defeq \spann \{ \hat{\m{v}}_{1}^{(b)}, \ldots, \hat{\m{v}}_{k}^{(b)} \}.
\end{align*}
We compare the reference $\hat{\c{V}}_{[D]}^{(1)}$ with each candidate $\hat{\c{V}}_{[k]}^{(2)}$ through
\begin{align}\label{eq:hd:split:criteria}
    \widehat{L}_{\Theta} (k ; \alpha) \defeq l_{\Theta} (\hat{\c{V}}_{[D]}^{(1)}, \hat{\c{V}}_{[k]}^{(2)})^{2}, \qquad
    \widehat{L}_{\c{W}} (k ; \alpha) \defeq l_{\c{W}} (\hat{\c{V}}_{[D]}^{(1)}, \hat{\c{V}}_{[k]}^{(2)})^{2}.
\end{align}
Keeping the reference dimension fixed allows candidates of different dimensions
to be compared against the same object and therefore create a rank ordering. If both subspaces instead had dimension $k$, their unmodified distance would be zero at $k = 0$, yielding a trivial minimum.
To analyze the limiting behavior of this loss, we define the joint squared alignment by
\begin{align}\label{eq:hd:joint:alignment}
    h (\lambda ; \gamma, \alpha) \defeq g (\lambda ; \gamma_{1}) \cdot g (\lambda ; \gamma_{2}) = g \left( \lambda ; \frac{\gamma}{\alpha} \right) g \left( \lambda ; \frac{\gamma}{1 - \alpha} \right).
\end{align}
This quantity is positive precisely when the spike is supercritical in both samples. It determines the limiting singular values of the cross matrix $(\hat{\m{V}}_{[D]}^{(1)})^{\top} \hat{\m{V}}_{[k]}^{(2)} \in \b{R}^{D \times k}$ between the empirical subspaces, including when population spikes are repeated.

\begin{theorem}[Distances between empirical subspaces]\label{thm:macro:bias:split}
Under (A1)--(A3), for each fixed $\alpha \in (0, 1)$ and every $k \in \{ 0, \ldots, D \}$,
\begin{align*}
    \widehat{L}_{\Theta} (k ; \alpha) &\xrightarrow{a.s.} L_{\Theta} (k ; \alpha) \defeq D + k - 2 \sum_{j = 1}^{k} h (\lambda_{j} ; \gamma, \alpha), \\
    \widehat{L}_{\c{W}} (k ; \alpha) &\xrightarrow{a.s.} L_{\c{W}} (k ; \alpha) \defeq D + k - 2 \sum_{j = 1}^{k} \sqrt{h (\lambda_{j} ; \gamma, \alpha)}.
\end{align*}
Both convergences are uniform over $k \in \{ 0, \ldots, D \}$.
\end{theorem}

The additional $D - d^{*}$ reference directions have no limiting alignment with
any fixed number of candidate directions. Once $D \ge d^{*}$, any change of the
fixed reference dimension, $D$, therefore shifts the limiting loss on the common candidate dimensions by an additive constant, while candidates beyond $d^{*}$ increase the unmodified loss by one per direction. Thus the limiting minimizing dimension is unchanged. The choice of $D$ can nevertheless affect finite-sample variability.
Compared to the population benchmark, splitting attenuates the alignment reward
because each subspace is estimated from fewer observations. In particular, $g
(\lambda ; \gamma)$ is nonincreasing in its aspect-ratio argument, and
$\gamma_{1}, \gamma_{2} > \gamma$, so the joint alignment in
\eqref{eq:hd:joint:alignment} satisfies $h (\lambda ; \gamma, \alpha) \le g
(\lambda ; \gamma)^{2} \le g (\lambda ; \gamma)$. In this regard, we must first
identify the split proportion that maximizes the limiting agreement between the
two empirical estimates.

\begin{proposition}[Optimal split proportion]\label{prop:half:ideal}
For every $\gamma > 0$, $\alpha \in (0, 1)$, and $\lambda \ge 1$,
\begin{align*}
    h (\lambda ; \gamma, \alpha) \le h (\lambda ; \gamma, 1/2) = g (\lambda ; 2 \gamma)^{2}.
\end{align*}
If $\lambda > 1 + \sqrt{2 \gamma}$, equality holds only at $\alpha = 1/2$.
\end{proposition}

Propositon~\ref{prop:half:ideal} demonstrates that equal splitting therefore maximizes every joint alignment and minimizes both limiting distance curves pointwise in the candidate dimension.
We henceforth take $n_{1} = \lfloor n/2 \rfloor$ and $n_{2} = n - n_{1}$, so both limiting aspect ratios equal $2 \gamma$. In this regard, we suppress $\alpha = 1/2$ by writing
\begin{align*}
    h (\lambda ; \gamma) \defeq h (\lambda ; \gamma, 1/2) = g (\lambda ; 2 \gamma)^{2}, \quad
    \widehat{L}_{\mathsf{G}} (k) \defeq \widehat{L}_{\mathsf{G}} (k ; 1/2), \quad L_{\mathsf{G}} (k) \defeq L_{\mathsf{G}} (k ; 1/2), \quad \mathsf{G} \in \{ \Theta, \c{W} \}.
\end{align*}

\subsubsection{Calibrating the Dimension Penalty}
Even with this optimal split, the alignment reward for including a direction
remains attenuated relative to the full-sample oracle. Minimizing this split
criterion can therefore exclude directions that would reduce the full-sample
subspace loss. To adjust this balance, we introduce a dimension penalty indexed by $\eta \in (-1, 1)$:
\begin{align*} 
    \widehat{L}_{\mathsf{G}, \eta} (k) &\defeq \widehat{L}_{\mathsf{G}} (k) - \eta k, \qquad \hat{d}_{\mathsf{G}, \eta} \defeq \min \argmin_{0 \le k \le D} \widehat{L}_{\mathsf{G}, \eta} (k), \qquad \mathsf{G} \in \{ \Theta, \c{W} \}.
\end{align*}
Here $\eta = 0$ gives the original distance, $\eta > 0$ reduces the dimension cost, and $\eta < 0$ increases it. As in \eqref{eq:hd:oracle:dimensions}, the smallest-minimizer convention resolves ties.
For either geometry $\mathsf{G} \in \{ \Theta, \c{W} \}$, the limiting criterion $L_{\mathsf{G}, \eta} (k) \defeq L_{\mathsf{G}} (k) - \eta k$ has increments
\begin{align}
    L_{\Theta, \eta} (k) - L_{\Theta, \eta} (k - 1) &= 1 - \eta - 2 g (\lambda_{k} ; 2 \gamma)^{2}, \label{eq:risk:curve:slope:chordal} \\
    L_{\c{W}, \eta} (k) - L_{\c{W}, \eta} (k - 1) &= 1 - \eta - 2 g (\lambda_{k} ; 2 \gamma). \label{eq:risk:curve:slope:wasserstein}
\end{align}
These increments are nondecreasing in $k$, so both limiting criteria remain
discretely convex. Beyond the spikes that are supercritical in both samples,
additional direction increments equal $1 - \eta > 0$. Increasing $\eta$ lowers the alignment required for retention while preserving an increasing limiting tail.
Hence, the natural follow-up question is how to choose $\eta$. The following result provides a general calibration: any prescribed population cutoff $c_{\gamma} > 1 + \sqrt{2 \gamma}$ above the half-sample detection threshold can be implemented through a suitable dimension penalty.

\begin{proposition}[Calibration to a prescribed cutoff]\label{cor:hd:selected:dimensions}
Assume (A1)--(A3), and fix a population cutoff $c_{\gamma} > 1 + \sqrt{2 \gamma}$. Define
\begin{align*}
    \eta_{\Theta} (c_{\gamma} ; \gamma) \defeq 1 - 2 g (c_{\gamma} ; 2 \gamma)^{2}, \qquad
    \eta_{\c{W}} (c_{\gamma} ; \gamma)
    \defeq 1 - 2 g (c_{\gamma} ; 2 \gamma).
\end{align*}
For $\mathsf{G} \in \{ \Theta, \c{W} \}$, let $\eta_{n} \in (-1, 1)$ satisfy $\eta_{n} \xrightarrow{a.s.} \eta_{\mathsf{G}} (c_{\gamma} ; \gamma)$.
If no population spike equals $c_{\gamma}$, then
\begin{align*}
    \hat{d}_{\mathsf{G}, \eta_{n}}
    \xrightarrow{a.s.}
    \# \{ 1 \le j \le d^{*} : \lambda_{j} > c_{\gamma} \}.
\end{align*}
\end{proposition}

This calibration identifies which population cutoff a penalty implements. If the cutoff rule $\gamma \mapsto c_{\gamma}$ is continuous at the limiting aspect ratio, the observed ratio $\gamma_{n} \defeq p/n$ can be substituted in \cref{cor:hd:selected:dimensions}. This uses an observed quantity and does not require estimating a population noise scale. A preferred penalty, however, still depends on the desired balance between retaining weak signals and selecting spurious directions.
The population benchmark provides one concrete choice when the objective is to minimize full-sample chordal loss. Its cutoff is
\begin{align}\label{eq:hd:chordal:cutoff}
    c_{\Theta} (\gamma) \defeq \tau (1/2 ; \gamma) = 1 + x_{\gamma}, \qquad x_{\gamma} \defeq \frac{\gamma + \sqrt{\gamma^{2} + 8 \gamma}}{2} > \sqrt{2 \gamma}.
\end{align}
In this case, we get
\begin{align}\label{eq:hd:optimal:penalty}
    \eta_{\Theta}^{\mathrm{or}} (\gamma) \defeq 1 - 2 g (c_{\Theta} (\gamma) ; 2 \gamma)^{2} = 1 - 2 \left( \frac{x_{\gamma}}{3 x_{\gamma} + 2} \right)^{2}.
\end{align}
If we then plug-in the observed ratio $\gamma_{n} \defeq p/n$, this gives rise to the calibrated selector by
\begin{align}\label{eq:hd:calibrated:selector}
    \hat{d}_{\Theta}^{\mathrm{cal}} \defeq \hat{d}_{\Theta, \eta_{\Theta}^{\mathrm{or}} (\gamma_{n})} = \min \argmin_{0 \le k \le D} \left\{ \widehat{L}_{\Theta} (k) - \eta_{\Theta}^{\mathrm{or}} (\gamma_{n}) k \right\}.
\end{align}

\begin{theorem}[Vanishing oracle regret for chordal loss]\label{thm:hd:oracle:calibration}
Under (A1)--(A3), the selector in \eqref{eq:hd:calibrated:selector} satisfies
\begin{align*} 
    l_{\Theta} (\c{U}_{[d^{*}]}, \hat{\c{V}}_{[\hat{d}_{\Theta}^{\mathrm{cal}}]})^{2} - \min_{0 \le k \le D} l_{\Theta} (\c{U}_{[d^{*}]}, \hat{\c{V}}_{[k]})^{2} \xrightarrow{a.s.} 0.
\end{align*}
If no population spike equals $c_{\Theta} (\gamma)$, then $\hat{d}_{\Theta}^{\mathrm{cal}} \xrightarrow{a.s.} d_{\Theta}^{\mathrm{or}}$.
\end{theorem}

This result shows that splitting need not incur asymptotic regret for the
full-sample chordal loss. Its target remains $d_{\Theta}^{\mathrm{or}}$, which
is smaller than or equal to $d_{\gamma}^{*}$. Exact recovery of the
supercritical spike count $d_{\gamma}^{*}$ requires a different signal-strength
condition, which is what we address next.

\subsection{Dimension Selection Guarantees}\label{ssec:phase:trans}

We now relate the minimizer of the observable criterion to the target
$d_{\gamma}^{*}$. The penalty, $\eta$, determines a positive alignment
threshold, while splitting determines which population spikes have positive
limiting agreement at all. Separating these two effects makes the conditions for
successful dimension selection explicit.

\begin{theorem}[Limiting selected dimension]\label{thm:hd:dimension:consistency}
Assume (A1)--(A3), and let $\eta_{n} \in (-1, 1)$ satisfy $\eta_{n} \xrightarrow{a.s.} \eta$ for a deterministic $\eta \in (-1, 1)$. Define
\begin{align*}
    d_{\Theta, \eta} &\defeq \# \{ 1 \le j \le d^{*} : g (\lambda_{j} ; 2 \gamma)^{2} > (1 - \eta)/2 \} = \# \{ 1 \le j \le d^{*} : \lambda_{j} > \tau (\sqrt{(1 - \eta)/2} ; 2 \gamma) \}, \\
     d_{\c{W}, \eta} &\defeq \# \{ 1 \le j \le d^{*} : g (\lambda_{j} ; 2 \gamma) > (1 - \eta)/2 \} = \# \{ 1 \le j \le d^{*} : \lambda_{j} > \tau ((1 - \eta)/2 ; 2 \gamma) \}.
\end{align*}
For chordal loss, suppose $g (\lambda_{j} ; 2 \gamma)^{2} \ne (1 - \eta)/2$ for every $1 \le j \le d^{*}$. For Bures--Wasserstein loss, suppose $g (\lambda_{j} ; 2 \gamma) \ne (1 - \eta)/2$ for every $1 \le j \le d^{*}$. Under the corresponding condition,
\begin{align*}
    \hat{d}_{\mathsf{G}, \eta_{n}} \xrightarrow{a.s.} d_{\mathsf{G}, \eta}, \qquad \mathsf{G} \in \{ \Theta, \c{W} \}.
\end{align*}
\end{theorem}

Repeated population spikes are counted according to their multiplicities. Since every direction in a repeated block has the same limiting alignment, the limiting criterion retains or excludes the entire block unless that alignment lies exactly at the selection threshold. The no-equality condition in \cref{thm:hd:dimension:consistency} excludes this latter case, in which the limiting objective alone does not determine the empirical minimizing dimension.
Every spike counted by $d_{\mathsf{G}, \eta}$ is supercritical in both half-samples and therefore also in the full sample in the following sense: $d_{\mathsf{G}, \eta} \le \# \{ 1 \le j \le d^{*} : \lambda_{j} > 1 + \sqrt{2 \gamma} \} \le d_{\gamma}^{*}$, which yields asymptotic control of over-selection:
\begin{align}\label{eq:hd:asymptotic:overselection}
    \b{P} \big( \hat{d}_{\mathsf{G}, \eta_{n}} > d_{\gamma}^{*} \big) \to 0.
\end{align}

\begin{corollary}[Consistent recovery of supercritical spikes]\label{cor:hd:rank:consistency}
Assume (A1)--(A3), and let $\eta_{n} \in (-1, 1)$ satisfy $\eta_{n} \xrightarrow{a.s.} \eta$ for a deterministic $\eta \in (-1, 1)$. Suppose either $d_{\gamma}^{*} = 0$ or, when $d_{\gamma}^{*} > 0$, the corresponding condition holds:
\begin{align*}
    \lambda_{d_{\gamma}^{*}} > \tau \left( \sqrt{\frac{1 - \eta}{2}} ; 2 \gamma \right) \qquad \text{for chordal loss}, \qquad
    \lambda_{d_{\gamma}^{*}} > \tau \left( \frac{1 - \eta}{2} ; 2 \gamma \right) \qquad \text{for B-W loss}.
\end{align*}
Then $\hat{d}_{\mathsf{G}, \eta_{n}} \xrightarrow{a.s.} d_{\gamma}^{*}$ and $\b{P} ( \hat{d}_{\mathsf{G}, \eta_{n}} \ne d_{\gamma}^{*} ) \to 0$ for the corresponding geometry. If all population spikes are supercritical, this gives consistent recovery of $d^{*}$.
\end{corollary}

For any fixed configuration whose target spikes all exceed $1 + \sqrt{2
  \gamma}$, their minimum joint alignment is positive. Therefore there exists an $\eta < 1$, sufficiently close to one, for which \cref{cor:hd:rank:consistency} applies.

\subsection{Type I-error Control under Gaussian Null}
Under the Gaussian pure null ($d^{*} = 0$), our criterion admits a calibration that controls the probability of selecting a positive dimension when no signal is present. We first derive a finite-sample guarantee from the exact marginal distribution of the projected lengths, and then obtain a sharper asymptotic calibration using their joint distribution. Throughout this subsubsection, we may assume $1 \le D \le \min \{p, n_{1}, n_{2} \}$, which will be trivially satisfied in the high-dimensional setup.
The squared Frobenius representation of chordal loss makes the relevant statistics explicit. Define
\begin{align*} 
    \hat{\m{\Pi}}_{[D]}^{(1)} \defeq \hat{\m{V}}_{[D]}^{(1)} (\hat{\m{V}}_{[D]}^{(1)})^{\top}, \qquad q_{j} \defeq \| \hat{\m{\Pi}}_{[D]}^{(1)} \hat{\m{v}}_{j}^{(2)} \|_{2}^{2}, \qquad 1 \le j \le D.
\end{align*}
Then, exactly in finite samples,
\begin{align}\label{eq:hd:cumulative:criterion}
    \widehat{L}_{\Theta, \eta} (k) = D + \sum_{j = 1}^{k} (1 - \eta - 2 q_{j}), \qquad 0 \le k \le D.
\end{align}
Consequently, $\hat{d}_{\Theta, \eta} \defeq \min \argmin_{0 \le k \le D} \widehat{L}_{\Theta, \eta} (k)$ is computed by evaluating the $D + 1$ cumulative sums and choosing the smallest minimizing index. The $q_{j}$ need not be ordered, so counting those exceeding $(1 - \eta)/2$ does not generally give the minimizer. More precisely, the smallest-minimizer convention gives the exact identity
\begin{align}\label{eq:hd:null:selection:event}
    \{ \hat{d}_{\Theta, \eta} > 0 \} = \left\{ \max_{1 \le k \le D} \frac{1}{k} \sum_{j = 1}^{k} q_{j} > \frac{1 - \eta}{2} \right\}.
\end{align}
Indeed, a positive dimension is selected precisely when at least one candidate has a strictly smaller criterion value than the value $D$ at $k = 0$.

\begin{proposition}[Finite-sample selection control]\label{prop:hd:null:calibration}
Suppose the observations are independent $\c{N}_{p} (\m{0}, \sigma^{2}
\m{I}_{p})$, where $\sigma^{2} > 0$ is arbitrary, and the sample split is independent of the observations. For every $1 \le j \le D$,
\begin{align*}
    q_{j} \sim \operatorname{Beta} \left( \frac{D}{2}, \frac{p - D}{2} \right).
\end{align*}
Let $F_{p,D}$ denote the cumulative distribution function (CDF) of this beta distribution. Then
\begin{align}\label{eq:hd:null:bound}
    \b{P}_{0} \big( \hat{d}_{\Theta, \eta} > 0 \big) \le \min \left\{ 1, D \left[ 1 - F_{p,D} \left( \frac{1 - \eta}{2} \right) \right] \right\}, \qquad \eta \in (-1, 1).
\end{align}
\end{proposition}

\cref{prop:hd:null:calibration} guarantees that, for a prescribed error level $\varepsilon \in (0, 1)$, we have $\b{P}_{0} \big( \hat{d}_{\Theta, \eta} > 0 \big) \le \varepsilon$ by letting $\eta \defeq 1 - 2 F_{p,D}^{-1} (1 - \varepsilon/D)$. In other words, the beta quantile provides finite-sample control of the type I error for $H_{0} : d^{*} = 0$, without estimating $\sigma^{2}$. However, this Bonferroni correction is extremely conservative, and a sharper asymptotic calibration follows by retaining the cumulative event in \eqref{eq:hd:null:selection:event}.

\begin{theorem}[Asymptotic calibration of type I-error]\label{thm:hd:null:joint}
Under the Gaussian null of \cref{prop:hd:null:calibration}, suppose that $D$ is fixed and $p \to \infty$. Then
\begin{align*}
    (p q_{1}, \ldots, p q_{D}) \xrightarrow{d} (Y_{1}, \ldots, Y_{D}), \qquad Y_{1}, \ldots, Y_{D} \stackrel{\mathrm{iid}}{\sim} \chi_{D}^{2}.
\end{align*}
Consequently,
\begin{align}\label{eq:hd:null:cumulative:limit}
    T_{p,D} \defeq \max_{1 \le k \le D} \frac{p}{k} \sum_{j = 1}^{k} q_{j} \xrightarrow{d} M_{D} \defeq \max_{1 \le k \le D} \frac{1}{k} \sum_{j = 1}^{k} Y_{j}.
\end{align}
Let $H_{D}$ denote the CDF of $M_{D}$. Then the penalty $\eta_{\Theta} \defeq 1 - 2 H_{D}^{-1} (1 - \varepsilon)/p$ belongs to $(-1, 1)$ for all sufficiently large $p$ and satisfies $\b{P}_{0} \big( \hat{d}_{\Theta, \eta_{\Theta}} > 0 \big) \to \varepsilon$.
\end{theorem}

\begin{table}[h!]
    \centering
    \caption{Quantiles $H_{D}^{-1} (1 - \varepsilon)$ for $\varepsilon \in \{ 0.05, 0.005 \}$ and different reference dimensions $D$.}
    \label{tab:hd:null:quantiles}
    \begin{tabular}{lrrrrrr}
        \hline
        $D$ & 1 & 5 & 10 & 20 & 50 & 100 \\
        \hline
        $H_{D}^{-1} (0.95)$ & 3.8415 & 11.3917 & 18.6606 & 31.8217 & 68.0439 & 125.0317 \\
        $H_{D}^{-1} (0.995)$ & 7.8794 & 16.7973 & 25.2364 & 40.0491 & 79.5538 & 140.2479 \\
        \hline
    \end{tabular}
\end{table}

\cref{thm:hd:null:joint} provides a penalty with asymptotically exact type I error under the Gaussian null. The same quantile also yields a calibration criterion for Bures--Wasserstein loss. Specifically, define
\begin{align}\label{eq:sim:gds:penalties}
    t_{p, D, \varepsilon} \defeq \frac{H_{D}^{-1} (1 - \varepsilon)}{p}, \qquad \eta_{\Theta} \defeq 1 - 2 t_{p, D, \varepsilon}, \qquad \eta_{\c{W}} \defeq 1 - 2 \sqrt{t_{p, D, \varepsilon}}.
\end{align}
For fixed $D$ and $\varepsilon$, we have $t_{p, D, \varepsilon} \to 0$ as $p \to \infty$. Thus both penalties belong to $(-1, 1)$ for all sufficiently large $p$. With these penalties,the penalized criteria satisfy
\begin{align*}
    \widehat{L}_{\Theta, \eta_{\Theta}} (k) - D = 2 k t_{p, D, \varepsilon} - 2 \vertii{\m{C}_{k}}_{2}^{2}, \qquad
    \widehat{L}_{\c{W}, \eta_{\c{W}}} (k) - D = 2 k \sqrt{t_{p, D, \varepsilon}} - 2 \vertii{\m{C}_{k}}_{1},
\end{align*}
where $\m{C}_{k} \defeq (\hat{\m{V}}_{[D]}^{(1)})^{\top} \hat{\m{V}}_{[k]}^{(2)} \in \b{R}^{D \times k}$. If the Bures--Wasserstein criterion is smaller than its value at zero for some $1 \le k \le D$, then the Cauchy--Schwarz inequality gives
\begin{align*}
    \sqrt{t_{p, D, \varepsilon}} < \frac{\vertii{\m{C}_{k}}_{1}}{k} \le \sqrt{\frac{\vertii{\m{C}_{k}}_{2}^{2}}{k}}.
\end{align*}
It follows that $\vertii{\m{C}_{k}}_{2}^{2} > k t_{p, D, \varepsilon}$, so the chordal criterion is also smaller than its value at zero. Consequently, the smallest-minimizer convention yields $\{ \hat{d}_{\c{W}, \eta_{\c{W}}} > 0 \} \subseteq \{ \hat{d}_{\Theta, \eta_{\Theta}} > 0 \}$. As a result, under the assumptions of \cref{thm:hd:null:joint}, we get
\begin{align*}
    \limsup_{p \to \infty} \b{P}_{0} \big( \hat{d}_{\c{W}, \eta_{\c{W}}} > 0 \big) \le \lim_{p \to \infty} \b{P}_{0} \big( \hat{d}_{\Theta, \eta_{\Theta}} > 0 \big) = \varepsilon,
\end{align*}
i.e., the Bures--Wasserstein calibration is therefore asymptotically conservative.

\section{Data Analysis}
\label{sec:simulations}

\subsection{Geometric Dimension Selection}\label{ssec:sim:gds}

We now present the algorithmic implementation of the regularized risk functions
studied above as the \emph{Geometric Dimension Selection} (GDS) procedure, with
GDS-$\Theta$ and GDS-$\c{W}$ denoting its chordal and Bures--Wasserstein
versions, respectively.
The user specifies a reference dimension $D$, a geometry $\mathsf{G} \in \{
\Theta, \c{W} \}$, and a nominal error level $\varepsilon$; the penalty $\eta_{\mathsf{G}}$ is then determined by \eqref{eq:sim:gds:penalties}.

Due to the maximal limiting agreement established in \cref{prop:half:ideal}, we
partition the observations independently of their values into two samples of
sizes $n_{1} = \lfloor n/2 \rfloor$ and $n_{2} = n - n_{1}$. We necessitate that
$D < p$ be chosen no larger than the available rank in either sample, and
require $t_{p, D, \varepsilon} < 1$ so that the calibrated penalties, $\eta_{\mathsf{G}}$, belong to $(-1, 1)$.
The nominal error level $\varepsilon$ in \cref{thm:hd:null:joint} specifies the asymptotic probability of positive selection under the Gaussian pure null: the chordal calibration attains this level, while the Bures--Wasserstein calibration is conservative. The resulting procedure is summarized in \cref{alg:gds}.

\begin{algorithm}[h!]
\caption{Geometric Dimension Selection (GDS-$\mathsf{G}$)}\label{alg:gds}
\begin{algorithmic}[1]
    \Require Data matrix $\m{X} \in \b{R}^{n \times p}$, reference dimension $D$, and nominal level $\varepsilon$.
    \Ensure Selected dimensions $\hat{d}_{\Theta, \eta_{\Theta}}$ and $\hat{d}_{\c{W}, \eta_{\c{W}}}$, their criterion curves, and the corresponding full-sample subspace estimates.
    \State Randomly partition the observations into two equal samples $\m{X}_{1}$ and $\m{X}_{2}$ independently of their values.
    \State Compute the leading $D$ right singular vectors $\hat{\m{V}}_{[D]}^{(1)}$ and $\hat{\m{V}}_{[D]}^{(2)}$ of the two data matrices.
    \State Form $\m{C}_{D} \gets (\hat{\m{V}}_{[D]}^{(1)})^{\top} \hat{\m{V}}_{[D]}^{(2)}$.
    \State Compute $\eta_{\Theta}$ and $\eta_{\c{W}}$ using \eqref{eq:sim:gds:penalties}.
    \State Set $\widehat{L}_{\Theta, \eta_{\Theta}} (0) \gets D$ and $\widehat{L}_{\c{W}, \eta_{\c{W}}} (0) \gets D$.
    \For{$k = 1, \ldots, D$}
        \State Let $\m{C}_{k}$ contain the first $k$ columns of $\m{C}_{D}$.
        \State $\widehat{L}_{\Theta, \eta_{\Theta}} (k) \gets D + (1 - \eta_{\Theta}) k - 2 \vertii{\m{C}_{k}}_{2}^{2}$.
        \State $\widehat{L}_{\c{W}, \eta_{\c{W}}} (k) \gets D + (1 - \eta_{\c{W}}) k - 2 \vertii{\m{C}_{k}}_{1}$.
    \EndFor
    \For{$\mathsf{G} \in \{ \Theta, \c{W} \}$}
        \State Select $\hat{d}_{\mathsf{G}, \eta_{\mathsf{G}}} \gets \min \argmin_{0 \le k \le D} \widehat{L}_{\mathsf{G}, \eta_{\mathsf{G}}} (k)$.
        \State Fit the leading $\hat{d}_{\mathsf{G}, \eta_{\mathsf{G}}}$ principal components using all observations.
    \EndFor
\end{algorithmic}
\end{algorithm}

Beyond the selected dimension, the criterion curve is useful on its own: it
displays the balance between additional dimensions and agreement across the two
samples. A pronounced minimum indicates a clear preference within the observed
criterion, whereas a relatively flat region around the minimizer indicates that several dimensions
receive similar scores. In particular, reduced agreement can arise from
high-dimensional estimation error or weak signals, so the shape of the curve
gives an idea of how much the data deviates from classical regime where $p$ is
fixed. The following simulations examine how the curves and selected dimensions
change with the aspect ratio, spike configuration, and noise distribution.

Section~\ref{ssec:sim:study} presents a detailed comparative study of our
proposed methods (implemented via Algorithm~\ref{alg:gds}) against three existing
methods in the literature. We study different data configurations
allowing us to provide a detailed overview of the strengths and drawbacks of
each method and display the effectiveness of the theory presented in Section~\ref{sec:high:dim}.
Section~\ref{ssec:real_data_study} contains results of the methods on a real
lung cancer dataset. In particular, we present results on the full
high-dimensional dataset that show the effectiveness and utility of our proposed
methods over existing tools in the literature.

All implementations and replication scripts for the results throughout the remainder of this section are provided at \url{https://github.com/HoYUN-stat/GeometricPCA.jl}.

\subsection{Simulation Study}\label{ssec:sim:study}

We compare dimension selection methods in four configurations (\cref{tab:sim:design}). The first holds strong spikes fixed while varying the size of rare noise entries; the next two examine what happens as the spike strengths decrease and the cost of splitting becomes more pronounced. The fourth removes the signal altogether. We fix $n = 1000$ and vary $p \in \{ 250, 500, 1000, 2000 \}$, so that $\gamma = p/n \in \{ 0.25, 0.5, 1, 2 \}$. In (C1)--(C3), all four spikes exceed the full-sample transition $1 + \sqrt{\gamma}$ throughout this range. The target dimension is therefore $d^{*} = 4$ in (C1)--(C3) and $d^{*} = 0$ in (C4).

\begin{table}[htbp]
    \centering
    \small
    \caption{Simulation configurations. All unlisted population eigenvalues equal one.}
    \label{tab:sim:design}
    \begin{tabular}{lll}
        \hline
        Configuration & Population spikes & Noise distributions \\
        \hline
        (C1) Strong spikes & $40, 30, 25, 20$ & Gaussian; mixtures with $a \in \{ 10, 20, 30 \}$ \\
        (C2) Moderately strong spikes & $12, 9, 7, 5$ & Gaussian; $t_{5}$; $t_{8}$ \\
        (C3) Weak, repeated spikes & $3.6, 3.6, 3.2, 3.2$ & Gaussian; $t_{5}$ \\
        (C4) Pure null & None & All six distributions above \\
        \hline
    \end{tabular}
\end{table}

As described in the beginning of \cref{sec:high:dim}, we generate $\m{X} = \m{Z} \m{\Sigma}^{1/2}$ with independent, identically distributed entries in $\m{Z}$. Gaussian entries are standard normal, and Student entries are distributed as $\sqrt{(\nu - 2)/\nu} \, t_{\nu}$. For mixture noise, each entry is an independent draw from
\begin{align*}
    \frac{A G}{\sqrt{1 + \rho (a^{2} - 1)}}, \qquad G \sim \c{N} (0, 1), \qquad \b{P} (A = a) = \rho = 10^{-4}, \qquad \b{P} (A = 1) = 1 - \rho,
\end{align*}
where $A$ and $G$ are independent. Hence, increasing $a$ amplifies rare entries while preserving unit variance. The population mean is known to be zero; we neither center the data nor standardize the columns.

Both GDS procedures use two independent halves of size $n_{1} = n_{2} = 500$, reference dimension $D = 20$, and nominal Gaussian-null level $\varepsilon = 0.005$. We use the same tuning choices throughout; \cref{ssec:sim:sensitivity} examines their effects. We compare GDS with three full-sample methods. Since the residual covariance is $\m{I}_{p}$ in every configuration, we use BEMA0, the common-variance version of bulk eigenvalue matching, at the same nominal level \cite{ke2023estimation}. Gavish--Donoho hard thresholding (GD) uses a cutoff derived for matrix reconstruction, with the noise scale estimated from the Marchenko--Pastur median \cite{gavish2014optimal}. We use deflated deterministic parallel analysis (DDPA) without inflating its estimated noise edge \cite{dobriban2019deterministic}. Each row summarizes $1000$ independent replicates. The mean selected rank indicates the direction of bias; the percentage selecting the true rank measures exact recovery. Boldface marks the highest observed recovery rate in each row, including ties.

\paragraph{(C1) Strong spikes and rare large noise entries.}
Here the four spikes are strong enough. Under Gaussian noise, both GDS procedures, BEMA0, and GD recovered all four in every replicate (\cref{tab:sim:c1}). However, for a mixture noise,
increasing $a$ has the effect that a few entries of $Z_{ij}$ take on extremal values, and even one such entry can inflate the empirical variance along a noise direction enough to make it appear to be an additional spike.
At $a = 20$ and $\gamma = 1$, BEMA0 selected $9.728$ directions on average, compared with $4.039$ for GDS-$\c{W}$, against a true rank of four. The difficulty here was distinguishing the four strong spikes from additional large sample eigenvalues. Splitting helps because a direction induced by an isolated entry in one half earns little alignment reward unless it is also represented in the leading subspace estimated from the other half. The strong signal directions, by contrast, tend to be represented in both halves.

The choice of distance mattered most under the strongest mixture. At $a = 30$ and $\gamma = 2$, GDS-$\Theta$ selected $6.177$ directions on average, compared with $4.317$ for GDS-$\c{W}$. Thus the main advantage of the Bures--Wasserstein criterion in this setting was a reduction in over-selection. DDPA that re-estimates the noise edge after deflation, was competitive under the stronger mixtures: GDS-$\c{W}$ led in five of the six displayed rows with $a \ge 20$, and DDPA led in the sixth. Finally, GD recovered the rank under the mild mixture ($a = 10$), but behaves poorly when $a \ge 20$.

\begin{table}[htbp]
    \centering
    \small
    \setlength{\tabcolsep}{3pt}
    \caption{(C1) Strong spikes under Gaussian and mixture noise. Each entry is the mean selected rank ($1000$ replicates), followed by the percentage selecting $d^{*} = 4$. Boldface marks the highest exact-recovery rate in each row, including ties.}
    \label{tab:sim:c1}
    \resizebox{\linewidth}{!}{%
    \begin{tabular}{lrrrrrr}
        \hline
        Noise & $\gamma$ & GDS-$\Theta$ & GDS-$\c{W}$ & BEMA0 & GD & DDPA \\
        \hline
        Gaussian & $0.5$ & \textbf{4.000 (100.0\%)} & \textbf{4.000 (100.0\%)} & \textbf{4.000 (100.0\%)} & \textbf{4.000 (100.0\%)} & 4.286 (73.7\%) \\
         & $1$ & \textbf{4.000 (100.0\%)} & \textbf{4.000 (100.0\%)} & \textbf{4.000 (100.0\%)} & \textbf{4.000 (100.0\%)} & 4.234 (78.1\%) \\
         & $2$ & \textbf{4.000 (100.0\%)} & \textbf{4.000 (100.0\%)} & \textbf{4.000 (100.0\%)} & \textbf{4.000 (100.0\%)} & 4.240 (77.3\%) \\
        \hline
        Mixture, $a = 10$ & $0.5$ & 4.001 (99.9\%) & 4.001 (99.9\%) & 4.072 (92.9\%) & \textbf{4.000 (100.0\%)} & 4.170 (84.6\%) \\
         & $1$ & 4.002 (99.8\%) & 4.002 (99.8\%) & 4.036 (96.4\%) & \textbf{4.000 (100.0\%)} & 4.221 (79.4\%) \\
         & $2$ & 4.001 (99.9\%) & \textbf{4.000 (100.0\%)} & 4.007 (99.3\%) & \textbf{4.000 (100.0\%)} & 4.227 (78.5\%) \\
        \hline
        Mixture, $a = 20$ & $0.5$ & 4.145 (93.4\%) & \textbf{4.077 (95.8\%)} & 8.699 (0.2\%) & 4.843 (41.2\%) & 4.089 (92.7\%) \\
         & $1$ & 4.143 (92.7\%) & \textbf{4.039 (97.7\%)} & 9.728 (0.0\%) & 4.557 (56.0\%) & 4.079 (93.1\%) \\
         & $2$ & 4.197 (89.4\%) & \textbf{4.045 (97.7\%)} & 10.139 (0.0\%) & 4.192 (82.6\%) & 4.147 (87.8\%) \\
        \hline
        Mixture, $a = 30$ & $0.5$ & 4.666 (70.9\%) & \textbf{4.183 (91.4\%)} & 15.741 (0.0\%) & 9.208 (0.6\%) & 4.312 (78.9\%) \\
         & $1$ & 5.338 (51.1\%) & \textbf{4.264 (86.3\%)} & 20.846 (0.0\%) & 10.116 (0.2\%) & 4.194 (84.9\%) \\
         & $2$ & 6.177 (29.5\%) & 4.317 (81.6\%) & 25.957 (0.0\%) & 9.170 (0.0\%) & \textbf{4.214 (82.5\%)} \\
        \hline
    \end{tabular}}
\end{table}

\paragraph{(C2) Moderately strong spikes.}
Reducing the weakest spike $\lambda_{4}$ from $20$ to $5$ leaves it above the half-sample transition even at $\gamma = 2$. Under Gaussian noise, GDS-$\c{W}$ recovered rank four in $99.7$--$99.9\%$ of replicates, compared with $100\%$ for BEMA0 (\cref{tab:sim:c2}). The ordering depended on the noise distribution: under $t_{5}$ noise, GD had the highest recovery rate at $\gamma = 0.5$ and $1$, but at $\gamma = 2$ its mean selected rank fell to $3.044$ and its exact-recovery rate to $4.2\%$. GDS-$\c{W}$ remained close to four at the same aspect ratio, averaging $4.037$ directions with $96.5\%$ exact recovery. Under the lighter-tailed $t_{8}$ noise, BEMA0 recovered the rank in every replicate, and GDS-$\c{W}$ did so in at least $99.9\%$.

\begin{table}[htbp]
    \centering
    \small
    \setlength{\tabcolsep}{3pt}
    \caption{(C2) Moderately strong spikes. Each entry is the mean selected rank ($1000$ replicates), followed by the percentage selecting $d^{*} = 4$. Boldface marks the highest exact-recovery rate in each row, including ties.}
    \label{tab:sim:c2}
    \resizebox{\linewidth}{!}{%
    \begin{tabular}{lrrrrrr}
        \hline
        Noise & $\gamma$ & GDS-$\Theta$ & GDS-$\c{W}$ & BEMA0 & GD & DDPA \\
        \hline
        Gaussian & $0.5$ & 4.001 (99.9\%) & 4.001 (99.9\%) & \textbf{4.000 (100.0\%)} & \textbf{4.000 (100.0\%)} & 4.272 (75.0\%) \\
         & $1$ & 4.001 (99.9\%) & 4.001 (99.9\%) & \textbf{4.000 (100.0\%)} & \textbf{4.000 (100.0\%)} & 4.204 (80.4\%) \\
         & $2$ & 4.028 (97.4\%) & 4.003 (99.7\%) & \textbf{4.000 (100.0\%)} & 3.032 (3.2\%) & 4.208 (80.0\%) \\
        \hline
        $t_{5}$ & $0.5$ & 4.026 (97.4\%) & 4.011 (98.9\%) & 4.104 (90.1\%) & \textbf{4.007 (99.3\%)} & 4.230 (79.3\%) \\
         & $1$ & 4.042 (95.9\%) & 4.026 (97.4\%) & 4.112 (89.1\%) & \textbf{4.008 (99.2\%)} & 4.236 (78.3\%) \\
         & $2$ & 4.079 (92.7\%) & \textbf{4.037 (96.5\%)} & 4.074 (93.0\%) & 3.044 (4.2\%) & 4.218 (79.0\%) \\
        \hline
        $t_{8}$ & $0.5$ & 4.003 (99.7\%) & 4.001 (99.9\%) & \textbf{4.000 (100.0\%)} & \textbf{4.000 (100.0\%)} & 4.252 (75.8\%) \\
         & $1$ & 4.005 (99.5\%) & \textbf{4.000 (100.0\%)} & \textbf{4.000 (100.0\%)} & \textbf{4.000 (100.0\%)} & 4.200 (81.2\%) \\
         & $2$ & 4.019 (98.2\%) & 4.001 (99.9\%) & \textbf{4.000 (100.0\%)} & 3.038 (3.8\%) & 4.221 (78.8\%) \\
        \hline
    \end{tabular}}
\end{table}

\paragraph{(C3) Weak, repeated spikes.}
With spikes $(3.6, 3.6, 3.2, 3.2)$, the cost of splitting becomes visible when the weakest spikes approach the half-sample transition (\cref{tab:sim:c3}). At $\gamma = 1$, the limiting squared agreement associated with the weakest spike block, $g (\lambda_{4} ; 2 \gamma)^{2}$, falls from $0.340$ at $\lambda_{4} = 5$ in (C2) to $0.094$ at $\lambda_{4} = 3.2$ in (C3). Thus the block remains detectable, but agreement between halves distinguishes it less clearly from chance overlap. Under Gaussian noise, BEMA0 recovered rank four in every replicate using the full sample. GDS-$\c{W}$ selected $4.025$ directions on average, yet erred in $88$ of the $1000$ replicates: $34$ ranks were too small and $54$ too large.
Additionally, we remark the poor behavior of GD under weak spikes. Even at $\gamma = 1$, its singular-value cutoff of $4/\sqrt{3}$ for $\m{X}/\sqrt{n}$ corresponds to a sample covariance eigenvalue cutoff of $16/3$ \cite{gavish2014optimal}. Since the limiting sample eigenvalue of a population spike $\lambda > 2$ is $\lambda^{2}/(\lambda - 1)$, the corresponding population cutoff is $\lambda = 4$, above even the largest spike in (C3), $\lambda_{1} = 3.6$.

\begin{table}[htbp]
    \centering
    \small
    \setlength{\tabcolsep}{3pt}
    \caption{(C3) Weak, repeated spikes. Each entry is the mean selected rank ($1000$ replicates), followed by the percentage selecting $d^{*} = 4$. Boldface marks the highest exact-recovery rate in each row, including ties.}
    \label{tab:sim:c3}
    \resizebox{\linewidth}{!}{%
    \begin{tabular}{lrrrrrr}
        \hline
        Noise & $\gamma$ & GDS-$\Theta$ & GDS-$\c{W}$ & BEMA0 & GD & DDPA \\
        \hline
        Gaussian & $0.25$ & \textbf{4.000 (100.0\%)} & \textbf{4.000 (100.0\%)} & \textbf{4.000 (100.0\%)} & 3.994 (99.4\%) & 4.347 (68.7\%) \\
         & $0.5$ & 4.009 (99.2\%) & \textbf{4.000 (100.0\%)} & \textbf{4.000 (100.0\%)} & 2.658 (4.9\%) & 4.242 (77.1\%) \\
         & $1$ & 4.342 (74.4\%) & 4.025 (91.2\%) & \textbf{4.000 (100.0\%)} & 0.078 (0.0\%) & 4.182 (82.6\%) \\
        \hline
        $t_{5}$ & $0.25$ & 4.031 (96.9\%) & 4.018 (98.2\%) & 4.088 (91.7\%) & \textbf{3.997 (98.7\%)} & 4.399 (67.7\%) \\
         & $0.5$ & 4.101 (91.1\%) & \textbf{4.056 (94.7\%)} & 4.100 (90.5\%) & 2.792 (7.2\%) & 4.242 (74.3\%) \\
         & $1$ & 4.576 (60.1\%) & 4.135 (82.8\%) & \textbf{4.075 (92.5\%)} & 0.127 (0.0\%) & 4.157 (80.9\%) \\
        \hline
    \end{tabular}}
\end{table}

\paragraph{(C4) Pure noise.}
Removing all four spikes asks whether a method can reject the noise directions that competed with the signal in (C1). At $\gamma = 1$, GDS-$\c{W}$ selected a positive rank in $0.6\%$ of Gaussian replicates, close to its nominal $0.5\%$, but in $8.2\%$ under the $a = 30$ mixture (\cref{tab:sim:c4}). The covariance is $\m{I}_{p}$ in both cases; only the Gaussian model guarantees the rotational invariance used for the null calibration. BEMA0 and GD selected positive ranks in every replicate of the latter mixture, averaging $17.790$ and $6.571$ directions even without signal. DDPA had the lowest false-positive rate in that row ($6.3\%$), despite selecting too often under Gaussian noise ($10.4\%$). Thus the advantage of split agreement in (C1) is relative: it rejects many noise directions that pass a spectral threshold, but it does not preserve the nominal Gaussian error rate under severe mixtures. A more comprehensive results for the other aspect ratios $\gamma = p/n \in \{ 0.25, 0.5, 1, 2 \}$ are given in \cref{tab:sim:c4:append}.

\begin{table}[htbp]
    \centering
    \small
    \setlength{\tabcolsep}{3pt}
    \caption{(C4) Pure noise at $\gamma = 1$. Each entry gives the mean selected rank followed by the percentage selecting zero. Its complement is the false-positive rate. Boldface marks the highest percentage selecting zero in each row, including ties.}
    \label{tab:sim:c4}
    \resizebox{\linewidth}{!}{%
    \begin{tabular}{lrrrrr}
        \hline
        Noise & GDS-$\Theta$ & GDS-$\c{W}$ & BEMA0 & GD & DDPA \\
        \hline
        Gaussian & 0.006 (99.4\%) & 0.006 (99.4\%) & 0.004 (99.6\%) & \textbf{0.000 (100.0\%)} & 0.108 (89.6\%) \\
        $t_{5}$ & \textbf{0.005 (99.5\%)} & \textbf{0.005 (99.5\%)} & 0.150 (85.6\%) & 0.008 (99.2\%) & 0.100 (90.3\%) \\
        $t_{8}$ & 0.004 (99.6\%) & 0.004 (99.6\%) & 0.006 (99.4\%) & \textbf{0.000 (100.0\%)} & 0.102 (90.2\%) \\
        Mixture, $a = 10$ & 0.011 (99.0\%) & 0.011 (99.0\%) & 0.070 (93.3\%) & \textbf{0.000 (100.0\%)} & 0.094 (91.2\%) \\
        Mixture, $a = 20$ & 0.154 (94.9\%) & \textbf{0.073 (96.2\%)} & 6.472 (0.0\%) & 0.684 (49.4\%) & 0.060 (94.3\%) \\
        Mixture, $a = 30$ & 0.722 (85.9\%) & 0.213 (91.8\%) & 17.790 (0.0\%) & 6.571 (0.0\%) & \textbf{0.065 (93.7\%)} \\
        \hline
    \end{tabular}}
\end{table}

\subsubsection{Sensitivity to the nominal level and reference dimension}\label{ssec:sim:sensitivity}

The choice of $\varepsilon$ matters differently for rank recovery and for false positives. At $D = 20$, raising $\varepsilon$ lowers the threshold $H_{D}^{-1} (1 - \varepsilon)/p$ and makes a positive rank easier to select. Across the $36$ non-null settings, raising $\varepsilon$ from $0.001$ to $0.01$ barely changed mean exact recovery for GDS-$\c{W}$ ($93.16\%$ to $93.05\%$), but its Gaussian-null selection rate rose from $0.06\%$ to $0.82\%$ (\cref{tab:sim:tuning}, Panel~A). Raising it further to $0.05$ increased mean rank and reduced exact recovery for both geometries, especially GDS-$\Theta$. These results support choosing $\varepsilon$ by the acceptable Gaussian false-positive rate: there was little average rank-recovery gain to trade against that choice near $0.005$.
Table~\ref{tab:sim:tuning}, Panel~B addresses the sensitivity regarding the choice of the reference dimension $D$. A larger reference subspace captures more of an unrelated direction by chance: at $p = 1000$, increasing $D$ from $5$ to $40$ raises the mean null projection $D/p$ eightfold. Recalibration offsets this increase for the first candidate $\b{P}_{0} (q_{1} > t_{p, D, \varepsilon})$ only from $0.473\%$ to $0.402\%$.

\begin{table}[htbp]
    \centering
    \small
    \caption{Sensitivity to $\varepsilon$ and $D$. Panel~A gives mean selected rank (exact recovery in percent) averaged equally over all 36 non-null configurations; boldface marks the higher exact-recovery rate. Gaussian-null positive selection averages the four $p$ values equally in the separate Haar experiment ($10000$ draws per $p$). Panel~B gives deterministic calibration values at $p = 1000$ and $\varepsilon = 0.005$; its last column is the exact single-direction probability in percent in \cref{prop:hd:null:calibration}, not the global selection rate.}
    \label{tab:sim:tuning}
    \begin{tabular}{lrrrr}
        \hline
        \multicolumn{5}{l}{\textit{A. Varying $\varepsilon$ at $D = 20$: observed performance}} \\
        & \multicolumn{2}{c}{Non-null: mean rank (correct \%)} & \multicolumn{2}{c}{Gaussian-null selection (\%)} \\
        $\varepsilon$ & GDS-$\Theta$ & GDS-$\c{W}$ & GDS-$\Theta$ & GDS-$\c{W}$ \\
        \hline
        $0.001$ & 4.105 (89.30\%) & \textbf{3.971 (93.16\%)} & 0.06 & 0.06 \\
        $0.005$ & 4.182 (87.91\%) & \textbf{3.998 (93.13\%)} & 0.39 & 0.39 \\
        $0.01$ & 4.227 (86.89\%) & \textbf{4.014 (93.05\%)} & 0.83 & 0.82 \\
        $0.025$ & 4.323 (84.64\%) & \textbf{4.044 (92.59\%)} & 2.13 & 2.09 \\
        $0.05$ & 4.437 (81.93\%) & \textbf{4.074 (91.63\%)} & 4.62 & 4.48 \\
        $0.1$ & 4.653 (77.64\%) & \textbf{4.123 (89.68\%)} & 9.71 & 9.17 \\
        $0.2$ & 5.070 (70.75\%) & \textbf{4.205 (86.07\%)} & 19.61 & 17.78 \\
        \hline
    \end{tabular}
    \par\medskip
    \begin{tabular}{rrrrrr}
        \hline
        \multicolumn{6}{l}{\textit{B. Varying $D$ at $p = 1000$ and $\varepsilon = 0.005$: calibration}} \\
        $D$ & $D/p$ & $t_{p, D, \varepsilon}$ & $\eta_{\Theta}$ & $\eta_{\c{W}}$ & $100 \, \b{P}_{0} (q_{1} > t_{p, D, \varepsilon})$ \\
        \hline
        5 & 0.005 & 0.01680 & 0.9664 & 0.7408 & 0.473 \\
        10 & 0.010 & 0.02524 & 0.9495 & 0.6823 & 0.463 \\
        20 & 0.020 & 0.04005 & 0.9199 & 0.5998 & 0.442 \\
        40 & 0.040 & 0.06683 & 0.8663 & 0.4830 & 0.402 \\
        \hline
    \end{tabular}
\end{table}

\subsection{Real Data Analysis: Lung Cancer}
\label{ssec:real_data_study}
The lung cancer dataset, collected by \cite{gordon2002translation} expression measurements for $p = 12{,}533$ genes from $n = 181$ subjects and has been studied using several dimension-reduction methods
\cite[amongst others]{ke2023estimation, jin2016influential, li2008selecting}.
The full dataset has a $\gamma \approx 69$, which constitutes a very
high-dimensional dataset. In fact, it can be prohibitively large for some PCA
methods. As a result, such datasets are typically pruned first. In this case,
the most commonly accepted pruning is done by removing genes that are
not \textit{differentially expressed} \cite{jin2016influential} based on Kolmogorov--Smirnov scores. This means the
genes that have limited variability (up to some pre-determined threshold) within
certain subgroups are
removed from consideration. Following the proposed pruning method of marginal
feature screening, we obtain a data matrix of $(n, p) = (181, 251)$ which is the
usual dataset used for analysis.
This give a significantly smaller $\gamma \approx 1.38$. For the pruned dataset,
it has been established the the number of spikes that should be estimated is
exactly equal to 1 \cite{ke2023estimation}. The full dataset, however, has not
been analysed previously and thus the \textit{correct} number of spikes to be
estimated on the larger feature space remains an open question, in part due to
the ineaxct pruning method that may omit some smaller signal spikes.

It remains unclear whether the screening method
proposed by \cite{jin2016influential} risks discarding critical biological
information by evaluating each gene strictly in isolation. It is possible that
secondary features of importance are being truncated and therefore weaker
background signals in the high-dimensional dataset are being lost before the PCA
approach can even attempt to detect it \cite{nadler2016discussion, fan2008sure, bourgon2010independent}. Thus, it is not directly evident that the
number of spikes estimated from the full dataset should match that of the pruned
dataset.

\begin{figure}[h!]
  \centering
  \begin{subfigure}[b]{0.45\textwidth}
    \centering
    \includegraphics[width=\textwidth]{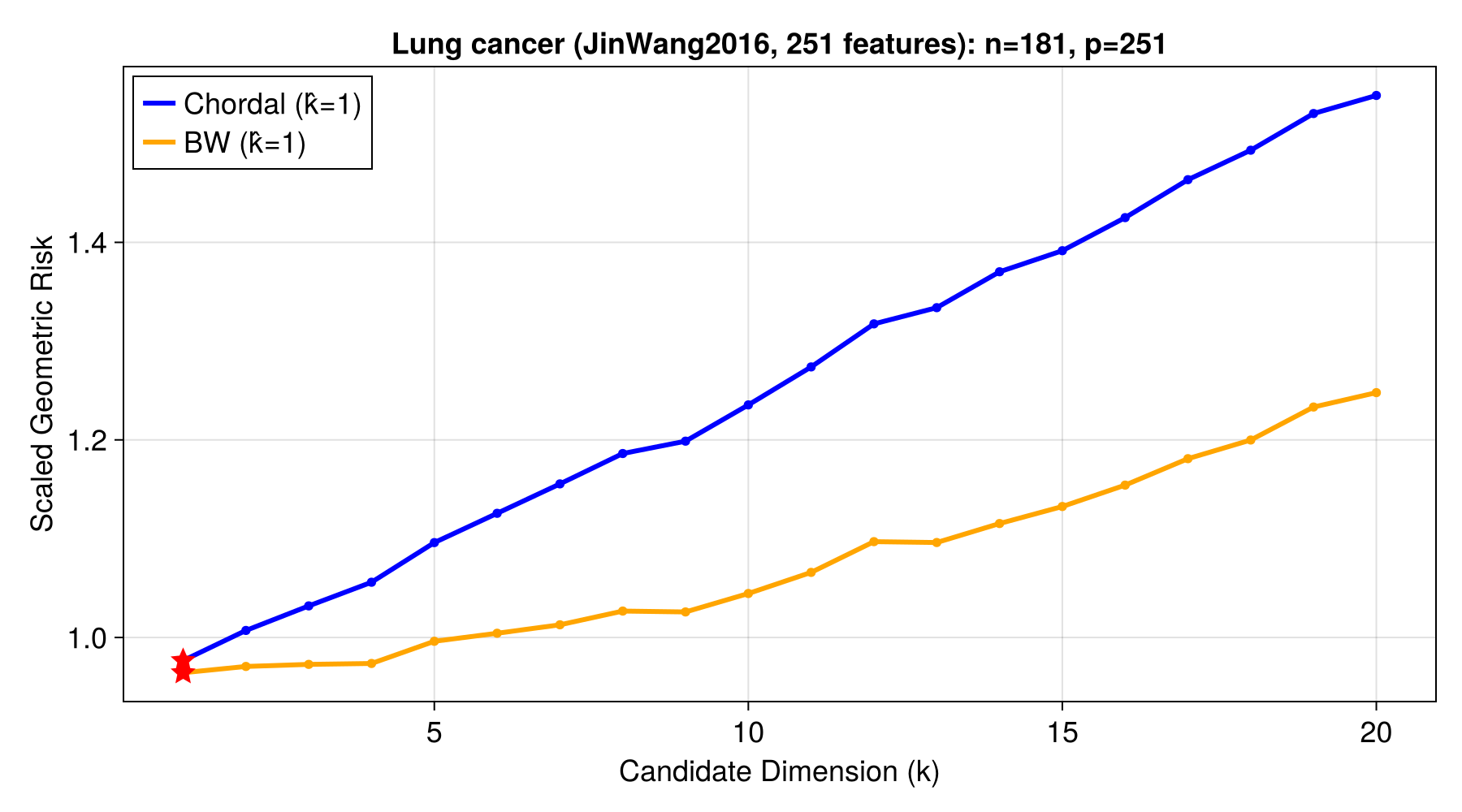}
    \caption{Scaled risk curves on the pruned dataset.}
    \label{fig:pruned-risk}
  \end{subfigure}
  \hfill
  \begin{subfigure}[b]{0.45\textwidth}
    \centering
    \includegraphics[width=\textwidth]{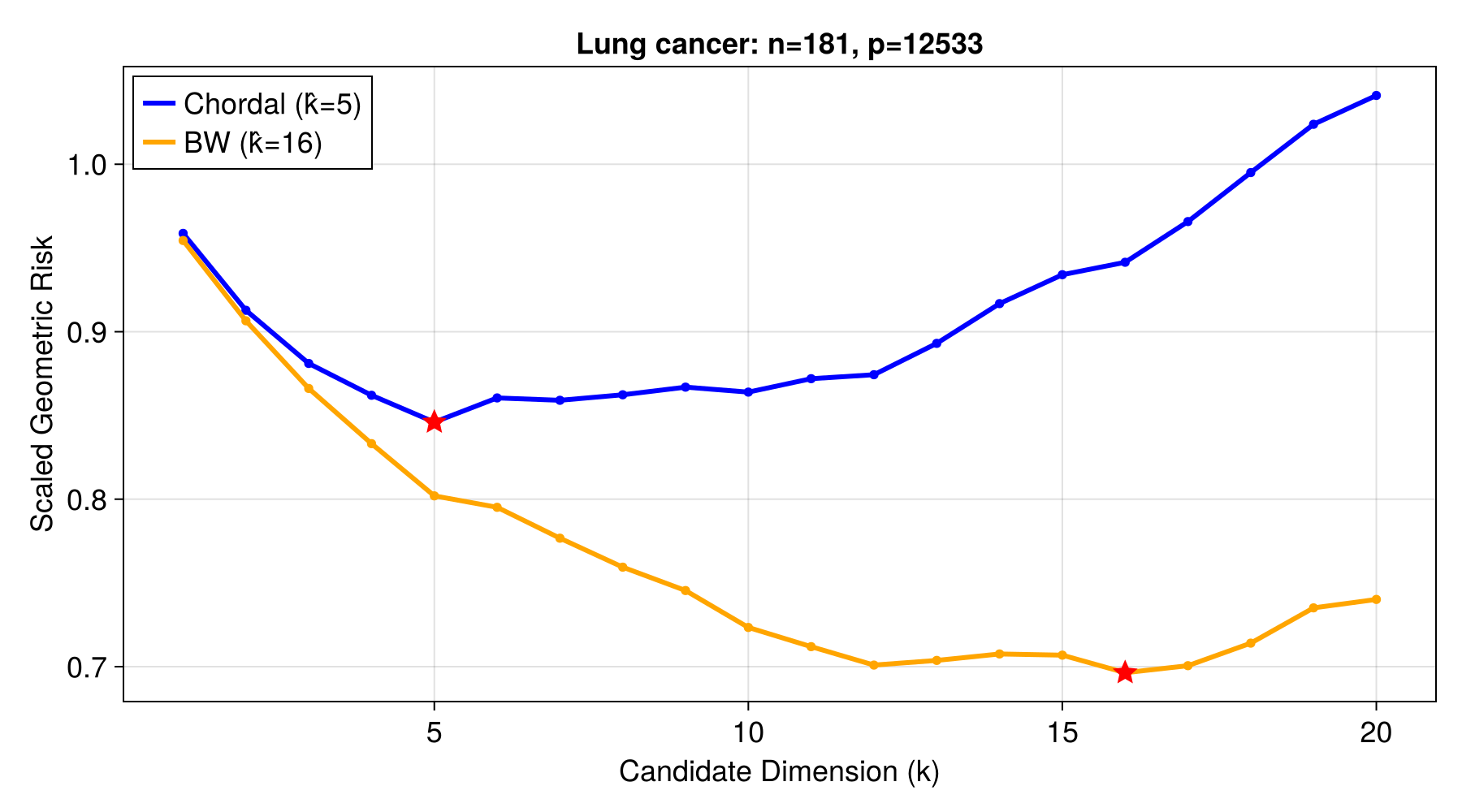}
    \caption{Scaled risk curves on the full dataset.}
    \label{fig:full-risk}
  \end{subfigure}
  \caption{Plots of the scaled risks $\widehat{L}_{\Theta} (k ; \alpha)/D$ and $\widehat{L}_{\c{W}} (k ; \alpha)/D$ on the (a) pruned and (b) full
    lung cancer datasets. $D=20$ and $\alpha=0.5$ in accordance with the theory
    (Proposition ~\ref{prop:half:ideal}).}
    \label{fig:lung_cancer_GDS_risk}
\end{figure}

We start by presenting the analysis of both the Chordal and Bures-Wasserstein
risk functions ($\widehat{L}_{\mathsf{G}} (k ; \alpha)$ for $\mathsf{G} \in \{ \Theta, \c{W} \}$ as defined in ~\eqref{eq:hd:split:criteria}), using an approximately equal split ($\alpha =0.5$) and reference dimension $D = 20$ for both datasets. Figure~\ref{fig:lung_cancer_GDS_risk} displays the criteria on both datasets divided by $D$ for scaling, a normalization that leaves their minimizers unchanged.
As expected, Bures-Wasserstein risk provides larger (or more conservative)
estimates for $d^{*}$ than the Chordal risk. We see that for
both risk curves the $\hat{d}_{\mathsf{G}} = \widehat{L}_{\mathsf{G}} (k ; 0.5)$ estimate increases on the full dataset compared
to the pruned one with $p=251$ only.
We point out here the advantage of a visual demonstration of the Chordal and B-W
methods in Figure~\ref{fig:lung_cancer_GDS_risk} that provides the researcher
with an additional degree of information on how the prediction quality may be
impacted when $\hat{d}_{\mathsf{G}}$ is deviated from the minimum-risk value, a feature that
is not available for comparison with other existing methods. Additionally, as
mentioned earlier, these plots show the fact that the shape of the curve
indicates whether a certain dataset falls into the high-dimensional regimes. The
curves in Figure~\ref{fig:full-risk} are clearly much flatter than the ones in Figure~\ref{fig:pruned-risk}, where such a flat region indicates limited separation between candidates under the observed criterion.

We now compare our proposed methods to existing approaches in the literature,
using the same collection of methods as in Section~\ref{ssec:sim:study}.
The results of the comparative studied are presented in
Table~\ref{tab:lung-cancer-comparison}.
\begin{table}[htbp]
  \centering
  \small
  \begin{tabular}{lrrrrr}
    \hline
    Dataset & GDS-$\Theta$ & GDS-$\c{W}$ & GD & DDPA & BEMA \\
    \hline
    Pruned & 1 & 1 & 39 & 180 & 1 \\
    Full & 5 & 14 & 61 & 180 & 29 \\
    \hline
  \end{tabular}
  \caption{Selected dimensions for the lung cancer datasets. Chordal and B-W use
    the split-sample ($\alpha=0.5$) risks against $D=20$, respectively.}
  \label{tab:lung-cancer-comparison}
\end{table}

We see from Table~\ref{tab:lung-cancer-comparison} that for the pruned dataset,
the Chordal,
B-W and BEMA methods all select the expected $\hat{d}=1$, while GD and DDPA
struggle to identify the true signal and significantly overestimate the number of
spikes. On the other hand, for the full dataset, Chordal and B-W are able to
identify, even under ultra high-dimensionality, the low dimensional noise to a
reasonable extent, while BEMA, GD and DDPA all struggle to identify a
low-dimensional signal space. In particular, as discussed earlier, the pruning
may remove some low level background signal that cannot be captured with only
$p=251$ features. Thus, while the Chordal and B-W approximation on the full
dataset is larger than the $\hat{d}=1$ expected for the pruned dataset, it is
possible that the additional dimensions that are being picked up provide some
smaller underlying signals that are lost or otherwise overwhelmed by noise in the other methods.

\section{Conclusion}
In this paper we have established rigorous theoretical control of two key risk
measures, the Chordal and Bures-Wasserstein for comparing subspaces of different
dimensions. We have shown asymptotically correct selection of the correct
subspace for the regularized empirical risk measures on sample-split data. The results hold
for high-dimensional datasets, giving a powerful tool for working with large
datasets without pre-processing for dimension reduction. Thus, our proposed risk
measures allow the data to speak for itself, even in ultra high-dimensional
settings, as demonstrated by our analysis on lung cancer data. There remain
interesting questions of interest that we leave for future work. In particular, the
applicability of such risk measures to the online learning framework, wherein
new data arrives in intervals is of interest. As it stands, the risk measures
proposed here are unable to identify specific subspace shifts in the data within the
online learning framework. Furthermore, application of the proposed metrics to downstream
estimation models like decisions trees, federated learning with communication
constraints remain interesting lines of future research.

\begin{acks}[Acknowledgments]
This work was primarily carried out while both authors were affiliated with the Institute of Mathematics at EPFL, where Chandak was supported by the Bernoulli Instructorship and Yun was supported by a Swiss National Science Foundation grant.
Chandak is currently supported in part by the Office of the Vice Chancellor for Research and Graduate Education at the University of Wisconsin-Madison, with funding from the Wisconsin Alumni Research Foundation. Yun is currently supported by the Founder's Postdoctoral Fellowship at Columbia University. The authors thank Victor Panaretos for helpful discussions and feedback on this manuscript.
\end{acks}

\bibliographystyle{imsart-number}
\bibliography{bibliography}

@misc{abiyounes2025lanczos,
  author        = {Abi Younes, Charbel and Ding, Xiucai and Trogdon, Thomas},
  title         = {A {Lanczos}-based algorithmic approach for spike detection in large sample covariance matrices},
  year          = {2025},
  howpublished  = {arXiv preprint arXiv:2504.03066},
  eprint        = {2504.03066},
  archivePrefix = {arXiv},
  primaryClass  = {math.ST},
  doi           = {10.48550/arXiv.2504.03066},
  url           = {https://arxiv.org/abs/2504.03066},
}

@article{benaych2011eigenvalues,
  author        = {Benaych-Georges, Florent and Nadakuditi, Raj Rao},
  title         = {The eigenvalues and eigenvectors of finite, low rank perturbations of large random matrices},
  journal       = {Advances in Mathematics},
  volume        = {227},
  number        = {1},
  pages         = {494--521},
  year          = {2011},
  doi           = {10.1016/j.aim.2011.02.007},
}

@book{bhatia1997matrix,
  author        = {Bhatia, Rajendra},
  title         = {Matrix Analysis},
  series        = {Graduate Texts in Mathematics},
  volume        = {169},
  year          = {1997},
  publisher     = {Springer},
  address       = {New York},
  isbn          = {978-0-387-94846-1},
  doi           = {10.1007/978-1-4612-0653-8},
}

@article{blanchard2007statistical,
  author        = {Blanchard, Gilles and Bousquet, Olivier and Zwald, Laurent},
  title         = {Statistical properties of kernel principal component analysis},
  journal       = {Machine Learning},
  volume        = {66},
  number        = {2--3},
  pages         = {259--294},
  year          = {2007},
  doi           = {10.1007/s10994-006-6895-9},
}

@article{bloemendal2016principal,
  author        = {Bloemendal, Alex and Knowles, Antti and Yau, Horng-Tzer and Yin, Jun},
  title         = {On the principal components of sample covariance matrices},
  journal       = {Probability Theory and Related Fields},
  volume        = {164},
  number        = {1--2},
  pages         = {459--552},
  year          = {2016},
  doi           = {10.1007/s00440-015-0616-x},
}

@article{bouwmans2014robust,
  author        = {Bouwmans, Thierry and Zahzah, El Hadi},
  title         = {Robust {PCA} via principal component pursuit: A review for a comparative evaluation in video surveillance},
  journal       = {Computer Vision and Image Understanding},
  volume        = {122},
  pages         = {22--34},
  year          = {2014},
  doi           = {10.1016/j.cviu.2013.11.009},
}

@article{candes2020phase,
  author        = {Cand{\`e}s, Emmanuel J. and Sur, Pragya},
  title         = {The phase transition for the existence of the maximum likelihood estimate in high-dimensional logistic regression},
  journal       = {The Annals of Statistics},
  volume        = {48},
  number        = {1},
  pages         = {27--42},
  year          = {2020},
  doi           = {10.1214/18-aos1789},
}

@article{cattaneo2024convergence,
  author        = {Cattaneo, Matias D. and Chandak, Rajita and Klusowski, Jason M.},
  title         = {Convergence rates of oblique regression trees for flexible function libraries},
  journal       = {The Annals of Statistics},
  volume        = {52},
  number        = {2},
  pages         = {466--490},
  year          = {2024},
  doi           = {10.1214/24-aos2354},
}

@article{cattell1966scree,
  author        = {Cattell, Raymond B.},
  title         = {The scree test for the number of factors},
  journal       = {Multivariate Behavioral Research},
  volume        = {1},
  number        = {2},
  pages         = {245--276},
  year          = {1966},
  doi           = {10.1207/s15327906mbr0102_10},
}

@article{davis1970rotation,
  author        = {Davis, Chandler and Kahan, William Morton},
  title         = {The rotation of eigenvectors by a perturbation. {III}},
  journal       = {SIAM Journal on Numerical Analysis},
  volume        = {7},
  number        = {1},
  pages         = {1--46},
  year          = {1970},
  doi           = {10.1137/0707001},
}

@article{dobriban2015efficient,
  author        = {Dobriban, Edgar},
  title         = {Efficient computation of limit spectra of sample covariance matrices},
  journal       = {Random Matrices: Theory and Applications},
  volume        = {4},
  number        = {4},
  pages         = {1550019},
  year          = {2015},
  doi           = {10.1142/s2010326315500197},
}

@article{dobriban2019deterministic,
  author        = {Dobriban, Edgar and Owen, Art B.},
  title         = {Deterministic parallel analysis: an improved method for selecting factors and principal components},
  journal       = {Journal of the Royal Statistical Society Series B: Statistical Methodology},
  volume        = {81},
  number        = {1},
  pages         = {163--183},
  year          = {2019},
  doi           = {10.1111/rssb.12301},
}

@inproceedings{donhauser2024certified,
  author        = {Donhauser, Konstantin and Lokna, Johan and Sanyal, Amartya and Boedihardjo, March and H{\"o}nig, Robert and Yang, Fanny},
  title         = {Certified private data release for sparse {Lipschitz} functions},
  editor        = {Dasgupta, Sanjoy and Mandt, Stephan and Li, Yingzhen},
  booktitle     = {Proceedings of the 27th International Conference on Artificial Intelligence and Statistics},
  series        = {Proceedings of Machine Learning Research},
  volume        = {238},
  pages         = {1396--1404},
  year          = {2024},
  publisher     = {PMLR},
  url           = {https://proceedings.mlr.press/v238/donhauser24a.html},
}

@article{dowson1982frechet,
  author        = {Dowson, D. C. and Landau, B. V.},
  title         = {The {Fr{\'e}chet} distance between multivariate normal distributions},
  journal       = {Journal of Multivariate Analysis},
  volume        = {12},
  number        = {3},
  pages         = {450--455},
  year          = {1982},
  doi           = {10.1016/0047-259x(82)90077-x},
}

@book{dryden2016statistical,
  author        = {Dryden, Ian L. and Mardia, Kanti V.},
  title         = {Statistical Shape Analysis, with Applications in {R}},
  edition       = {2},
  series        = {Wiley Series in Probability and Statistics},
  year          = {2016},
  publisher     = {John Wiley \& Sons},
  isbn          = {978-0-470-69962-1},
  doi           = {10.1002/9781119072492},
}

@article{gavish2014optimal,
  author        = {Gavish, Matan and Donoho, David L.},
  title         = {The optimal hard threshold for singular values is {$4/\sqrt{3}$}},
  journal       = {IEEE Transactions on Information Theory},
  volume        = {60},
  number        = {8},
  pages         = {5040--5053},
  year          = {2014},
  doi           = {10.1109/tit.2014.2323359},
}

@article{gewers2022principal,
  author        = {Gewers, Felipe L. and Ferreira, Gustavo R. and de Arruda, Henrique F. and Silva, Filipi N. and Comin, Cesar H. and Amancio, Diego R. and Costa, Luciano da F.},
  title         = {Principal component analysis: A natural approach to data exploration},
  journal       = {ACM Computing Surveys},
  volume        = {54},
  number        = {4},
  pages         = {70:1--70:34},
  year          = {2022},
  doi           = {10.1145/3447755},
}

@article{givens1984class,
  author        = {Givens, Clark R. and Shortt, Rae Michael},
  title         = {A class of {Wasserstein} metrics for probability distributions},
  journal       = {Michigan Mathematical Journal},
  volume        = {31},
  number        = {2},
  pages         = {231--240},
  year          = {1984},
  doi           = {10.1307/mmj/1029003026},
}

@inproceedings{guan2009sparse,
  author        = {Guan, Yue and Dy, Jennifer},
  title         = {Sparse probabilistic principal component analysis},
  editor        = {van Dyk, David and Welling, Max},
  booktitle     = {Proceedings of the Twelfth International Conference on Artificial Intelligence and Statistics},
  series        = {Proceedings of Machine Learning Research},
  volume        = {5},
  pages         = {185--192},
  year          = {2009},
  publisher     = {PMLR},
  url           = {https://proceedings.mlr.press/v5/guan09a.html},
}

@article{hasan2021review,
  author        = {Salih Hasan, Basna Mohammed and Abdulazeez, Adnan Mohsin},
  title         = {A review of principal component analysis algorithm for dimensionality reduction},
  journal       = {Journal of Soft Computing and Data Mining},
  volume        = {2},
  number        = {1},
  pages         = {20--30},
  year          = {2021},
  doi           = {10.30880/jscdm.2021.02.01.003},
}

@article{he2025differentially,
  author        = {He, Yiyun and Strohmer, Thomas and Vershynin, Roman and Zhu, Yizhe},
  title         = {Differentially private low-dimensional synthetic data from high-dimensional datasets},
  journal       = {Information and Inference: A Journal of the IMA},
  volume        = {14},
  number        = {1},
  pages         = {iaae034},
  year          = {2025},
  doi           = {10.1093/imaiai/iaae034},
}

@article{horn1965rationale,
  author        = {Horn, John L.},
  title         = {A rationale and test for the number of factors in factor analysis},
  journal       = {Psychometrika},
  volume        = {30},
  number        = {2},
  pages         = {179--185},
  year          = {1965},
  doi           = {10.1007/bf02289447},
}

@article{hotelling1933analysis,
  author        = {Hotelling, Harold},
  title         = {Analysis of a complex of statistical variables into principal components},
  journal       = {Journal of Educational Psychology},
  volume        = {24},
  number        = {6},
  pages         = {417--441},
  year          = {1933},
  doi           = {10.1037/h0071325},
}

@article{hucker2023note,
  author        = {Hucker, Laura and Wahl, Martin},
  title         = {A note on the prediction error of principal component regression in high dimensions},
  journal       = {Theory of Probability and Mathematical Statistics},
  volume        = {109},
  pages         = {37--53},
  year          = {2023},
  doi           = {10.1090/tpms/1196},
}

@article{jewson2020alternative,
  author        = {Jewson, Stephen},
  title         = {An alternative to {PCA} for estimating dominant patterns of climate variability and extremes, with application to {U.S.} and {China} seasonal rainfall},
  journal       = {Atmosphere},
  volume        = {11},
  number        = {4},
  pages         = {354},
  year          = {2020},
  doi           = {10.3390/atmos11040354},
}

@article{jiang2006many,
  author        = {Jiang, Tiefeng},
  title         = {How many entries of a typical orthogonal matrix can be approximated by independent normals?},
  journal       = {The Annals of Probability},
  volume        = {34},
  number        = {4},
  pages         = {1497--1529},
  year          = {2006},
  doi           = {10.1214/009117906000000205},
}

@article{jirak2020perturbation,
  author        = {Jirak, Moritz and Wahl, Martin},
  title         = {Perturbation bounds for eigenspaces under a relative gap condition},
  journal       = {Proceedings of the American Mathematical Society},
  volume        = {148},
  number        = {2},
  pages         = {479--494},
  year          = {2020},
  doi           = {10.1090/proc/14714},
}

@article{mirsky1975trace,
  author        = {Mirsky, L.},
  title         = {A trace inequality of {John von Neumann}},
  journal       = {Monatshefte für Mathematik},
  volume        = {79},
  number        = {4},
  pages         = {303--306},
  year          = {1975},
  doi           = {10.1007/BF01647331},
}

@book{jolliffe2002principal,
  author        = {Jolliffe, Ian T.},
  title         = {Principal Component Analysis},
  edition       = {2},
  series        = {Springer Series in Statistics},
  year          = {2002},
  publisher     = {Springer},
  isbn          = {978-0-387-95442-4},
  doi           = {10.1007/b98835},
}

@article{jolliffe2016principal,
  author        = {Jolliffe, Ian T. and Cadima, Jorge},
  title         = {Principal component analysis: a review and recent developments},
  journal       = {Philosophical Transactions of the Royal Society A: Mathematical, Physical and Engineering Sciences},
  volume        = {374},
  number        = {2065},
  pages         = {20150202},
  year          = {2016},
  doi           = {10.1098/rsta.2015.0202},
}

@article{jombart2010discriminant,
  author        = {Jombart, Thibaut and Devillard, S{\'e}bastien and Balloux, Fran{\c{c}}ois},
  title         = {Discriminant analysis of principal components: a new method for the analysis of genetically structured populations},
  journal       = {BMC Genetics},
  volume        = {11},
  number        = {1},
  pages         = {94},
  year          = {2010},
  doi           = {10.1186/1471-2156-11-94},
}

@article{journee2010generalized,
  author        = {Journ{\'e}e, Michel and Nesterov, Yurii and Richt{\'a}rik, Peter and Sepulchre, Rodolphe},
  title         = {Generalized power method for sparse principal component analysis},
  journal       = {Journal of Machine Learning Research},
  volume        = {11},
  number        = {15},
  pages         = {517--553},
  year          = {2010},
  url           = {https://jmlr.org/papers/v11/journee10a.html},
}

@article{ke2023estimation,
  author        = {Ke, Zheng Tracy and Ma, Yucong and Lin, Xihong},
  title         = {Estimation of the number of spiked eigenvalues in a covariance matrix by bulk eigenvalue matching analysis},
  journal       = {Journal of the American Statistical Association},
  volume        = {118},
  number        = {541},
  pages         = {374--392},
  year          = {2023},
  doi           = {10.1080/01621459.2021.1933497},
}

@inproceedings{menze2011oblique,
  author        = {Menze, Bjoern H. and Kelm, B. Michael and Splitthoff, Daniel N. and Koethe, Ullrich and Hamprecht, Fred A.},
  title         = {On oblique random forests},
  editor        = {Gunopulos, Dimitrios and Hofmann, Thomas and Malerba, Donato and Vazirgiannis, Michalis},
  booktitle     = {Machine Learning and Knowledge Discovery in Databases, Part II},
  series        = {Lecture Notes in Computer Science},
  volume        = {6912},
  pages         = {453--469},
  year          = {2011},
  publisher     = {Springer},
  doi           = {10.1007/978-3-642-23783-6_29},
}

@article{milbradt2020high,
  author        = {Milbradt, Cassandra and Wahl, Martin},
  title         = {High-probability bounds for the reconstruction error of {PCA}},
  journal       = {Statistics \& Probability Letters},
  volume        = {161},
  pages         = {108741},
  year          = {2020},
  doi           = {10.1016/j.spl.2020.108741},
}

@article{muirhead1980asymptotic,
  author        = {Muirhead, Robb J. and Waternaux, Christine M.},
  title         = {Asymptotic distributions in canonical correlation analysis and other multivariate procedures for nonnormal populations},
  journal       = {Biometrika},
  volume        = {67},
  number        = {1},
  pages         = {31--43},
  year          = {1980},
  doi           = {10.1093/biomet/67.1.31},
}

@book{muirhead1982aspects,
  author        = {Muirhead, Robb J.},
  title         = {Aspects of Multivariate Statistical Theory},
  series        = {Wiley Series in Probability and Mathematical Statistics},
  year          = {1982},
  publisher     = {John Wiley \& Sons},
  isbn          = {978-0-471-09442-5},
  doi           = {10.1002/9780470316559},
}

@article{nadler2008finite,
  author        = {Nadler, Boaz},
  title         = {Finite sample approximation results for principal component analysis: A matrix perturbation approach},
  journal       = {The Annals of Statistics},
  volume        = {36},
  number        = {6},
  pages         = {2791--2817},
  year          = {2008},
  doi           = {10.1214/08-aos618},
}

@article{pasini2017principal,
  author        = {Pasini, Giorgia},
  title         = {Principal component analysis for stock portfolio management},
  journal       = {International Journal of Pure and Applied Mathematics},
  volume        = {115},
  number        = {1},
  pages         = {153--167},
  year          = {2017},
  doi           = {10.12732/ijpam.v115i1.12},
}

@article{paul2007asymptotics,
  author        = {Paul, Debashis},
  title         = {Asymptotics of sample eigenstructure for a large dimensional spiked covariance model},
  journal       = {Statistica Sinica},
  volume        = {17},
  number        = {4},
  pages         = {1617--1642},
  year          = {2007},
  url           = {https://www3.stat.sinica.edu.tw/statistica/j17n4/j17n418/j17n418.html},
}

@article{pearson1901lines,
  author        = {Pearson, Karl},
  title         = {{LIII}. On lines and planes of closest fit to systems of points in space},
  journal       = {The London, Edinburgh, and Dublin Philosophical Magazine and Journal of Science},
  volume        = {2},
  number        = {11},
  pages         = {559--572},
  year          = {1901},
  doi           = {10.1080/14786440109462720},
}

@article{reiss2020nonasymptotic,
  author        = {Rei{\ss}, Markus and Wahl, Martin},
  title         = {Nonasymptotic upper bounds for the reconstruction error of {PCA}},
  journal       = {The Annals of Statistics},
  volume        = {48},
  number        = {2},
  pages         = {1098--1123},
  year          = {2020},
  doi           = {10.1214/19-aos1839},
}

@article{shawe2005eigenspectrum,
  author        = {Shawe-Taylor, John and Williams, Christopher K. I. and Cristianini, Nello and Kandola, Jaz},
  title         = {On the eigenspectrum of the {Gram} matrix and the generalization error of kernel-{PCA}},
  journal       = {IEEE Transactions on Information Theory},
  volume        = {51},
  number        = {7},
  pages         = {2510--2522},
  year          = {2005},
  doi           = {10.1109/tit.2005.850052},
}

@inproceedings{simonov2019refined,
  author        = {Simonov, Kirill and Fomin, Fedor and Golovach, Petr and Panolan, Fahad},
  title         = {Refined complexity of {PCA} with outliers},
  editor        = {Chaudhuri, Kamalika and Salakhutdinov, Ruslan},
  booktitle     = {Proceedings of the 36th International Conference on Machine Learning},
  series        = {Proceedings of Machine Learning Research},
  volume        = {97},
  pages         = {5818--5826},
  year          = {2019},
  publisher     = {PMLR},
  url           = {https://proceedings.mlr.press/v97/simonov19a.html},
}

@article{sur2019modern,
  author        = {Sur, Pragya and Cand{\`e}s, Emmanuel J.},
  title         = {A modern maximum-likelihood theory for high-dimensional logistic regression},
  journal       = {Proceedings of the National Academy of Sciences},
  volume        = {116},
  number        = {29},
  pages         = {14516--14525},
  year          = {2019},
  doi           = {10.1073/pnas.1810420116},
}

@article{tomita2020sparse,
  author        = {Tomita, Tyler M. and Browne, James and Shen, Cencheng and Chung, Jaewon and Patsolic, Jesse L. and Falk, Benjamin and Priebe, Carey E. and Yim, Jason and Burns, Randal and Maggioni, Mauro and Vogelstein, Joshua T.},
  title         = {Sparse projection oblique randomer forests},
  journal       = {Journal of Machine Learning Research},
  volume        = {21},
  number        = {104},
  pages         = {1--39},
  year          = {2020},
  url           = {https://jmlr.org/papers/v21/18-664.html},
}

@article{tyler1981asymptotic,
  author        = {Tyler, David E.},
  title         = {Asymptotic inference for eigenvectors},
  journal       = {The Annals of Statistics},
  volume        = {9},
  number        = {4},
  pages         = {725--736},
  year          = {1981},
  doi           = {10.1214/aos/1176345514},
}

@article{vitale2017selecting,
  author        = {Vitale, Raffaele and Westerhuis, Johan A. and N{\ae}s, Tormod and Smilde, Age K. and de Noord, Onno E. and Ferrer, Alberto},
  title         = {Selecting the number of factors in principal component analysis by permutation testing---Numerical and practical aspects},
  journal       = {Journal of Chemometrics},
  volume        = {31},
  number        = {12},
  pages         = {e2937},
  year          = {2017},
  doi           = {10.1002/cem.2937},
}

@article{wang2015singular,
  author        = {Wang, Rongrong},
  title         = {Singular vector perturbation under {Gaussian} noise},
  journal       = {SIAM Journal on Matrix Analysis and Applications},
  volume        = {36},
  number        = {1},
  pages         = {158--177},
  year          = {2015},
  doi           = {10.1137/130938177},
}

@article{wang2026analysis,
  author        = {Wang, Ke},
  title         = {Analysis of singular subspaces under random perturbations},
  journal       = {The Annals of Statistics},
  volume        = {54},
  number        = {2},
  pages         = {667--691},
  year          = {2026},
  doi           = {10.1214/25-aos2582},
}

@article{warmuth2008randomized,
  author        = {Warmuth, Manfred K. and Kuzmin, Dima},
  title         = {Randomized online {PCA} algorithms with regret bounds that are logarithmic in the dimension},
  journal       = {Journal of Machine Learning Research},
  volume        = {9},
  number        = {75},
  pages         = {2287--2320},
  year          = {2008},
  url           = {https://jmlr.org/papers/v9/warmuth08a.html},
}

@article{wedin1972perturbation,
  author        = {Wedin, Per-{\AA}ke},
  title         = {Perturbation bounds in connection with singular value decomposition},
  journal       = {BIT},
  volume        = {12},
  number        = {1},
  pages         = {99--111},
  year          = {1972},
  doi           = {10.1007/bf01932678},
}

@article{yu2015useful,
  author        = {Yu, Yi and Wang, Tengyao and Samworth, Richard J.},
  title         = {A useful variant of the {Davis--Kahan} theorem for statisticians},
  journal       = {Biometrika},
  volume        = {102},
  number        = {2},
  pages         = {315--323},
  year          = {2015},
  doi           = {10.1093/biomet/asv008},
}

@article{zhou2018eigenvalue,
  author        = {Zhou, Yi-Hui and Marron, J. S. and Wright, Fred A.},
  title         = {Eigenvalue significance testing for genetic association},
  journal       = {Biometrics},
  volume        = {74},
  number        = {2},
  pages         = {439--447},
  year          = {2018},
  doi           = {10.1111/biom.12767},
}

@article{zou2006sparse,
  author        = {Zou, Hui and Hastie, Trevor and Tibshirani, Robert},
  title         = {Sparse principal component analysis},
  journal       = {Journal of Computational and Graphical Statistics},
  volume        = {15},
  number        = {2},
  pages         = {265--286},
  year          = {2006},
  doi           = {10.1198/106186006x113430},
}

@article{yun2025high,
  title={High-Dimensional Asymptotics of Differentially Private PCA},
  author={Yun, Youngjoo and Dudeja, Rishabh},
  journal={arXiv preprint arXiv:2511.07270},
  year={2025}
}

@article{gordon2002translation,
  title={Translation of microarray data into clinically relevant cancer diagnostic tests using gene expression ratios in lung cancer and mesothelioma},
  author={Gordon, Gavin J and Jensen, Roderick V and Hsiao, Li-Li and Gullans, Steven R and Blumenstock, Joshua E and Ramaswamy, Sridhar and Richards, William G and Sugarbaker, David J and Bueno, Raphael},
  journal={Cancer research},
  volume={62},
  number={17},
  pages={4963--4967},
  year={2002},
  publisher={American Association for Cancer Research}
}

@article{jin2016influential,
  title={Influential features PCA for high dimensional clustering},
  author={Jin, Jiashun and Wang, Wanjie},
  year={2016}
}

@article{li2008selecting,
  title={Selecting subsets of newly extracted features from PCA and PLS in microarray data analysis},
  author={Li, Guo-Zheng and Bu, Hua-Long and Yang, Mary Qu and Zeng, Xue-Qiang and Yang, Jack Y},
  journal={BMC genomics},
  volume={9},
  number={Suppl 2},
  pages={S24},
  year={2008},
  publisher={Springer}
}

@article{nadler2016discussion,
 ISSN = {00905364},
 URL = {http://www.jstor.org/stable/44245753},
 author = {Boaz Nadler},
 journal = {The Annals of Statistics},
 number = {6},
 pages = {2366--2371},
 publisher = {Institute of Mathematical Statistics},
 title = {DISCUSSION OF "INFLUENTIAL FEATURES PCA FOR HIGH DIMENSIONAL CLUSTERING"},
 urldate = {2026-11-23},
 volume = {44},
 year = {2016}
}

@article{fan2008sure,
  title={Sure independence screening for ultrahigh dimensional feature space},
  author={Fan, Jianqing and Lv, Jinchi},
  journal={Journal of the Royal Statistical Society Series B: Statistical Methodology},
  volume={70},
  number={5},
  pages={849--911},
  year={2008},
  publisher={Oxford University Press}
}

@article{bourgon2010independent,
  title={Independent filtering increases detection power for high-throughput experiments},
  author={Bourgon, Richard and Gentleman, Robert and Huber, Wolfgang},
  journal={Proceedings of the national academy of sciences},
  volume={107},
  number={21},
  pages={9546--9551},
  year={2010},
  publisher={National Academy of Sciences}
}

\newpage

\begin{appendix}

\section{Pure-null comparison}

\begin{table}[htbp]
    \centering
    \small
    \setlength{\tabcolsep}{3pt}
    \caption{(C4) Pure null. Each entry is the mean selected rank ($1000$ replicates), followed by the percentage selecting $d^{*} = 0$. Its complement is the false-positive rate. Boldface marks the highest exact-recovery rate in each row, including ties.}
    \label{tab:sim:c4:append}
    \resizebox{\linewidth}{!}{%
    \begin{tabular}{lrrrrrr}
        \hline
        Noise & $\gamma$ & GDS-$\Theta$ & GDS-$\c{W}$ & BEMA0 & GD & DDPA \\
        \hline
        Gaussian & $0.25$ & 0.002 (99.8\%) & 0.002 (99.8\%) & 0.001 (99.9\%) & \textbf{0.000 (100.0\%)} & 0.090 (91.1\%) \\
         & $0.5$ & 0.003 (99.7\%) & 0.003 (99.7\%) & 0.004 (99.6\%) & \textbf{0.000 (100.0\%)} & 0.104 (89.7\%) \\
         & $1$ & 0.006 (99.4\%) & 0.006 (99.4\%) & 0.004 (99.6\%) & \textbf{0.000 (100.0\%)} & 0.108 (89.6\%) \\
         & $2$ & 0.007 (99.3\%) & 0.007 (99.3\%) & 0.001 (99.9\%) & \textbf{0.000 (100.0\%)} & 0.104 (89.6\%) \\
        \hline
        $t_{5}$ & $0.25$ & \textbf{0.002 (99.8\%)} & \textbf{0.002 (99.8\%)} & 0.139 (86.7\%) & 0.012 (98.8\%) & 0.076 (92.7\%) \\
         & $0.5$ & 0.007 (99.4\%) & \textbf{0.005 (99.5\%)} & 0.148 (85.5\%) & 0.008 (99.2\%) & 0.087 (91.8\%) \\
         & $1$ & \textbf{0.005 (99.5\%)} & \textbf{0.005 (99.5\%)} & 0.150 (85.6\%) & 0.008 (99.2\%) & 0.100 (90.3\%) \\
         & $2$ & \textbf{0.003 (99.7\%)} & \textbf{0.003 (99.7\%)} & 0.109 (89.8\%) & 0.007 (99.3\%) & 0.097 (90.3\%) \\
        \hline
        $t_{8}$ & $0.25$ & 0.003 (99.7\%) & 0.003 (99.7\%) & 0.002 (99.8\%) & \textbf{0.000 (100.0\%)} & 0.107 (89.9\%) \\
         & $0.5$ & 0.001 (99.9\%) & 0.001 (99.9\%) & 0.003 (99.7\%) & \textbf{0.000 (100.0\%)} & 0.091 (91.3\%) \\
         & $1$ & 0.004 (99.6\%) & 0.004 (99.6\%) & 0.006 (99.4\%) & \textbf{0.000 (100.0\%)} & 0.102 (90.2\%) \\
         & $2$ & 0.001 (99.9\%) & 0.001 (99.9\%) & 0.006 (99.4\%) & \textbf{0.000 (100.0\%)} & 0.106 (89.6\%) \\
        \hline
        Mixture, $a = 10$ & $0.25$ & 0.005 (99.6\%) & 0.004 (99.6\%) & 0.266 (75.4\%) & \textbf{0.001 (99.9\%)} & 0.067 (93.7\%) \\
         & $0.5$ & 0.007 (99.3\%) & 0.007 (99.3\%) & 0.167 (84.3\%) & \textbf{0.001 (99.9\%)} & 0.067 (93.4\%) \\
         & $1$ & 0.011 (99.0\%) & 0.011 (99.0\%) & 0.070 (93.3\%) & \textbf{0.000 (100.0\%)} & 0.094 (91.2\%) \\
         & $2$ & \textbf{0.000 (100.0\%)} & \textbf{0.000 (100.0\%)} & 0.018 (98.2\%) & \textbf{0.000 (100.0\%)} & 0.100 (90.4\%) \\
        \hline
        Mixture, $a = 20$ & $0.25$ & 0.069 (95.7\%) & 0.055 (96.1\%) & 4.002 (0.6\%) & 0.995 (36.3\%) & \textbf{0.022 (97.9\%)} \\
         & $0.5$ & 0.093 (94.9\%) & 0.069 (95.5\%) & 5.348 (0.0\%) & 1.000 (36.1\%) & \textbf{0.040 (96.6\%)} \\
         & $1$ & 0.154 (94.9\%) & \textbf{0.073 (96.2\%)} & 6.472 (0.0\%) & 0.684 (49.4\%) & 0.060 (94.3\%) \\
         & $2$ & 0.120 (96.4\%) & \textbf{0.047 (97.3\%)} & 6.803 (0.0\%) & 0.209 (81.0\%) & 0.121 (89.7\%) \\
        \hline
        Mixture, $a = 30$ & $0.25$ & 0.249 (90.5\%) & 0.117 (93.8\%) & 7.850 (0.0\%) & 3.954 (1.3\%) & \textbf{0.021 (98.0\%)} \\
         & $0.5$ & 0.477 (87.9\%) & 0.210 (92.1\%) & 12.257 (0.0\%) & 5.654 (0.1\%) & \textbf{0.033 (96.8\%)} \\
         & $1$ & 0.722 (85.9\%) & 0.213 (91.8\%) & 17.790 (0.0\%) & 6.571 (0.0\%) & \textbf{0.065 (93.7\%)} \\
         & $2$ & 0.887 (86.8\%) & \textbf{0.247 (91.7\%)} & 22.879 (0.0\%) & 5.441 (0.4\%) & 0.137 (87.7\%) \\
        \hline
    \end{tabular}}
\end{table}

\section{Proofs for Section~\ref{sec:prelim}}
\label{app:prelim}

\begin{proof}[Proof of \cref{prop:matrix_product}]
Write $\m{C} \defeq \m{U}^{\top}\m{V}$.
Because $\m{\Pi}_{\c{U}}\m{\Pi}_{\c{V}} = \m{U}\m{C}\m{V}^{\top}$,
multiplication by the semi-orthogonal matrices preserves the
nonzero singular values of $\m{C}$.
Note that their upper bound trivially follows from
$\vertii{\m{\Pi}_{\c{U}}\m{\Pi}_{\c{V}}}_{\infty} \le 1$.
Furthermore,
\begin{align}\label{eq:uni:equi}
    \m{U}^{\top}\m{\Pi}_{\c{V}}\m{U} = \m{C}\m{C}^{\top}, \quad
    \m{\Pi}_{\c{U}}\m{\Pi}_{\c{V}}\m{\Pi}_{\c{U}} = \m{U}(\m{C}\m{C}^{\top})\m{U}^{\top}, \quad
    \m{\Pi}_{\c{V}}\m{\Pi}_{\c{U}}\m{\Pi}_{\c{V}} = \m{V}(\m{C}^{\top}\m{C})\m{V}^{\top}.
\end{align}
These matrices have the same nonzero eigenvalues, namely the
squared nonzero singular values of $\m{C}$.

Suppose these singular values all equal one.
Then $\m{\Pi}_{\c{U}}\m{\Pi}_{\c{V}}\m{\Pi}_{\c{U}}$ is idempotent, and
\begin{align*}
    \vertiii{(\m{I}_{p}-\m{\Pi}_{\c{U}})\m{\Pi}_{\c{V}}\m{\Pi}_{\c{U}}}_{2}^{2}
    = \tr\big[
        \m{\Pi}_{\c{U}}\m{\Pi}_{\c{V}}(\m{I}_{p}-\m{\Pi}_{\c{U}})\m{\Pi}_{\c{V}}\m{\Pi}_{\c{U}}
        \big]
    = \tr\big[
        \m{\Pi}_{\c{U}}\m{\Pi}_{\c{V}}\m{\Pi}_{\c{U}}-(\m{\Pi}_{\c{U}}\m{\Pi}_{\c{V}}\m{\Pi}_{\c{U}})^{2}
        \big]
     = 0.
\end{align*}
Hence $\m{\Pi}_{\c{V}}\m{\Pi}_{\c{U}}=\m{\Pi}_{\c{U}}\m{\Pi}_{\c{V}}\m{\Pi}_{\c{U}}$.
Taking transposes gives $\m{\Pi}_{\c{U}}\m{\Pi}_{\c{V}}=\m{\Pi}_{\c{U}}\m{\Pi}_{\c{V}}\m{\Pi}_{\c{U}}$, so $\m{\Pi}_{\c{U}}$ and $\m{\Pi}_{\c{V}}$ commute. Their product is an orthogonal projection whose range is $\c{U}\cap\c{V}$, proving \eqref{eq:SVD:commute}.
Conversely, if the projections commute, their product is the orthogonal projection onto $\c{U}\cap\c{V}$.
Its nonzero singular values, and therefore those of $\m{C}$, all equal one.
\end{proof}

\begin{proof}[Proof of \cref{thm:procrustes}]
The case $d_{\c{V}}=0$ follows directly from the definitions.
Assume $d_{\c{V}}>0$ and write
$\m{C}\defeq\m{U}^{\top}\m{V}$.
For any $\m{R}\in\c{O}(d_{\c{U}},d_{\c{V}})$,
orthonormality of the columns of $\m{U}$, $\m{V}$, and $\m{R}$
gives
\begin{align*}
    \vertiii{\m{U}-\m{V}\m{R}^{\top}}_{2}^{2}
    =
    d_{\c{U}}+d_{\c{V}}-2\tr(\m{R}^{\top}\m{C}), \qquad
    \vertiii{\m{U}\m{R}-\m{V}}_{2}^{2}
    =
    2d_{\c{V}}-2\tr(\m{R}^{\top}\m{C}).
\end{align*}
Thus both minimizations reduce to maximizing
$\tr(\m{R}^{\top}\m{C})$.
Take a thin singular value decomposition
$\m{C}=\m{L}\m{S}\m{H}^{\top}$, where
$\m{L}\in\c{O}(d_{\c{U}},d_{\c{V}})$,
$\m{H}\in\c{O}(d_{\c{V}})$, and $\m{S}\in\c{D}(d_{\c{V}})$ is diagonal
with nonnegative entries.
The von-Neumann trace inequality \cite[Eq.~(1)]{mirsky1975trace} for singular values yields $\tr(\m{R}^{\top}\m{C})\le \vertii{\m{C}}_{1} = \tr(\m{S})$,
with equality attained at $\m{R}=\m{L}\m{H}^{\top}$.
Consequently,
\begin{align*}
    \max_{\m{R}\in\c{O}(d_{\c{U}},d_{\c{V}})}
    \tr(\m{R}^{\top}\m{C})
    = \vertii{\m{C}}_{1}.
\end{align*}
Substituting this maximum into the two expansions and using
$l_{\c{W}}(\c{U},\c{V})^{2}
=d_{\c{U}}+d_{\c{V}}-2\vertii{\m{C}}_{1}$
proves both identities.
\end{proof}

\section{Proofs for Section~\ref{sec:main_results}}
\label{app:main_proofs}

\begin{proof}[Proof of Proposition~\ref{prop:kahan:v1}]
If $\c{I}^{c}=\varnothing$ or $\c{J}=\varnothing$, the inequality
holds by convention of $0\cdot\infty=0$ and $0^{0}=1$.
Assume henceforth that both sets are nonempty, i.e., $m\defeq\min(p-d_{\c{U}},d_{\c{V}}) \in \b{N}$.
Let $\m{E}\defeq\hat{\m{\Sigma}}-\m{\Sigma} \in \b{R}^{p \times p}$ and $\m{M}\defeq\m{U}_{\c{I}^{c}}^{\top}\hat{\m{V}}_{\c{J}}$.
Since $\hat{\m{\Sigma}} \hat{\m{V}}_{\c{J}} = \hat{\m{V}}_{\c{J}} \hat{\m{\Omega}}_{\c{J}}$ and $\m{U}_{\c{I}^{c}}^{\top} \m{\Sigma} = \m{\Lambda}_{\c{I}^{c}} \m{U}_{\c{I}^{c}}^{\top}$, it holds that
\begin{align*}
    \m{U}_{\c{I}^{c}}^{\top}\m{E}\hat{\m{V}}_{\c{J}}
  = \m{U}_{\c{I}^{c}}^{\top} (\hat{\m{\Sigma}} \hat{\m{V}}_{\c{J}}) - (\m{U}_{\c{I}^{c}}^{\top} \m{\Sigma}) \m{U}_{\c{J}}
  = (\m{U}_{\c{I}^{c}}^{\top} \hat{\m{V}}_{\c{J}}) \hat{\m{\Omega}}_{\c{J}}
  - \m{\Lambda}_{\c{I}^{c}} (\m{U}_{\c{I}^{c}}^{\top} \hat{\m{V}}_{\c{J}})
  = \m{M} \hat{\m{\Omega}}_{\c{J}} - \m{\Lambda}_{\c{I}^{c}} \m{M},
\end{align*}
Taking the squared Frobenius norm gives
\begin{align*}
    \vertii{\m{U}_{\c{I}^{c}}^{\top} \m{E} \hat{\m{V}}_{\c{J}}}_{2}^{2} =
    \sum_{i\in\c{I}^{c}}\sum_{j\in\c{J}}
    (\hat{\omega}_j-\lambda_i)^{2}\m{M}_{ij}^{2} \ge
    \delta^{2}\vertii{\m{M}}_{2}^{2},
\end{align*}
where $\m{M}_{ij}=\m{u}_{i}^{\top}\hat{\m{v}}_{j}$.
To bound the perturbation from above, we apply the generalized Hölder's inequality for Schatten norms \cite[Exercise IV.2.7]{bhatia1997matrix}. For any triplet $q, r, s \in [2, \infty]$ satisfying $q^{-1}+r^{-1}+s^{-1}=1/2$, we have:
\begin{align*}
    \delta \vertii{\m{M}}_{2} \le \vertii{\m{U}_{\c{I}^{c}}^{\top} \m{E} \hat{\m{V}}_{\c{J}}}_{2} \le \vertii{\m{U}_{\c{I}^{c}}^{\top}}_{r} \vertii{\m{E}}_{q} \vertii{\hat{\m{V}}_{\c{J}}}_{s} = (p-d_{\c{U}})^{1/r} d_{\c{V}}^{1/s} \vertii{\m{E}}_{q}.
\end{align*}
Taking either $r = \infty$ or $s = \infty$ provides the tightest bound: $\delta \vertiii{\m{M}}_{2} \le \min(p-d_{\c{U}}, d_{\c{V}})^{1/2-1/q} \vertii{\m{E}}_{q}$.
\end{proof}

\begin{proof}[Proof of Theorem~\ref{thm:kahan:v2_generalized}]
It is trivial that $\delta_{\m{\Sigma}} > 0$ implies $\c{J}\subseteq\c{I}$. If $\c{J}=\varnothing$ or $\c{I}^{c}=\varnothing$, the excess squared distance and the upper bound are both zero. We therefore assume $0<d_{\c{V}}\le d_{\c{U}}<p$.
Consider
\begin{align*}
    \m{Q} \defeq \hat{\m{V}}_{\c{J}}\m{\Lambda}_{\c{J}}-\m{\Sigma}\hat{\m{V}}_{\c{J}}
    = \hat{\m{V}}_{\c{J}}(\m{\Lambda}_{\c{J}}-\hat{\m{\Omega}}_{\c{J}})+(\hat{\m{\Sigma}}-\m{\Sigma})\hat{\m{V}}_{\c{J}}.
\end{align*}
We bound the Frobenius norm of this residual from above by the perturbation norm and from below by the population gap.

\noindent\textbf{Upper bound.}
Applying the triangle inequality and Cauchy-Schwarz yields:
\begin{align*}
    \vertii{\m{Q}}_{2}^{2}
    \le 2\vertii{\hat{\m{V}}_{\c{J}}(\m{\Lambda}_{\c{J}}-\hat{\m{\Omega}}_{\c{J}})}_{2}^{2}
        +2\vertii{(\hat{\m{\Sigma}}-\m{\Sigma})\hat{\m{V}}_{\c{J}}}_{2}^{2}
    \le 2d_{\c{V}}^{1-2/q}\left(\vertii{\m{\Lambda}_{\c{J}}-\hat{\m{\Omega}}_{\c{J}}}_{q}^{2}
        +\vertii{\hat{\m{\Sigma}}-\m{\Sigma}}_{q}^{2}\right).
\end{align*}
where we utilized the generalized Hölder's inequality for unitarily invariant norms \cite[Exercise IV.2.7]{bhatia1997matrix}.
Using the $\ell_{q}$-norm as a symmetric gauge function in Lidskii's Theorem \cite[Theorem III.4.4]{bhatia1997matrix}, we bound the eigenvalue perturbations as $\vertii{\m{\Lambda}_{\c{J}}-\hat{\m{\Omega}}_{\c{J}}}_{q} \le \vertii{\m{\Lambda}-\hat{\m{\Omega}}}_{q} \le \vertii{\hat{\m{\Sigma}}-\m{\Sigma}}_{q}$. Substituting this provides the upper bound $\vertii{\m{Q}}_{2}^{2}\le 4d_{\c{V}}^{1-2/q}\vertii{\hat{\m{\Sigma}}-\m{\Sigma}}_{q}^{2}$.

\noindent\textbf{Lower bound.}
As in the proof of Proposition~\ref{prop:kahan:v1}, defining $\m{M} :=
\m{U}_{\c{I}^{c}}^{\top} \widehat{\m{V}}_{\c{J}} \in \b{R}^{(p-d_{\c{U}}) \times
d_{\c{V}}}$ allows for expressing
$\m{U}_{\c{I}^{c}}^{\top}\m{Q} = \m{M}\m{\Lambda}_{\c{J}}-\m{\Lambda}_{\c{I}^{c}}\m{M}$.
Using the orthogonal decomposition $\b{R}^{p}=\c{U}_{\c{I}}\oplus\c{U}_{\c{I}^{c}}$, we obtain
\begin{align*}
    \vertii{\m{Q}}_{2}^{2}
    = \vertii{\m{U}_{\c{I}}^{\top}\m{Q}}_{2}^{2}+\vertii{\m{U}_{\c{I}^{c}}^{\top}\m{Q}}_{2}^{2} \ge \sum_{i\in\c{I}^{c}}\sum_{j\in\c{J}}(\lambda_j-\lambda_i)^{2}\m{M}_{ij}^{2}
    \ge \delta_{\m{\Sigma}}^{2}\vertii{\m{M}}_{2}^{2}.
\end{align*}
Combining these bounds with the identity \eqref{eq:chordal:iden}
completes the proof.
\end{proof}

\begin{proof}[Proof of Theorem~\ref{thm:general_asymptotic_dist}]
If $\c{I}=[p]$, $\c{I}^{c} = \varnothing$ and both the excess distance and the limiting sum are zero. Assume henceforth that $\c{I}^{c}\ne\varnothing$, and write $\m{E}_n\defeq\hat{\m{\Sigma}}_n-\m{\Sigma}$ and $\m{Z}_n\defeq n^{\zeta}\m{E}_n\xrightarrow{d}\m{Z}$. In particular, $\m{Z}_n=O_{\b{P}}(1)$ and $\m{E}_n=O_{\b{P}}(n^{-\zeta})$.
We first align the empirical basis with $\m{U}_{\c{J}}$, which has the same number of columns. Since $\tilde{\delta}_{\m{\Sigma}}>0$, we may apply Theorem~\ref{thm:kahan:v2_generalized} with both index sets equal to $\c{J}$ and $q=2$. Together with Theorem~\ref{thm:procrustes}, this gives an orthogonal matrix $\m{R}_n\in\c{O}(d_{\c{V}})$ such that
\begin{align*}
    \tilde{\m{V}}_{\c{J}}\defeq\hat{\m{V}}_{\c{J}}\m{R}_n \in \c{O}(\hat{\c{V}}_{\c{J}}),
    \qquad
    \vertii{\tilde{\m{V}}_{\c{J}}-\m{U}_{\c{J}}}_{2}
    =l_{\c{W}}(\c{U}_{\c{J}},\hat{\c{V}}_{\c{J}})
    \le l_{\Theta}(\c{U}_{\c{J}},\hat{\c{V}}_{\c{J}})
    =O_{\b{P}}(n^{-\zeta}).
\end{align*}

Consider the Sylvester-type linear operator
\begin{align*}
    \c{L}:\b{R}^{(p-d_{\c{U}})\times d_{\c{V}}}\to\b{R}^{(p-d_{\c{U}})\times d_{\c{V}}},
    \qquad
    \m{M}\mapsto\m{M}\m{\Lambda}_{\c{J}}-\m{\Lambda}_{\c{I}^{c}}\m{M}.
\end{align*}
Since $\c{I}^{c}\subseteq\c{J}^{c}$, we have $|\lambda_j-\lambda_i|\ge\tilde{\delta}_{\m{\Sigma}}>0$ for every $i\in\c{I}^{c}$ and $j\in\c{J}$. Thus $\c{L}$ is bijective, with bounded inverse given component-wise by $(\c{L}^{-1}[\m{M}])_{ij} = \m{M}_{ij}/(\lambda_j-\lambda_i)$ for $i\in\c{I}^{c}, j\in\c{J}$.
We claim that it suffices to prove
\begin{align}\label{eq:sylver:goal}
    \m{U}_{\c{I}^{c}}^{\top}\m{Z}_n\m{U}_{\c{J}}
    -n^{\zeta}\c{L}[\widehat{\m{M}}]
    =o_{\b{P}}(1), \qquad \widehat{\m{M}}\defeq\m{U}_{\c{I}^{c}}^{\top}\tilde{\m{V}}_{\c{J}}
\end{align}
Indeed, applying $\c{L}^{-1}$ and using the continuous mapping theorem and Slutsky's theorem would give
\begin{align*}
    n^{\zeta}\widehat{\m{M}}
    =\c{L}^{-1}\left[\m{U}_{\c{I}^{c}}^{\top}\m{Z}_n\m{U}_{\c{J}}\right]+o_{\b{P}}(1)
    \xrightarrow{d}\c{L}^{-1}[\tilde{\m{Z}}],
    \qquad
    \tilde{\m{Z}}\defeq\m{U}_{\c{I}^{c}}^{\top}\m{Z}\m{U}_{\c{J}}.
\end{align*}
Because right multiplication by $\m{R}_n$ in the definition of $\widehat{\m{M}}$ preserves the Frobenius norm, we get
\begin{align*}
    n^{2\zeta}\left[l_{\Theta}(\c{U}_{\c{I}},\hat{\c{V}}_{\c{J}})^2-(d_{\c{U}}-d_{\c{V}})\right]
    =2\vertii{n^{\zeta}\widehat{\m{M}}}_{2}^{2}
    \xrightarrow{d}2\vertii{\c{L}^{-1}[\tilde{\m{Z}}]}_{2}^{2}
    =2\sum_{i\in\c{I}^{c}}\sum_{j\in\c{J}}
    \left(\frac{\tilde{\m{Z}}_{ij}}{\lambda_j-\lambda_i}\right)^2.
\end{align*}

To prove \eqref{eq:sylver:goal}, define $\m{H}_{\c{J}}\defeq\m{R}_n^{\top}\hat{\m{\Omega}}_{\c{J}}\m{R}_n =\tilde{\m{V}}_{\c{J}}^{\top}\hat{\m{\Sigma}}_n\tilde{\m{V}}_{\c{J}}$. We then use
\begin{align*}
    \hat{\m{\Sigma}}_n \tilde{\m{V}}_{\c{J}}
    = \hat{\m{V}}_{\c{J}} \hat{\m{\Omega}}_{\c{J}} \m{R}_{n} = \tilde{\m{V}}_{\c{J}} (\m{R}_n^{\top}\hat{\m{\Omega}}_{\c{J}}\m{R}_n)
    =\tilde{\m{V}}_{\c{J}}\m{H}_{\c{J}}, \qquad \m{U}_{\c{I}^{c}}^{\top}\m{\Sigma}=\m{\Lambda}_{\c{I}^{c}}\m{U}_{\c{I}^{c}}^{\top},
\end{align*}
to obtain
\begin{align*}
    \m{U}_{\c{I}^{c}}^{\top}\m{Z}_n\tilde{\m{V}}_{\c{J}}
    &=n^{\zeta}\m{U}_{\c{I}^{c}}^{\top}(\hat{\m{\Sigma}}_n-\m{\Sigma})\tilde{\m{V}}_{\c{J}} =n^{\zeta}\left(\m{U}_{\c{I}^{c}}^{\top}\tilde{\m{V}}_{\c{J}}\m{H}_{\c{J}}
    -\m{\Lambda}_{\c{I}^{c}}\m{U}_{\c{I}^{c}}^{\top}\tilde{\m{V}}_{\c{J}}\right) =n^{\zeta}\left(\widehat{\m{M}}\m{H}_{\c{J}}-\m{\Lambda}_{\c{I}^{c}}\widehat{\m{M}}\right).
\end{align*}
Consequently, the left hand side of \eqref{eq:sylver:goal} decomposes as
\begin{align*}
    \m{U}_{\c{I}^{c}}^{\top}\m{Z}_n\m{U}_{\c{J}}-n^{\zeta}\c{L}[\widehat{\m{M}}]
    &=
    \underbrace{\widehat{\m{M}}\left[n^{\zeta}(\m{H}_{\c{J}}-\m{\Lambda}_{\c{J}})\right]}_{(A)}
    -\underbrace{\m{U}_{\c{I}^{c}}^{\top}\m{Z}_n(\tilde{\m{V}}_{\c{J}}-\m{U}_{\c{J}})}_{(B)}.
\end{align*}
We bound these terms separately.
\begin{itemize}[leftmargin=*]
\item[(A)] Since $\c{J}\subseteq\c{I}$, we have $\m{U}_{\c{I}^{c}}^{\top}\m{U}_{\c{J}}=\m{0}$. Hence the alignment bound gives $\widehat{\m{M}}
    =\m{U}_{\c{I}^{c}}^{\top}(\tilde{\m{V}}_{\c{J}}-\m{U}_{\c{J}})
    =O_{\b{P}}(n^{-\zeta})$.
On the other hand, using $\m{\Lambda}_{\c{J}}=\m{U}_{\c{J}}^{\top}\m{\Sigma}\m{U}_{\c{J}}$, we obtain
\begin{align*}
    \m{H}_{\c{J}}-\m{\Lambda}_{\c{J}}
    &=\tilde{\m{V}}_{\c{J}}^{\top}\m{E}_n\tilde{\m{V}}_{\c{J}}
    +\tilde{\m{V}}_{\c{J}}^{\top}\m{\Sigma}\tilde{\m{V}}_{\c{J}}
    -\m{U}_{\c{J}}^{\top}\m{\Sigma}\m{U}_{\c{J}}\\
    &=\tilde{\m{V}}_{\c{J}}^{\top}\m{E}_n\tilde{\m{V}}_{\c{J}}
    +(\tilde{\m{V}}_{\c{J}}-\m{U}_{\c{J}})^{\top}\m{\Sigma}\tilde{\m{V}}_{\c{J}}
    +\m{U}_{\c{J}}^{\top}\m{\Sigma}(\tilde{\m{V}}_{\c{J}}-\m{U}_{\c{J}}).
\end{align*}
Since both bases have orthonormal columns, the triangle inequality yields
\begin{align*}
    \vertii{\m{H}_{\c{J}}-\m{\Lambda}_{\c{J}}}_{2}
    \le\vertii{\m{E}_n}_{2}
    +2\vertii{\m{\Sigma}}_{\infty}\vertii{\tilde{\m{V}}_{\c{J}}-\m{U}_{\c{J}}}_{2}
    =O_{\b{P}}(n^{-\zeta}).
\end{align*}
Therefore, $(A)=O_{\b{P}}(n^{-\zeta})O_{\b{P}}(1)=o_{\b{P}}(1)$.

\item[(B)] Since $\m{U}_{\c{I}^{c}}^{\top}\m{Z}_n=O_{\b{P}}(1)$ and $\tilde{\m{V}}_{\c{J}}-\m{U}_{\c{J}}=O_{\b{P}}(n^{-\zeta})$, we have $(B)=O_{\b{P}}(n^{-\zeta})=o_{\b{P}}(1)$.
\end{itemize}
Combining these bounds establishes \eqref{eq:sylver:goal} and completes the proof.
\end{proof}

\begin{proof}[Proof of Corollary~\ref{cor:asymp_sample_cov}]
The finite fourth moment assumption and the multivariate CLT imply $\sqrt{n}(\hat{\m{\Sigma}}_n-\m{\Sigma})\xrightarrow{d}\m{Z}$, where $\operatorname{vec}(\m{Z})\sim\c{N}(\m{0},\m{C})$.
Applying Theorem~\ref{thm:general_asymptotic_dist} with $\zeta=1/2$ gives the asserted limit with $\tilde{\m{Z}}=\m{U}_{\c{I}^{c}}^{\top}\m{Z}\m{U}_{\c{J}}$. Its covariance follows from $\operatorname{vec}(\tilde{\m{Z}})=(\m{U}_{\c{J}}^{\top}\otimes\m{U}_{\c{I}^{c}}^{\top})\operatorname{vec}(\m{Z})$.
\end{proof}

\begin{proof}[Proof of Corollary~\ref{cor:asymp:ellip}]
For an elliptical distribution with covariance $\m{\Sigma}$ and kurtosis parameter $\kappa$, the covariance matrix $\m{C}$ in Corollary~\ref{cor:asymp_sample_cov} satisfies:
\begin{align*}
    \m{C}
    =(1+\kappa)(\m{I}_{p^2}+\m{K}_{p,p})(\m{\Sigma}\otimes\m{\Sigma})
    +\kappa\operatorname{vec}(\m{\Sigma})\operatorname{vec}(\m{\Sigma})^{\top},
\end{align*}
where $\m{K}_{p,p}$ is the commutation matrix \cite{muirhead1980asymptotic}.
Since $\c{J}\subseteq\c{I}$, we have $\m{U}_{\c{I}^{c}}^{\top}\m{\Sigma}\m{U}_{\c{J}}=\m{0}$.
Because $\c{J} \subseteq \c{I}$, the orthogonality $\m{U}_{\c{I}^{c}}^{\top} \m{\Sigma} \m{U}_{\c{J}} = \m{0}$ holds, causing the cross-terms in $\operatorname{Cov}[\operatorname{vec}(\tilde{\m{Z}})] = (\m{U}_{\c{J}}^{\top} \otimes \m{U}_{\c{I}^{c}}^{\top}) \m{C} (\m{U}_{\c{J}} \otimes \m{U}_{\c{I}^{c}})$ to vanish:
\begin{align*}
    (\m{U}_{\c{J}}^{\top} \otimes \m{U}_{\c{I}^{c}}^{\top}) \vect(\m{\Sigma})\vect(\m{\Sigma})^{\top} (\m{U}_{\c{J}} \otimes \m{U}_{\c{I}^{c}}) &= \vect(\m{U}_{\c{I}^{c}}^{\top} \m{\Sigma} \m{U}_{\c{J}}) \vect(\m{U}_{\c{I}^{c}}^{\top} \m{\Sigma} \m{U}_{\c{J}})^{\top} = \m{0}, \\
    (\m{U}_{\c{J}}^{\top} \otimes \m{U}_{\c{I}^{c}}^{\top}) (\m{\Sigma} \otimes \m{\Sigma}) \m{K}_{p,p} (\m{U}_{\c{J}} \otimes \m{U}_{\c{I}^{c}}) &= (\m{U}_{\c{J}}^{\top} \m{\Sigma} \otimes \m{U}_{\c{I}^{c}}^{\top} \m{\Sigma}) (\m{U}_{\c{I}^{c}} \otimes \m{U}_{\c{J}}) \m{K}_{d_{\c{V}}, p-d_{\c{U}}} \\
    &= (\m{U}_{\c{J}}^{\top} \m{\Sigma} \m{U}_{\c{I}^{c}} \otimes \m{U}_{\c{I}^{c}}^{\top} \m{\Sigma} \m{U}_{\c{J}}) \m{K}_{d_{\c{V}}, p-d_{\c{U}}} = \m{0}.
\end{align*}
The remaining term evaluates explicitly to:
\begin{align*}
    (\m{U}_{\c{J}}^{\top} \otimes \m{U}_{\c{I}^{c}}^{\top}) (\m{\Sigma} \otimes \m{\Sigma}) (\m{U}_{\c{J}} \otimes \m{U}_{\c{I}^{c}}) = (\m{U}_{\c{J}}^{\top} \m{\Sigma} \m{U}_{\c{J}}) \otimes (\m{U}_{\c{I}^{c}}^{\top} \m{\Sigma} \m{U}_{\c{I}^{c}}) = \m{\Lambda}_{\c{J}} \otimes \m{\Lambda}_{\c{I}^{c}}.
\end{align*}
Consequently, the covariance matrix of $\operatorname{vec}(\tilde{\m{Z}})$ reduces strictly to the diagonal matrix $(1+\kappa) (\m{\Lambda}_{\c{J}} \otimes \m{\Lambda}_{\c{I}^{c}})$. This implies that the entries $\tilde{\m{Z}}_{ij}$ are mutually independent centered Gaussians with variances $\operatorname{Var}(\tilde{\m{Z}}_{ij}) = (1+\kappa) \lambda_{i} \lambda_{j}$ for all $i \in \c{I}^{c}$ and $j \in \c{J}$. Standardizing $\tilde{\m{Z}}_{ij}$ gives:
\begin{align*}
    \left( \frac{\tilde{\m{Z}}_{ij}}{\lambda_{j} - \lambda_{i}} \right)^{2} = (1+\kappa) \frac{\lambda_{j} \lambda_{i}}{(\lambda_{j} - \lambda_{i})^{2}} \chi^{2}_{1, ij}, \quad \chi^{2}_{1, ij} \stackrel{\text{i.i.d.}}{\sim} \chi^{2}_{1}.
\end{align*}
\end{proof}

\section{Proofs for Section~\ref{sec:high:dim}}
\label{app:high:dim}

\begin{proof}[Proof of \cref{thm:macro:bias}]
If $d^{*} = 0$ or $k = 0$, the identities hold exactly. Otherwise, let $\m{C}_{k} \defeq \m{U}_{[d^{*}]}^{\top} \hat{\m{V}}_{[k]}$. Summing \eqref{eq:hd:population:overlap} over the distinct spike blocks gives
\begin{align*}
    \m{C}_{k}^{\top} \m{C}_{k} = \hat{\m{V}}_{[k]}^{\top} \left( \sum_{\c{B}} \m{\Pi}_{\c{B}} \right) \hat{\m{V}}_{[k]} \xrightarrow{a.s.} \diag \big( g (\lambda_{j} ; \gamma) \big)_{j = 1}^{k},
\end{align*}
where $g (\lambda_{j} ; \gamma) = 0$ for $j > d^{*}$. Continuity of the eigenvalues of $\m{C}_{k}^{\top} \m{C}_{k}$ proves both limits. Finally, convergence over any fixed finite collection of dimensions is uniform.

\end{proof}

\begin{proof}[Proof of \cref{thm:macro:bias:split}]
The case $k = 0$ is immediate. For $k \ge 1$, define the population signal projector by $\m{\Pi}_{[d^{*}]} \defeq \m{U}_{[d^{*}]} \m{U}_{[d^{*}]}^{\top}$ and its orthogonal complement by $\m{\Pi}_{[d^{*}]}^{\perp} \defeq \m{I}_{p} - \m{\Pi}_{[d^{*}]}$.
For $1 \le i, j \le D$, define $\m{e}_{i} \defeq \m{\Pi}_{[d^{*}]}^{\perp} \hat{\m{v}}_{i}^{(1)}$ and $\m{f}_{j} \defeq \m{\Pi}_{[d^{*}]}^{\perp} \hat{\m{v}}_{j}^{(2)}$.
We first claim that their inner products vanish almost surely:
\begin{align}\label{eq:hd:residual:orthogonality}
    \m{e}_{i}^{\top} \m{f}_{j} \equiv \hat{\m{v}}^{(1)\top}_{i} \m{\Pi}_{[d^{*}]}^{\perp} \hat{\m{v}}_{j}^{(2)} \xrightarrow{a.s.} 0.
\end{align}
The isotropic eigenvector delocalization estimates (see \cite[Theorems~2.11, and~2.17]{bloemendal2016principal}) imply that, for every fixed $\varepsilon \in (0, 1/2)$ and $A > 1$, and all sufficiently large $n$,
\begin{align*}
    \b{P} \left( | \m{e}_{i}^{\top} \m{f}_{j} | > p^{-1/2 + \varepsilon} \right) = \b{E} \left[ \b{P} \left( | \m{e}_{i}^{\top} \m{f}_{j} | > p^{-1/2 + \varepsilon} \mid \m{e}_{i} \right) \right] \le \sup_{\substack{\m{w} \in \c{U}_{[d^{*}]}^{\perp} \\ \| \m{w} \|_{2} = 1}} \b{P} \left( | \m{w}^{\top} \hat{\m{v}}_{j}^{(2)} | > p^{-1/2 + \varepsilon} \right) \le p^{-A}.
\end{align*}
Consequently, for every fixed $\tilde{\varepsilon} \in (0, 1)$, we get $\sum_{n
= 1}^{\infty} \b{P} ( | \m{e}_{i}^{\top} \m{f}_{j} | > \tilde{\varepsilon} )<
\infty$. The first Borel--Cantelli lemma implies that $| \m{e}_{i}^{\top}
\m{f}_{j} | > \tilde{\varepsilon}$ occurs only finitely often, almost surely.
Applying this conclusion to $\tilde{\varepsilon} = 1/r$, $r = 2, 3, \ldots$,
gives $\m{e}_{i}^{\top} \m{f}_{j} \xrightarrow{a.s.} 0$. Since $D$ is fixed,
this convergence holds simultaneously for all $1 \le i, j \le D$.

Now define
\begin{align*}
    \m{A}_{1} \defeq \m{U}_{[d^{*}]}^{\top} \hat{\m{V}}_{[D]}^{(1)}, \qquad \m{A}_{2} \defeq \m{U}_{[d^{*}]}^{\top} \hat{\m{V}}_{[k]}^{(2)}, \qquad \m{C}_{k} \defeq (\hat{\m{V}}_{[D]}^{(1)})^{\top} \hat{\m{V}}_{[k]}^{(2)}.
\end{align*}
Since
\begin{align*}
    [\m{C}_{k}]_{ij} = (\hat{\m{v}}_{i}^{(1)})^{\top} \hat{\m{v}}_{j}^{(2)}
    =
    \sum_{\ell = 1}^{d^{*}}
    \big( (\hat{\m{v}}_{i}^{(1)})^{\top} \m{u}_{\ell} \big)
    \big( \m{u}_{\ell}^{\top} \hat{\m{v}}_{j}^{(2)} \big)
    + \m{e}_{i}^{\top} \m{f}_{j} = [\m{A}_{1}^{\top} \m{A}_{2}]_{ij} + \m{e}_{i}^{\top} \m{f}_{j},
\end{align*}
we obtain from \eqref{eq:hd:residual:orthogonality} that $\vertii{\m{C}_{k} - \m{A}_{1}^{\top} \m{A}_{2}}_{2} \xrightarrow{a.s.} 0$.
Also, as in the proof of \cref{thm:macro:bias}, we get
\begin{align*}
    \m{A}_{1} \m{A}_{1}^{\top} \xrightarrow{a.s.} \m{G}_{1} \defeq \diag \big( g (\lambda_{\ell} ; \gamma_{1}) \big)_{\ell = 1}^{d^{*}}.
\end{align*}
All matrices $\m{A}_{1}$ and $\m{A}_{2}$ have operator norm at most one. Hence, we have $\vertii{\m{C}_{k}^{\top} \m{C}_{k} - \m{A}_{2}^{\top} \m{G}_{1} \m{A}_{2}}_{2} \xrightarrow{a.s.} 0$.
The diagonal entries of $\m{G}_{1}$ are constant within each repeated spike block $\c{B} \subseteq [d^{*}]$. Applying \eqref{eq:hd:population:overlap} to the second sample therefore yields
\begin{align*}
    \m{A}_{2}^{\top} \m{G}_{1} \m{A}_{2} &= \sum_{\c{B}} g (\lambda_{\c{B}} ; \gamma_{1}) (\hat{\m{V}}_{[k]}^{(2)})^{\top} \m{\Pi}_{\c{B}} \hat{\m{V}}_{[k]}^{(2)} \xrightarrow{a.s.} \diag \big( h (\lambda_{j} ; \gamma, \alpha) \big)_{j = 1}^{k},
\end{align*}
which results in $\m{C}_{k}^{\top} \m{C}_{k} \xrightarrow{a.s.} \diag \big( h (\lambda_{j} ; \gamma, \alpha) \big)_{j = 1}^{k}$.
When $d^{*} = 0$, the signal terms are absent and the same conclusion follows directly from \eqref{eq:hd:residual:orthogonality}.
Finally, the two squared distances are $D + k - 2 \vertii{\m{C}_{k}}_{2}^{2}$ and $D + k - 2 \vertii{\m{C}_{k}}_{1}$, respectively. Continuity of the eigenvalues of $\m{C}_{k}^{\top} \m{C}_{k}$ proves the asserted limits.
\end{proof}

\begin{proof}[Proof of \cref{prop:half:ideal}]
If either factor defining $h (\lambda ; \gamma, \alpha)$ is zero, the inequality is immediate. Otherwise, letting $x \defeq \lambda - 1$, both $\alpha$ and $1 - \alpha$ exceed $\gamma / x^{2}$. On the interval $t > \gamma / x^{2}$,
\begin{align*}
    \tilde{g} (t) \defeq \log g (\lambda ; \gamma/t) = \log (t - \gamma / x^{2}) - \log (t + \gamma / x), \quad
    \tilde{g}'' (t) = - \frac{1}{(t - \gamma / x^{2})^{2}} + \frac{1}{(t + \gamma / x)^{2}} < 0.
\end{align*}
Strict concavity gives $\tilde{g} (\alpha) + \tilde{g} (1 - \alpha) \le 2 \tilde{g} (1/2)$, with equality only at $\alpha = 1/2$. Exponentiating proves the claim. If $\lambda > 1 + \sqrt{2 \gamma}$, the right-hand side is positive, so equality holds only at $\alpha = 1/2$.
\end{proof}

\begin{proof}[Proof of \cref{cor:hd:selected:dimensions}]
Under the prescribed calibration, the limiting increments are
\begin{align*}
    L_{\Theta, \eta_{\Theta} (c_{\gamma} ; \gamma)} (k)
    - L_{\Theta, \eta_{\Theta} (c_{\gamma} ; \gamma)} (k - 1)
    &= 2 \{ g (c_{\gamma} ; 2 \gamma)^{2} - g (\lambda_{k} ; 2 \gamma)^{2} \}, \\
    L_{\c{W}, \eta_{\c{W}} (c_{\gamma} ; \gamma)} (k)
    - L_{\c{W}, \eta_{\c{W}} (c_{\gamma} ; \gamma)} (k - 1)
    &= 2 \{ g (c_{\gamma} ; 2 \gamma) - g (\lambda_{k} ; 2 \gamma) \}.
\end{align*}
Because $c_{\gamma} > 1 + \sqrt{2 \gamma}$ and $g (\lambda ; 2 \gamma)$ is strictly increasing above this threshold, both increments are negative exactly when $\lambda_{k} > c_{\gamma}$ and positive exactly when $\lambda_{k} < c_{\gamma}$. The absence of spikes at the cutoff therefore gives the stated unique limiting minimizer.
By Theorem~\ref{thm:macro:bias:split} and the assumed convergence of $\eta_{n}$,
\begin{align*}
    \max_{0 \le k \le D}
    \left|
        \widehat{L}_{\mathsf{G}, \eta_{n}} (k)
        - L_{\mathsf{G}, \eta_{\mathsf{G}} (c_{\gamma} ; \gamma)} (k)
    \right|
    \xrightarrow{a.s.} 0.
\end{align*}
On the finite candidate set, the unique limiting minimum is separated from all other values by a positive gap. Consequently, the empirical minimizer eventually equals the limiting minimizer almost surely.
\end{proof}

\begin{proof}[Proof of \cref{thm:hd:oracle:calibration}]
The map $\gamma \mapsto \eta_{\Theta}^{\mathrm{or}} (\gamma)$ in \eqref{eq:hd:optimal:penalty} is continuous on $(0, \infty)$. Since $\gamma_{n} \to \gamma$, Theorem~\ref{thm:macro:bias:split} implies that the objective in \eqref{eq:hd:calibrated:selector} converges almost surely, uniformly over $0 \le k \le D$, to
\begin{align*}
    C (k) \defeq D + 2 k g (c_{\Theta} (\gamma) ; 2 \gamma)^{2} - 2 \sum_{j = 1}^{k} g (\lambda_{j} ; 2 \gamma)^{2}.
\end{align*}

To compare $C$ with the full-sample oracle loss, consider the cost of adding direction $j$. For $1 \le j \le D$, the increment is $C (j) - C (j - 1) = 2 \left[ g (c_{\Theta} (\gamma) ; 2 \gamma)^{2} - g (\lambda_{j} ; 2 \gamma)^{2} \right]$.
Because $c_{\Theta} (\gamma) > 1 + \sqrt{2 \gamma}$, the increment $C (j) - C (j - 1) = 2[ g (c_{\Theta} (\gamma) ; 2 \gamma)^{2} - g (\lambda_{j} ; 2 \gamma)^{2} ]$ has the same sign as $c_{\Theta} (\gamma) - \lambda_{j}$, with equality exactly at the cutoff. The same sign rule holds for the oracle increment $1 - 2 g (\lambda_{j} ; \gamma)$ in \eqref{eq:oracle:incre}, since $1 - 2 g (\lambda_{j} ; \gamma) > 0$ if and only if $c_{\Theta} (\gamma) - \lambda_{j} > 0$. Both sequences of increments are nondecreasing, so
\begin{align*}
    \argmin_{0 \le k \le D} C (k) = \argmin_{0 \le k \le D} L_{\Theta}^{\mathrm{or}} (k).
\end{align*}
Since the candidate set is finite, every value outside this common minimizer set is separated from the minimum of $C$ by a positive gap. Uniform convergence consequently implies that $\hat{d}_{\Theta}^{\mathrm{cal}}$ eventually belongs to the common minimizer set almost surely.
Finally, in order to establish vanishing regret, define
\begin{align*}
    \varepsilon_{n} \defeq \max_{0 \le k \le D} \left| \widehat{L}_{\Theta}^{\mathrm{or}} (k) - L_{\Theta}^{\mathrm{or}} (k) \right|.
\end{align*}
Proposition~\ref{thm:macro:bias} gives $\varepsilon_{n} \xrightarrow{a.s.} 0$. Once $\hat{d}_{\Theta}^{\mathrm{cal}}$ belongs to the common minimizer set,
\begin{align*}
    0 \le \widehat{L}_{\Theta}^{\mathrm{or}} (\hat{d}_{\Theta}^{\mathrm{cal}}) - \min_{0 \le k \le D} \widehat{L}_{\Theta}^{\mathrm{or}} (k) \le 2 \varepsilon_{n}.
\end{align*}
This proves the asserted almost-sure convergence of the regret, including cases with spikes at the cutoff. Finally, if no population spike equals $c_{\Theta} (\gamma)$, the common minimizer is uniquely $d_{\Theta}^{\mathrm{or}}$. Hence $\hat{d}_{\Theta}^{\mathrm{cal}} = d_{\Theta}^{\mathrm{or}}$ eventually almost surely, proving the dimension convergence.
\end{proof}

\begin{proof}[Proof of \cref{thm:hd:dimension:consistency}]
The increments in \eqref{eq:risk:curve:slope:chordal} and \eqref{eq:risk:curve:slope:wasserstein} are nondecreasing. Under the relevant no-equality condition, they are strictly negative up to $d_{\mathsf{G}, \eta}$ and strictly positive thereafter. Thus $d_{\mathsf{G}, \eta}$ is the unique minimizer of $L_{\mathsf{G}, \eta}$.
Moreover,
\begin{align*}
    \max_{0 \le k \le D} | \widehat{L}_{\mathsf{G}, \eta_{n}} (k) - L_{\mathsf{G}, \eta} (k) | \le \max_{0 \le k \le D} | \widehat{L}_{\mathsf{G}} (k) - L_{\mathsf{G}} (k) | + D | \eta_{n} - \eta | \xrightarrow{a.s.} 0.
\end{align*}
On the finite candidate set, the unique limiting minimum is separated from every other value by a positive gap. Once the uniform error is smaller than half this gap, the empirical minimizer equals $d_{\mathsf{G}, \eta}$.
\end{proof}

\begin{proof}[Proof of \cref{cor:hd:rank:consistency}]
The stated condition ensures that every supercritical population spike is included in $d_{\mathsf{G}, \eta}$. All remaining spikes have zero limiting half-sample alignment and are excluded because $\eta < 1$. Hence $d_{\mathsf{G}, \eta} = d_{\gamma}^{*}$, and the conclusion follows from \cref{thm:hd:dimension:consistency}.
\end{proof}

\begin{proof}[Proof of \cref{prop:hd:null:calibration}]
Under the Gaussian null, $\m{X}_{2} \m{O} \overset{d}{=} \m{X}_{2}$ for every deterministic $\m{O} \in \c{O} (p)$, and hence $\m{O}^{\top} \hat{\m{\Sigma}}_{2} \m{O} \overset{d}{=} \hat{\m{\Sigma}}_{2}$. Since the selected eigenvalues are simple almost surely, their ordered eigenvector lines have a jointly rotation-invariant distribution.
To remove the sign ambiguity, let $s_{1}, \ldots, s_{D}$ be independent random signs, independent of both samples, with $\b{P} (s_{j} = 1) = \b{P} (s_{j} = -1) = 1/2$, and define
\begin{align*}
    \tilde{\m{V}}_{[D]}^{(2)} \defeq \hat{\m{V}}_{[D]}^{(2)} \diag (s_{1}, \ldots, s_{D}).
\end{align*}
The signed frame $\tilde{\m{V}}_{[D]}^{(2)}$ has the uniform distribution on $\c{O} (p, D)$: its law is invariant under every orthogonal transformation of $\b{R}^{p}$. In particular, each column $\tilde{\m{v}}_{j}^{(2)}$ is uniformly distributed on $\mathbb{S}^{p - 1}$ and remains independent of the first sample.

Conditionally on $\hat{\m{\Pi}}_{[D]}^{(1)}$, choose $\m{Q} \in \c{O} (p)$ such that $\hat{\m{\Pi}}_{[D]}^{(1)} = \m{Q} \diag (\m{I}_{D}, \m{0}_{p - D}) \m{Q}^{\top}$. Since signs leave $q_{j}$ unchanged and $\m{Q}^{\top} \tilde{\m{v}}_{j}^{(2)} \sim \operatorname{Unif}(\mathbb{S}^{p - 1})$, the normalized Gaussian representation gives
\begin{align*}
    q_{j} = \sum_{\ell = 1}^{D} \big[ \m{Q}^{\top} \tilde{\m{v}}_{j}^{(2)} \big]_{\ell}^{2} \overset{d}{=} \frac{\sum_{\ell = 1}^{D} G_{\ell}^{2}}{\sum_{\ell = 1}^{p} G_{\ell}^{2}}, \qquad G_{1}, \ldots, G_{p} \stackrel{\mathrm{iid}}{\sim} \c{N} (0, 1),
\end{align*}
and this establishes the asserted beta law.
Finally, if every $q_{j} \le (1 - \eta)/2$, all increments in \eqref{eq:hd:cumulative:criterion} are nonnegative, and the smallest minimizer is zero. Thus
\begin{align*}
    \{ \hat{d}_{\Theta, \eta} > 0 \} \subseteq \bigcup_{j = 1}^{D} \left\{ q_{j} > \frac{1 - \eta}{2} \right\}.
\end{align*}
The union bound and the beta marginals prove \eqref{eq:hd:null:bound}.
\end{proof}

\begin{proof}[Proof of \cref{thm:hd:null:joint}]
Use the signed frame $\tilde{\m{V}}_{[D]}^{(2)}$ from the proof of \cref{prop:hd:null:calibration}, and define the cross-matrix $\tilde{\m{C}}_{D} \defeq (\hat{\m{V}}_{[D]}^{(1)})^{\top} \tilde{\m{V}}_{[D]}^{(2)}$.
Conditionally on the first-sample basis, rotation invariance implies that $\tilde{\m{C}}_{D}$ has the distribution of the upper-left $D \times D$ block of a Haar-distributed orthogonal matrix in $\c{O} (p)$. The Gaussian approximation \cite[Theorem~1]{jiang2006many} therefore gives the convergence in total variation, thus implies the convergence in law:
\begin{align*}
    \sqrt{p} \, \tilde{\m{C}}_{D} \xrightarrow{d} \m{G}, \qquad G_{ij} \stackrel{\mathrm{iid}}{\sim} \c{N} (0, 1), \qquad 1 \le i, j \le D.
\end{align*}
Since $p q_{j} = \sum_{i = 1}^{D} (\sqrt{p} \, \tilde{C}_{ij})^{2} \xrightarrow{d} Y_{j} \overset{d}{=} \chi_{D}^{2}$, continuous mapping yields $(p q_{1}, \ldots, p q_{D}) \xrightarrow{d} (Y_{1}, \ldots, Y_{D})$. The limiting variables are independent because they are squared norms of distinct columns of $\m{G}$. Applying continuous mapping again proves \eqref{eq:hd:null:cumulative:limit}.
Finally, the distribution of $M_{D}$ is continuous, hence $H_{D} (H_{D}^{-1} (1 - \varepsilon)) = 1 - \varepsilon$. By \eqref{eq:hd:null:selection:event}, the result follows.
\end{proof}

\end{appendix}
\end{document}